\documentclass[12pt]{article} 
\usepackage[utf8]{inputenc}
\usepackage{soul} 
\usepackage{setspace} 
\usepackage{shellesc}
\usepackage{mathrsfs,amsthm,xcolor,verbatim,bbm,amsmath,amsfonts,amssymb,nicefrac,enumitem,bm,mathtools,xparse,etoolbox,textcomp}
\usepackage[margin=1in,a4paper]{geometry}
\usepackage{xurl}
\usepackage{frontespizio}
\usepackage{caption}
\usepackage{booktabs} 
\usepackage{makecell}          
\usepackage{array}             
\usepackage{graphicx}          
\usepackage{float}             
\usepackage{colortbl}
\usepackage[table]{xcolor}

\usepackage{tikz}
\usetikzlibrary{matrix, arrows.meta}
\usepackage{tikz-3dplot}
\usepackage{pst-solides3d}
\usetikzlibrary{automata,chains}
\usepackage{subcaption}
\usepackage{footnote}
\usepackage[most]{tcolorbox}
\usepackage{cite}
\usepackage[hypertexnames=false]{hyperref}
\hypersetup{
	breaklinks=true,
    colorlinks,
    linkcolor={blue!90!black},  
    citecolor={green},
    urlcolor={blue!99!black}   
}
\usepackage[capitalise,nameinlink,sort]{cleveref}

\newcommand\summaryname{Abstract}
    {\cleardoublepage\null\vfill \begin{center}%
    \bfseries{\summaryname} \end{center}}
    {\vfill\null}
    {\cleardoublepage\null\vfill\begin{center}%
    \bfseries Ringraziamenti\end{center}}%
    {\vfill\null}

\crefname{enumi}{item}{items}

\crefname{equation}{}{}
\crefname{subsection}{Subsection}{Subsections}

\usepackage{thmtools}

\theoremstyle{plain}
\newtheorem{theorem}{Theorem}[section]

\newtheorem{prop}[theorem]{Proposition}

\newtheorem{definition}[theorem]{Definition}
\newtheorem{setting}[theorem]{Setting}

\theoremstyle{definition}

\newcommand{\lrexpo}{\nu}

\DeclareMathAlphabet{\mathpzc}{OT1}{pzc}{m}{it}

\numberwithin{equation}{section}

\DeclareFontEncoding{LS1}{}{}
\DeclareFontSubstitution{LS1}{stix}{m}{n}
\DeclareMathAlphabet{\mathscr}{LS1}{stixscr}{m}{n}

\newcommand{\R}{\mathbb{R}}
\newcommand{\N}{\mathbb{N}}

\newcommand{\bbL}{\mathbb{L}}
\newcommand{\bbI}{\mathbb{I}}
\newcommand{\bbM}{\mathbb{M}}
\newcommand{\bbA}{\mathbb{A}}

\renewcommand{\c}[1]{\mathfrak{c}^{#1}}

\DeclarePairedDelimiterX{\infdivx}[2]{(}{)}{%
  #1\;\delimsize\|\;#2%
}

\newcommand{\m}{\mathbf{m}}

\newcommand{\smallsum}{\textstyle\sum}

\newcommand{\with}{\curvearrowleft}

\newcommand{\cG}{\mathcal{G}}

\newcommand{\cL}{\mathcal{L}}
\newcommand{\cM}{\mathcal{M}}
\newcommand{\cN}{\mathcal{N}}

\newcommand{\act}{\mathscr{a}}

\newcommand{\scrC}{\mathscr{C}}

\newcommand{\scrL}{\mathscr{L}}

\newcommand{\scrN}{\mathscr{N}}

\newcommand{\scrc}{\mathscr{c}}

\newcommand{\fC}{\mathfrak{C}}

\newcommand{\fM}{\mathfrak{M}}

\newcommand{\fc}{\mathfrak{c}}
\newcommand{\fd}{\mathfrak{d}}

\newcommand{\fm}{\mathfrak{m}}

\newcommand{\fx}{\mathfrak{x}}
\newcommand{\fy}{\mathfrak{y}}

\renewcommand{\emptyset}{\varnothing}

\definecolor{darkgreen}{rgb}{0.0, 0.6, 0.0}
\definecolor{cyan(process)}{rgb}{0.0, 0.72, 0.92}
\definecolor{cinnamon}{rgb}{0.82, 0.41, 0.12}
\definecolor{darkred}{rgb}{0.55, 0.0, 0.0} 
\definecolor{persianblue}{rgb}{0.11, 0.22, 0.73}
\definecolor{resnamecolor}{rgb}{0.0, 0.5, 0.0}

\newcommand{\resname}[1]{\textcolor{red}{#1}}

\DeclarePairedDelimiter{\norm}{\lVert}{\rVert}
\DeclarePairedDelimiter{\abs}{\lvert}{\rvert}

\DeclarePairedDelimiter{\spro}{\langle}{\rangle}

\newcommand{\Norm}[1]{\bigl\Vert #1 \bigr\Vert}

\newcommand{\NormM}[1]{{\left\vert\kern-0.25ex\left\vert\kern-0.25ex\left\vert #1
    \right\vert\kern-0.25ex\right\vert\kern-0.25ex\right\vert}}

\newcommand{\Abs}[1]{\bigl| #1 \bigr|}

\newcommand{\Forall}{\forall \,}

\renewcommand{\d}{\, \mathrm{d}}
\newcommand{\eps}{\varepsilon}

\newcommand{\qandq}{\quad\text{and}\quad}
\newcommand{\andq}{\text{and}\quad}

\newcommand{\qqandqq}{\qquad\text{and}\qquad}
\newcommand{\andqq}{\text{and}\qquad}

\NewDocumentCommand{\Index}{ m o m o }{^{#1\IfValueT{#2}{, #2}}_{#3\IfValueT{#4}{, #4}}}

\NewDocumentEnvironment{cproof}{m}
{\begin{proof}[Proof of \cref{#1}]}%
{\noindent The proof of~\cref{#1} is thus complete.
\end{proof}}

\NewDocumentEnvironment{cproof2}{m}
{\begin{proof}[Proof of \cref{#1}]}%
{\noindent This completes the proof of~\cref{#1}.
\end{proof}}

\ExplSyntaxOn

\seq_new:N \g_cflist_loaded
\seq_new:N \g_cflist_pending

\NewDocumentCommand{\cfadd} { m } {
  \seq_if_in:NnF \g_cflist_loaded { #1 } {
    \seq_if_in:NnF \g_cflist_pending { #1 } {
      \seq_gput_right:Nn \g_cflist_pending { #1 }
    }
  }
}

\NewDocumentCommand{\cfconsiderloaded} { m } {
  \seq_gput_right:Nn \g_cflist_loaded {#1}
}

\NewDocumentCommand{\cfremove} { m } {
  \seq_gremove_all:Nn \g_cflist_pending { #1 }
}

\NewDocumentCommand{\cfload} { o } {
  \seq_if_empty:NTF \g_cflist_pending {
    \IfValueTF{#1}{\ignorespaces}{\unskip}
  } {
    (cf.\ \cref{\seq_use:Nn \g_cflist_pending {,}})\IfValueTF{#1}{#1~}{\unskip}
    \seq_gconcat:NNN \g_cflist_loaded \g_cflist_loaded \g_cflist_pending
    \seq_gclear:N \g_cflist_pending
    \IfValueT{#1}{\ignorespaces}
  }
}

\NewDocumentCommand{\cfclear} {} {
  \seq_gclear:N \g_cflist_loaded
  \seq_gclear:N \g_cflist_pending
}

\NewDocumentCommand{\cfout} { o } {
  \seq_if_empty:NTF \g_cflist_pending {\unskip\IfValueT{#1}{\ignorespaces}} {
    (cf.\ \cref{\seq_use:Nn \g_cflist_pending {,}})\IfValueTF{#1}{#1~}{\unskip}
    \seq_gclear:N \g_cflist_pending
    \IfValueT{#1}{\ignorespaces}
  }
}

\NewDocumentCommand{\ifnocf} { m } {
  \seq_if_empty:NT \g_cflist_pending { #1 }
}

\ExplSyntaxOff

\ExplSyntaxOn

\int_new:N \g_furthermore

\NewDocumentCommand{\Moreover}{ o o }{
  \IfValueT{#1}{
    \str_case:nn {#1} {
      {Furthermore} {\int_set:Nn {\g_furthermore} {0}}
      {Moreover} {\int_set:Nn {\g_furthermore} {1}}
      {In~addition} {\int_set:Nn {\g_furthermore} {2}}
      {note} {\bool_gset_true:N \g_noteobserve}
      {observe} {\bool_gset_false:N \g_noteobserve}
    }
    \IfValueT{#2}{
      \str_case:nn {#2} {
        {Furthermore} {\int_set:Nn {\g_furthermore} {0}}
        {Moreover} {\int_set:Nn {\g_furthermore} {1}}
        {In~addition} {\int_set:Nn {\g_furthermore} {2}}
        {note} {\bool_gset_true:N \g_noteobserve}
        {observe} {\bool_gset_false:N \g_noteobserve}
      }
    }
  }
  \int_case:nn { \int_mod:nn {\g_furthermore} {3} } {
    { 0 } { Furthermore,~\nobs that}
    { 1 } { Moreover,~\nobs that}
    { 2 } { In~addition,~\nobs that}
  }
  \int_incr:N \g_furthermore
  \IfValueF{#1}{~}
}

\ExplSyntaxOff

\ExplSyntaxOn

\bool_new:N \g_noteobserve

\NewDocumentCommand{\setnote}{}{
  \bool_gset_true:N \g_noteobserve
}

\NewDocumentCommand{\setobserve}{}{
  \bool_gset_false:N \g_noteobserve
}

\NewDocumentCommand{\nobs}{ o }{
  \IfValueT{#1}{
    \str_if_eq:noTF {note} {#1} {
      \bool_gset_true:N \g_noteobserve
    } {
      \str_if_eq:noTF {Note} {#1} {
        \bool_gset_true:N \g_noteobserve
      } {
        \bool_gset_false:N \g_noteobserve
      }
    }
  }
  \bool_if:nTF { \g_noteobserve } {
    \bool_gset_false:N \g_noteobserve
    note
  } {
    \bool_gset_true:N \g_noteobserve
    observe
  }
  \IfValueF{#1}{~}
}

\NewDocumentCommand{\Nobs}{ o }{
  \IfValueT{#1}{
    \str_if_eq:noTF {note} {#1} {
      \bool_gset_true:N \g_noteobserve
    } {
      \str_if_eq:noTF {Note} {#1} {
        \bool_gset_true:N \g_noteobserve
      } {
        \bool_gset_false:N \g_noteobserve
      }
    }
  }
  \bool_if:nTF { \g_noteobserve } {
    \bool_gset_false:N \g_noteobserve
    Note
  } {
    \bool_gset_true:N \g_noteobserve
    Observe
  }
  \IfValueF{#1}{~}
}

\ExplSyntaxOff

\ExplSyntaxOn

\bool_new:N \g_hencetherefore

\NewDocumentCommand{\hence}{ o }{
  \IfValueT{#1}{
    \str_if_eq:noTF {hence} {#1} {
      \bool_gset_true:N \g_hencetherefore
    } {
      \str_if_eq:noTF {Hence} {#1} {
        \bool_gset_true:N \g_hencetherefore
      } {
        \bool_gset_false:N \g_hencetherefore
      }
    }
  }
  \bool_if:nTF { \g_hencetherefore } {
    \bool_gset_false:N \g_hencetherefore
    hence
  } {
    \bool_gset_true:N \g_hencetherefore
    therefore
  }
  \IfValueF{#1}{~}
}

\NewDocumentCommand{\Hence}{ o }{
  \IfValueT{#1}{
    \str_if_eq:noTF {hence} {#1} {
      \bool_gset_true:N \g_hencetherefore
    } {
      \str_if_eq:noTF {Hence} {#1} {
        \bool_gset_true:N \g_hencetherefore
      } {
        \bool_gset_false:N \g_hencetherefore
      }
    }
  }
  \bool_if:nTF { \g_hencetherefore } {
    \bool_gset_false:N \g_hencetherefore
    Hence,~we~obtain
  } {
    \bool_gset_true:N \g_hencetherefore
    Therefore,~we~obtain
  }
  \IfValueF{#1}{~}
}

\ExplSyntaxOff

\ExplSyntaxOn

\seq_const_from_clist:Nn \g_prove_mru {
  establish,
  demonstrate,
  prove,
  show,
  imply,
  ensure
}

\prop_new:N \l__verbs
\prop_put:Nnn \l__verbs {show} {shows}
\prop_put:Nnn \l__verbs {imply} {implies}
\prop_put:Nnn \l__verbs {demonstrate} {demonstrates}
\prop_put:Nnn \l__verbs {prove} {proves}
\prop_put:Nnn \l__verbs {establish} {establishes}
\prop_put:Nnn \l__verbs {ensure} {ensures}
\prop_put:Nnn \l__verbs {assure} {assures}

\tl_new:N \g_wordtmp
\seq_new:N \l_mytmps

\cs_generate_variant:Nn \str_if_in:nnTF { nVTF }
\cs_generate_variant:Nn \str_if_in:nnTF { xVTF }

\NewDocumentCommand{\prove}{ o }{
  \IfValueTF{#1}{
    \seq_clear:N \l_mytmps
    \seq_map_inline:Nn \g_prove_mru {
      \str_if_eq:nnTF {##1} {ensure} {
        \str_set:Nn \l_temps {n}
      } {
        \str_set:Nx \l_temps {\str_head_ignore_spaces:n {##1}}
      }
      \str_if_in:xVTF {#1} \l_temps {
        \seq_put_right:Nn \l_mytmps {##1}
      } { }
    }
    \seq_get_right:NN \l_mytmps \g_wordtmp
  } {
    \seq_get_right:NN \g_prove_mru \g_wordtmp
  }
  \tl_use:N \g_wordtmp
  \IfValueTF{#1}{}{~}
  \seq_gput_left:NV \g_prove_mru \g_wordtmp
  \seq_gremove_duplicates:N \g_prove_mru
}

\NewDocumentCommand{\proves}{ o }{
  \IfValueTF{#1}{
    \seq_clear:N \l_mytmps
    \seq_map_inline:Nn \g_prove_mru {
      \str_if_eq:nnTF {##1} {ensure} {
        \str_set:Nn \l_temps {n}
      } {
        \str_set:Nx \l_temps {\str_head_ignore_spaces:n {##1}}
      }
      \str_if_in:xVTF {#1} \l_temps {
        \seq_put_right:Nn \l_mytmps {##1}
      } { }
    }
    \seq_get_right:NN \l_mytmps \g_wordtmp
  } {
    \seq_get_right:NN \g_prove_mru \g_wordtmp
  }
  \str_set:NV \l_tmpa_str \g_wordtmp
  \prop_get:NVN \l__verbs \l_tmpa_str \l_tmpa_tl
  \tl_use:N \l_tmpa_tl
  \IfValueTF{#1}{}{~}
  \seq_gput_left:NV \g_prove_mru \g_wordtmp
  \seq_gremove_duplicates:N \g_prove_mru
}

\newcommand{\llabel}[1]{\savelabel{#1}\label{\loc.#1}\ignorespaces}

\clist_new:N \l_localreflist
\clist_new:N \l_reflist

\NewDocumentCommand{\lref} { m } {
  \clist_set:No \l_localreflist {#1}
  \clist_clear:N \l_reflist
  \clist_map_inline:Nn \l_localreflist { \clist_put_right:Nn \l_reflist {\loc.##1} }
  \cref{\l_reflist}
}

\NewDocumentCommand{\Lref} { m } {
  \clist_set:No \l_localreflist {#1}
  \clist_clear:N \l_reflist
  \clist_map_inline:Nn \l_localreflist { \clist_put_right:Nn \l_reflist {\loc.##1} }
  \Cref{\l_reflist}
}

\NewDocumentCommand{\itref}{ m m }{
  \clist_set:No \l_localreflist {#2}
  \clist_clear:N \l_reflist
  \clist_map_inline:Nn \l_localreflist { \clist_put_right:Nn \l_reflist {#1.##1} }
  \cref{\l_reflist}~in~\cref{#1}
}

\seq_new:N \l_enum_seq
\int_new:N \l_num_items

\bool_new:N \g_commaused_bool

\providecommand{\comma}{}

\cs_new:Nn \enum_it:nn {
  \int_case:nnF {\l_num_items - #1} {
    {0} {
      \renewcommand{\comma}{}
      #2
      \space
    }
    {1} {
      \bool_gset_false:N \g_commaused_bool
      \renewcommand{\comma}{,~\bool_gset_true:N \g_commaused_bool}
      #2
      \bool_if:NTF \g_commaused_bool {} {,~}
      and~
    }
  } {
    \bool_gset_false:N \g_commaused_bool
    \renewcommand{\comma}{,~\bool_gset_true:N \g_commaused_bool}
    #2
    \bool_if:NTF \g_commaused_bool {} {,~}
  }
}

\cs_new:Nn \enum_it_U:nn {
  \int_case:nnF {\l_num_items - #1} {
    {0} {
      \renewcommand{\comma}{}
      #2
      \space
    }
    {1} {
      \bool_gset_false:N \g_commaused_bool
      \renewcommand{\comma}{,~\bool_gset_true:N \g_commaused_bool}
      #2
      \bool_if:NTF \g_commaused_bool {} {,~}
      and~
    }
  } {
    \bool_gset_false:N \g_commaused_bool
    \renewcommand{\comma}{,~\bool_gset_true:N \g_commaused_bool}
    \int_compare:nTF {#1=1} {
      \text_titlecase_first:n {#2}
    } {
      #2
    }
    \bool_if:NTF \g_commaused_bool {} {,~}
  }
}

\cs_generate_variant:Nn \tl_if_eq:nnTF {onTF}

\cs_new:Nn \enum:nnnn {
  \seq_set_split:Nnn \l_enum_seq ; {#1}
  \seq_remove_all:Nn \l_enum_seq { }
  \seq_remove_all:Nn \l_enum_seq {#2}
  \seq_log:N \l_enum_seq
  \int_set:Nn \l_num_items {\seq_count:N \l_enum_seq}
  \int_log:N \l_num_items
  \int_case:nnF {\l_num_items} {
    { 0 } { 0 }
    { 1 } {
      \IfBooleanTF{#4} {
        \tl_set:Nn \l_text_case_exclude_arg_tl {\cref}
        \bool_set_false:N \l_text_titlecase_check_letter_bool
        \text_titlecase_first:n {\seq_use:Nn \l_enum_seq {}}
      } {
        \seq_use:Nn \l_enum_seq {}
      }
      \space
      \tl_if_eq:onTF{#3}{-}{}{
        \bool_if:NTF \l_plural_bool {
          \prove[#3]~
        } {
          \proves[#3]~
        }
      }
      \ignorespaces
    }
    { 2 } {
      \IfBooleanTF{#4} {
        \tl_set:Nn \l_text_case_exclude_arg_tl {\cref}
        \bool_set_false:N \l_text_titlecase_check_letter_bool
        \text_titlecase_first:n {\seq_use:Nn \l_enum_seq {~and~}}
      } {
        \seq_use:Nn \l_enum_seq {~and~}
      }
      \space
      \tl_if_eq:onTF{#3}{-}{}{
        \prove[#3]~
      }
      \ignorespaces
    }
  } {
    \IfBooleanTF{#4} {
      \tl_set:Nn \l_text_case_exclude_arg_tl {\cref}
      \seq_indexed_map_function:NN \l_enum_seq \enum_it_U:nn
    } {
      \seq_indexed_map_function:NN \l_enum_seq \enum_it:nn
    }
    \tl_if_eq:onTF{#3}{-}{}{
      \prove[#3]~
    }
    \ignorespaces
  }
}

\cs_generate_variant:Nn \enum:nnnn {nxnn}
\cs_generate_variant:Nn \enum:nnnn {nxxn}

\NewDocumentCommand{\enum}{O{} m O{-} s}{
  \IfBooleanTF{#4}{
    \enum:nxnn {#2} {#1} {sindep} \BooleanFalse
  } {
    \enum:nxxn {#2} {#1} {#3} \BooleanFalse
  }
}

\NewDocumentCommand{\dott}{}{\ifnocf{.}\space}

\bool_new:N \g_arg_start_bool
\bool_gset_true:N \g_arg_start_bool

\NewDocumentCommand{\startnewargseq}{}{\bool_gset_true:N \g_arg_start_bool \tl_set:Nn \g_label_tl {}}

\cs_generate_variant:Nn \seq_if_in:NnTF {NxTF}
\cs_generate_variant:Nn \seq_remove_all:Nn {Nx}

\int_new:N \l_random_int

\cs_generate_variant:Nn \tl_if_head_eq_catcode:nNTF {oNTF}

\cs_generate_variant:Nn \tl_if_head_eq_catcode:nNTF {eNTF}
\cs_generate_variant:Nn \tl_log:n {o}
\cs_generate_variant:Nn \tl_log:n {f}
\cs_generate_variant:Nn \tl_log:n {x}
\cs_generate_variant:Nn \tl_log:n {e}

\cs_generate_variant:Nn \tl_if_in:nnTF {onTF}
\cs_generate_variant:Nn \tl_if_in:NnTF {NeTF}

\cs_generate_variant:Nn \tl_if_head_eq_meaning:nNTF {VNTF}

\seq_const_from_clist:Nn \g_arg_mru_this {
  Ahpr,
  Tapr,
  Ctapr,
  H
}

\seq_const_from_clist:Nn \g_arg_mru_nothis {
  Ia,
  Nwc,
  N,
  Itns,
  Fm,
  Itnswc,
  Mo
}

\seq_new:N \l_arg_seq
\tl_new:N \l_cons_tl
\tl_new:N \l_dummy_tl

\bool_new:N \g_debug_bool
\bool_gset_false:N \g_debug_bool

\bool_new:N \l_insidearg_bool

\bool_new:N \g_firstargletter_bool

\sys_gset_rand_seed:n {0903}

\bool_new:N \l_plural_bool
\tl_new:N \l_arg_verbs_tl

\definecolor{cyan(process)}{rgb}{0.0, 0.72, 0.92}
\definecolor{new_color}{rgb}{0.1, 0.6, 0.8}

\NewDocumentCommand{\argument}{mom}{
\color{black}
  \bool_set_false:N \l_plural_bool
  \tl_set:Nn \l_arg_verbs_tl {sindep}
  \keys_define:nn { benno/argument } {
    plural .value_forbidden:n = true,
    plural .code:n = {\bool_set_true:N \l_plural_bool},
    verbs .value_required:n = false,
    verbs .tl_set:N = \l_arg_verbs_tl,
  }
  \IfValueT{#2}{
    \keys_set:nn { benno/argument } {#2}
  }
  \bool_log:N \l_plural_bool
  \bool_gset_true:N \l_insidearg_bool
  \seq_set_split:Nnn \l_arg_seq ; {#1}
  \seq_remove_all:Nn \l_arg_seq { }
  \seq_log:N \l_arg_seq
  \tl_set:Nn \l_cons_tl {#3}
  \tl_trim_spaces:N \l_cons_tl
  \seq_if_in:NxTF \l_arg_seq {\lref{\g_label_tl}} {
    \seq_remove_all:Nx \l_arg_seq {\lref{\g_label_tl}}
    \seq_get_left:NNTF \l_arg_seq \l_dummy_tl {
      \tl_trim_spaces:N \l_dummy_tl
      \bool_gset_false:N \g_firstargletter_bool
      \tl_if_head_eq_catcode:VNTF \l_dummy_tl a {
        \bool_gset_true:N \g_firstargletter_bool
      } {
        \tl_if_head_eq_meaning:VNTF \l_dummy_tl {\cref} {
          \tl_set:Nx \l_tmpa_tl {\tl_tail:N \l_dummy_tl}
          \tl_set:Nx \l_tmpb_tl {\tl_head:N \l_tmpa_tl}
          \bool_gset_true:N \g_firstargletter_bool
          \tl_if_in:NeTF \l_tmpb_tl {lem\c_colon_str} {} {
            \tl_if_in:NeTF \l_tmpb_tl {thm\c_colon_str} {} {
              \tl_if_in:NeTF \l_tmpb_tl {prop\c_colon_str} {} {
                \tl_if_in:NeTF \l_tmpb_tl {cor\c_colon_str} {} {
                  \bool_gset_false:N \g_firstargletter_bool
                }
              }
            }
          }
        } {
        }
      }
      \bool_if:NTF \g_firstargletter_bool {
        \seq_set_eq:NN \l_tmpa_seq \g_arg_mru_this
        \seq_remove_all:Nn \l_tmpa_seq {H}
        \seq_get_right:NN \l_tmpa_seq \l_tmpa_tl
        \int_case:nnF {\seq_count:N \l_arg_seq} {
          {1} {
            \str_case:VnF {\l_tmpa_tl} {
              {Ahpr} {
                \bool_if:NT \g_debug_bool {C1.1}
                \seq_gput_left:Nn \g_arg_mru_this {Ahpr}
                \seq_gremove_duplicates:N \g_arg_mru_this
                \enum:nxnn {#1} {\lref{\g_label_tl}} {-} {\BooleanTrue}
                \hence~
                \bool_if:NTF \l_plural_bool {
                  \prove[\l_arg_verbs_tl]~\ignorespaces #3
                } {
                  \proves[\l_arg_verbs_tl]~\ignorespaces #3
                }
              }
              {Tapr} {
                \bool_if:NT \g_debug_bool {C1.2}
                \seq_gput_left:Nn \g_arg_mru_this {Tapr}
                \seq_gremove_duplicates:N \g_arg_mru_this
                \enum[\lref{\g_label_tl}]{
                  This;
                  #1
                }[\l_arg_verbs_tl]\ignorespaces #3
              }
              {Ctapr} {
                \bool_if:NT \g_debug_bool {C1.3}
                \seq_gput_left:Nn \g_arg_mru_this {Ctapr}
                \seq_gremove_duplicates:N \g_arg_mru_this
                Combining~
                \enum[\lref{\g_label_tl}]{
                  this;
                  #1
                } \proves[\l_arg_verbs_tl]~\ignorespaces #3
              }
            } {}
          }
        } {
          \str_case:VnF {\l_tmpa_tl} {
             {Ahpr} {
              \bool_if:NT \g_debug_bool {C2.1}
              \seq_gput_left:Nn \g_arg_mru_this {Ahpr}
              \seq_gremove_duplicates:N \g_arg_mru_this
              \enum:nxnn {#1} {\lref{\g_label_tl}} {-} {\BooleanTrue}
              \hence~
              \prove[\l_arg_verbs_tl]~\ignorespaces #3
            }
            {Tapr} {
              \bool_if:NT \g_debug_bool {C2.2}
              \seq_gput_left:Nn \g_arg_mru_this {Tapr}
              \seq_gremove_duplicates:N \g_arg_mru_this
              \enum[\lref{\g_label_tl}]{
                This;
                #1
              }[\l_arg_verbs_tl]\ignorespaces #3
            }
            {Ctapr} {
              \int_case:nn {\int_rand:nn {0} {1}} {
                {0} {
                  \bool_if:NT \g_debug_bool {C2.3}
                  \seq_gput_left:Nn \g_arg_mru_this {Ctapr}
                  \seq_gremove_duplicates:N \g_arg_mru_this
                  Combining~
                  \enum[\lref{\g_label_tl}]{
                    this;
                    #1
                  } \proves[\l_arg_verbs_tl]~\ignorespaces #3
                }
                {1} {
                  \bool_if:NT \g_debug_bool {C2.4}
                  \seq_gput_left:Nn \g_arg_mru_this {Ctapr}
                  \seq_gremove_duplicates:N \g_arg_mru_this
                  Combining~
                  \enum:nxnn {#1} {\lref{\g_label_tl}} {-} {\BooleanFalse}
                  \hence~
                  \proves[\l_arg_verbs_tl]~\ignorespaces #3
                }
              }
            }
          } {}
        }
      } {
        \seq_set_eq:NN \l_tmpa_seq \g_arg_mru_this
        \seq_remove_all:Nn \l_tmpa_seq {H}
        \seq_remove_all:Nn \l_tmpa_seq {Ahpr}
        \seq_get_right:NN \l_tmpa_seq \l_tmpa_tl
        \int_case:nnF {\seq_count:N \l_arg_seq} {
          {1} {
            \str_case:VnF {\l_tmpa_tl} {
              {Tapr} {
                \bool_if:NT \g_debug_bool {C3.1}
                \seq_gput_left:Nn \g_arg_mru_this {Tapr}
                \seq_gremove_duplicates:N \g_arg_mru_this
                \enum[\lref{\g_label_tl}]{
                  This;
                  #1
                }[\l_arg_verbs_tl]\ignorespaces #3
              }
              {Ctapr} {
                \bool_if:NT \g_debug_bool {C3.2}
                \seq_gput_left:Nn \g_arg_mru_this {Ctapr}
                \seq_gremove_duplicates:N \g_arg_mru_this
                Combining~
                \enum[\lref{\g_label_tl}]{
                  this;
                  #1
                } \proves[\l_arg_verbs_tl]~\ignorespaces #3
              }
            } {}
          }
        } {
          \str_case:VnF {\l_tmpa_tl} {
            {Tapr} {
              \bool_if:NT \g_debug_bool {C4.1}
              \seq_gput_left:Nn \g_arg_mru_this {Tapr}
              \seq_gremove_duplicates:N \g_arg_mru_this
              \enum[\lref{\g_label_tl}]{
                This;
                #1
              }[\l_arg_verbs_tl]\ignorespaces #3		
            }
            {Ctapr} {
              \int_case:nn {\int_rand:nn {0} {1}} {
                {0} {
                  \bool_if:NT \g_debug_bool {C4.2}
                  \seq_gput_left:Nn \g_arg_mru_this {Ctapr}
                  \seq_gremove_duplicates:N \g_arg_mru_this
                  Combining~
                  \enum[\lref{\g_label_tl}]{
                    this;
                    #1
                  } \proves[\l_arg_verbs_tl]~\ignorespaces #3		
                }
                {1} {
                  \bool_if:NT \g_debug_bool {C4.3}
                  \seq_gput_left:Nn \g_arg_mru_this {Ctapr}
                  \seq_gremove_duplicates:N \g_arg_mru_this
                  Combining~
                  \enum:nxnn {#1} {\lref{\g_label_tl}} {-} {\BooleanFalse}
                  \hence~
                  \proves[\l_arg_verbs_tl]~\ignorespaces #3    
                }
              }
            }
          } {}
        }
      }
    } {
      \tl_if_head_eq_catcode:oNTF \l_cons_tl a {
        \seq_set_eq:NN \l_tmpa_seq \g_arg_mru_this
        \seq_remove_all:Nn \l_tmpa_seq {Ctapr}
        \seq_remove_all:Nn \l_tmpa_seq {Ahpr}
        \seq_get_right:NN \l_tmpa_seq \l_tmpa_tl
        \str_case:VnF {\l_tmpa_tl} {
          {H} {
            \bool_if:NT \g_debug_bool {C5.1}
            \seq_gput_left:Nn \g_arg_mru_this {H}
            \seq_gremove_duplicates:N \g_arg_mru_this
            Hence,~we~obtain~\ignorespaces #3
          }
          {Tapr} {
            \bool_if:NT \g_debug_bool {C5.2}
            \seq_gput_left:Nn \g_arg_mru_this {Tapr}
            \seq_gremove_duplicates:N \g_arg_mru_this
            This~\proves[\l_arg_verbs_tl]~\ignorespaces #3
          }
        } {}
      } {
        \bool_if:NT \g_debug_bool {C6.1}
        \seq_gput_left:Nn \g_arg_mru_this {Tapr}
        \seq_gremove_duplicates:N \g_arg_mru_this
        This~\proves[\l_arg_verbs_tl]~\ignorespaces #3
      }
    } 
  } {
    \int_compare:nNnTF {\seq_count:N \l_arg_seq} = {0} {
      \bool_if:NTF \g_arg_start_bool {
        \bool_if:NT \g_debug_bool {C7.1}
        \Nobs\unskip
        #3
      } {
        \bool_if:NT \g_debug_bool {C7.2}
        \Moreover~
        #3
      }
    } {
      \bool_if:NTF \g_arg_start_bool {
        \bool_if:NT \g_debug_bool {C8.1}
        \tl_log:N \l_arg_verbs_tl
        \Nobs~that~
        \enum{
          #1
        }[\l_arg_verbs_tl]\ignorespaces #3
      } {
        \int_compare:nNnTF {\seq_count:N \l_arg_seq} = {1} {
          \seq_set_eq:NN \l_tmpa_seq \g_arg_mru_nothis
          \seq_remove_all:Nn \l_tmpa_seq {Nwc}
          \seq_remove_all:Nn \l_tmpa_seq {Itnswc}
          \seq_get_right:NN \l_tmpa_seq \l_tmpa_tl
        } {
          \seq_get_right:NN \g_arg_mru_nothis \l_tmpa_tl
        }
        \str_case:VnF {\l_tmpa_tl} {
          {Mo} {
            \bool_if:NT \g_debug_bool {C9.1}
            \seq_gput_left:Nn \g_arg_mru_nothis {Mo}
            \seq_gremove_duplicates:N \g_arg_mru_nothis
            Moreover,~\nobs~that~
            \enum{
              #1
            }[\l_arg_verbs_tl]\ignorespaces #3		
          }
          {Fm} {
            \bool_if:NT \g_debug_bool {C9.2}
            \seq_gput_left:Nn \g_arg_mru_nothis {Fm}
            \seq_gremove_duplicates:N \g_arg_mru_nothis
            Furthermore,~\nobs~that~
            \enum{
              #1
            }[\l_arg_verbs_tl]\ignorespaces #3		
          }
          {Ia} {
            \bool_if:NT \g_debug_bool {C9.3}
            \seq_gput_left:Nn \g_arg_mru_nothis {Ia}
            \seq_gremove_duplicates:N \g_arg_mru_nothis
            In~addition,~\nobs~that~
            \enum{
              #1
            }[\l_arg_verbs_tl]\ignorespaces #3		
          }
          {N} {
            \bool_if:NT \g_debug_bool {C9.4}
            \seq_gput_left:Nn \g_arg_mru_nothis {N}
            \seq_gremove_duplicates:N \g_arg_mru_nothis
            Next,~\nobs~that~
            \enum{
              #1
            }[\l_arg_verbs_tl]\ignorespaces #3		
          }
          {Itns} {
            \bool_if:NT \g_debug_bool {C9.5}
            \seq_gput_left:Nn \g_arg_mru_nothis {Itnswc}
            \seq_gput_left:Nn \g_arg_mru_nothis {Itns}
            \seq_gremove_duplicates:N \g_arg_mru_nothis
            In~the~next~step~we~\nobs~that~
            \enum{
              #1
            }[\l_arg_verbs_tl]\ignorespaces #3		
          }
          {Nwc} {
            \bool_if:NT \g_debug_bool {C9.6}
            \seq_gput_left:Nn \g_arg_mru_nothis {Nwc}
            \seq_gremove_duplicates:N \g_arg_mru_nothis
            Next~we~combine~
            \enum{
              #1
            }to~obtain~\ignorespaces #3
          }
          {Itnswc} {
            \bool_if:NT \g_debug_bool {C9.7}
            \seq_gput_left:Nn \g_arg_mru_nothis {Itns}
            \seq_gput_left:Nn \g_arg_mru_nothis {Itnswc}
            \seq_gremove_duplicates:N \g_arg_mru_nothis
            In~the~next~step~we~combine~
            \enum{
              #1
            }to~obtain~\ignorespaces #3
          }
        } {}
      }
    }
  }
  \bool_gset_false:N \g_arg_start_bool
  \bool_gset_false:N \l_insidearg_bool
  \cfload[.]
  \color{black}
}

\tl_new:N \g_label_tl
\tl_gset:Nn \g_label_tl { }

\NewDocumentCommand{\savelabel}{m}{
  \bool_if:NTF \l_insidearg_bool {
    \tl_gset:Nn \g_label_tl {#1}
  } {
    \tl_gset:Nn \g_label_tl { }
  }
}

\ExplSyntaxOff

\ExplSyntaxOn

\NewDocumentEnvironment {athm} {m m o} {
\str_if_eq:noTF {example} {#1} {
  \bool_gset_true:N \g_example_bool
} {
  \bool_gset_false:N \g_example_bool
}
\cfclear
\IfNoValueTF{#3}{
\begin{#1}\label{#2}\global\def\loc{#2}
}{
\begin{#1}[#3]\label{#2}\global\def\loc{#2}
}
}{
\end{#1}
}

\NewDocumentEnvironment {adef} {m} {
\begin{definition}\label{#1}\global\def\loc{#1}
}{
\end{definition}
}

\NewDocumentEnvironment{aproof} {} {
\bool_if:NTF \g_example_bool {
  \bool_gset_true:N \g_arg_start_bool
  \begin{proof}[Proof~for~\cref{\loc}]
} {
  \bool_gset_true:N \g_arg_start_bool
  \begin{proof}[Proof~of~\cref{\loc}]
}
\bool_gset_false:N \g_finishproof_bool
}{
\bool_if:NTF \g_finishproof_bool {}
{\finishproofthus}
\end{proof}
}

\NewDocumentCommand{\finishproofthus} {} {
  \bool_gset_true:N \g_finishproof_bool 
  \bool_if:NTF \g_example_bool {
    The~proof~for~\cref{\loc}~is~thus~complete.
  } {
    The~proof~of~\cref{\loc}~is~thus~complete.
  }
}
\NewDocumentCommand{\finishproofthis} {} {
  \bool_gset_true:N \g_finishproof_bool 
  \bool_if:NTF \g_example_bool {
    This~completes~the~proof~for~\cref{\loc}.
  } {
    This~completes~the~proof~of~\cref{\loc}.
  }
}

\ExplSyntaxOff

\ExplSyntaxOn

\bool_new:N \g_forexample

\NewDocumentCommand{\eg}{ o }{
\IfValueT{#1}{
\str_if_eq:noTF {fe} {#1} {
\bool_gset_true:N \g_forexample
} {\bool_gset_false:N \g_forexample}
}
\bool_if:nTF { \g_forexample } {
\bool_gset_false:N \g_forexample
for~example
}{
\bool_gset_true:N \g_forexample
for~instance
}
}

\ExplSyntaxOff

\makeatletter
\ExplSyntaxOn
\seq_new:N \g_abbrs
\prop_new:N \g_abbr_counts
\tl_new:N \l_abbr_count_tl

\ExplSyntaxOn

\NewDocumentCommand{\abbr}{m m O{#1} m m O{#4} m}{
	\expandafter\newcommand\csname#3\endcsname[1][]{
		\seq_if_in:NnTF \g_abbrs {#1} {
			\prop_get:NnN \g_abbr_counts {#1} \l_abbr_count_tl
			\prop_gput:Nnx \g_abbr_counts {#1} {\int_eval:n {\l_abbr_count_tl + 1}}
			\hyperref[#1]{#7}
		} {
			\seq_gput_left:Nn \g_abbrs {#1}
			\prop_gput:Nnn \g_abbr_counts {#1} {1}
			\expandafter\gdef\csname#1@def\endcsname{#2}
			\phantomsection\label{#1}
			\str_if_eq:nnTF{##1}{}{\emph{#2}}{##1}~(\hyperref[#1]{#7})
		}
	}
	\expandafter\newcommand\csname#6\endcsname[1][]{
		\seq_if_in:NnTF \g_abbrs {#1} {
			\prop_get:NnN \g_abbr_counts {#1} \l_abbr_count_tl
			\prop_gput:Nnx \g_abbr_counts {#1} {\int_eval:n {\l_abbr_count_tl + 1}}
			\hyperref[#1]{#4}
		} {
			\expandafter\gdef\csname#1@def\endcsname{#5}
			\seq_gput_left:Nn \g_abbrs {#1}
			\prop_gput:Nnn \g_abbr_counts {#1} {1}
			\phantomsection\label{#1}
			\str_if_eq:nnTF{##1}{}{\emph{#5}}{##1}~(\hyperref[#1]{#4})
		}
	}
}

\ExplSyntaxOff

\makeatother
\abbr{KL}{Kurdyka-{\L}ojasiewicz}{KLs}{Kurdyka-{\L}ojasiewicz}{KL}
\abbr{PL}{Polyak-{\L}ojasiewicz}{PLs}{Polyak-{\L}ojasiewicz}{PL}
\abbr{HB}{Heavy-Ball}{HBs}{Heavy-Ball}{HB}
\abbr{OP}{optimization problem}{OPs}{optimization problems}{OP}
\abbr{ODE}{ordinary differential equation}{ODEs}{ordinary differential equations}{ODE}
\abbr{SOP}{stochastic optimization problem}{SOPs}{stochastic optimization problems}{SOP}
\abbr{NAG}{Nesterov accelerated gradient}{NAGs}{Nesterov accelerated gradient}{NAG}
\abbr{Adam}{adaptive moment estimation}{Adams}{adaptive moment estimation}{Adam}
\abbr{Adamax}{adaptive moment estimation maximum}{Adamaxs}{Adaptive moment estimation maximum}{Adamax}
\abbr{AdaBelief}{adaptive belief}{AdaBeliefs}{adaptive belief}{AdaBelief}
\abbr{UGD}{unified gradient descend}{UGDs}{unified gradient descend}{UGD}
\abbr{Nadam}{Nesterov accelerated adaptive moment estimation}{Nadams}{Nesterov accelerated adaptive moment estimation}{Nadam}
\abbr{Nadamax}{Nesterov accelerated adaptive moment estimation maximum}{Nadamaxs}{Nesterov accelerated adaptive moment estimation maximum}{Nadamax}
\abbr{Adan}{adaptive Nesterov momentum}{Adans}{Adaptive Nesterov momentum}{Adan}
\abbr{RMSprop}{root mean square propagation}{RMSprops}{the root mean square propagation}{RMSprop}
\abbr{ANN}{artificial neural network}{ANNs}{artificial neural networks}{ANN}
\abbr{DNN}{deep neural network}{DNNs}{deep neural networks}{DNN}
\abbr{GD}{Gradient descent}{GDs}{gradient descent}{GD}
\abbr{GF}{gradient flow}{GFs}{gradient flow}{GF}
\abbr{SGD}{stochastic gradient descent}{SGDs}{stochastic gradient descent}{SGD}
\abbr{iid}{independent and identically distributed}{i.i.d}{independent and identically distributed}{i.i.d.}
\abbr{DL}{Deep learning}{DLs}{Deep learning}{DL}
\abbr{AI}{artificial intelligence}{AIs}{artificial intelligence}{AI}
\abbr{LLM}{large language model}{LLMs}{large language models}{LLM}
\abbr{PDE}{partial differential equation}{PDEs}{partial differential equations}{PDE}
\abbr{as}{almost surely}{ass}{almost surely}{a.s.}
\abbr{pas}{$\P$-almost surely}{pass}{almost surely}{$\P$-a.s.}
\abbr{ReLUs}{rectified linear unit}{ReLU}{rectified linear unit}
\makeatother

\begin{document}

\title{Unified convergence analysis for \\
gradient descent optimization methods \\
in the training of deep neural networks}

\author{Shokhrukh Ibragimov$^{1}$ and Arnulf Jentzen$^{2,3}$
\bigskip\\
\small{$^1$ Applied Mathematics: Institute for Analysis and Numerics,}\vspace{-0.1cm}\\
\small{University of M\"unster, Germany; e-mail: \texttt{sibragim}\textcircled{\texttt{a}}\texttt{uni-muenster.de}}\smallskip\\
\small{$^2$ School of Data Science and School of Artificial Intelligence,}\vspace{-0.1cm}\\
\small{The Chinese University of Hong Kong, Shenzhen (CUHK-Shenzhen),}\vspace{-0.1cm}\\
\small{China; e-mail: \texttt{ajentzen}\textcircled{\texttt{a}}\texttt{cuhk.edu.cn}}\smallskip\\
\small{$^3$ Applied Mathematics: Institute for Analysis and Numerics,}\vspace{-0.1cm}\\
\small{University of M\"unster, Germany; e-mail: \texttt{ajentzen}\textcircled{\texttt{a}}\texttt{uni-muenster.de}}
}

\date{\today}

\maketitle

\begin{abstract}
{\footnotesize Gradient based optimization methods are nowadays the methods of choice for training deep neural networks (DNNs) in artificial intelligence (AI) systems. In practically relevant DNN training problems, one does usually not apply the standard gradient descent (GD) optimization method but instead one employs suitable sophisticated GD optimization
methods, which incorporate adaptivity and/or acceleration techniques, such as the famous Adam optimizer. It is a key contribution of this work to provide a general unified convergence analysis for GD optimization methods in the training of DNNs with analytic activations such as the softplus and the popular Gaussian error linear unit (GeLU) activation. Our general unified convergence result applies to a large class of gradient based optimization methods such as the standard GD, the momentum, the Nesterov accelerated gradient (NAG), the RMSprop, the Adam, the Adamax, the Nadam, the Nadamax, the Adan, the AdaBelief, the AMSGrad, and the Yogi optimizers. Our analysis employs the theory of Kurdyka-{\L}ojasiewicz (KL) inequalities to establish convergence to critical points in the training of DNNs. To the best of our knowledge, the generality of our convergence analysis is also just in the special situation of the Adam optimizer a new contribution to the literature on the analysis of AI optimization algorithms.}
\end{abstract}

\newpage
\tableofcontents

\pagebreak
\section{Introduction}
\label{sec:introduction}

\GD[\emph{Gradient descent}]\ based optimization methods are nowadays the standard methods for the training of \DNNs\ in \AI\ systems. In the training of large scale \AI\ models, the standard \GD\ method and its stochastic counterpart \cite{GarrigosGowerSGD}, respectively, are often not the employed optimization methods but instead suitable sophisticated variants of the standard \GD\ method involving adaptivity and acceleration are considered \cite{JentzenBookDeepLearning2023,SRuderGDOverview}. Prominent adaptive and/or accelerated \GD\ based optimization methods are, \eg,
\begin{enumerate}[label=(\roman*)]
\item\label{intro:momentum}
the momentum \cite{POLYAK19641},

\item\label{intro:nagd}
the \NAG\ \cite{Nesterov1983AMF},

\item\label{intro:RMSprop}
the \RMSprop\ \cite{HintonRMSProp},

\item\label{intro:Adam}
the \Adam\ \cite{KingmaBaAdam},

\item\label{intro:Adamax}
the \Adamax\ \cite{KingmaBaAdam},

\item\label{intro:Nadam}
the \Nadam\ \cite{DozatNadam2},

\item\label{intro:Nadamax}
the \Nadamax\ \cite{DozatNadam2},

\item\label{intro:Adan}
the vanilla \Adan\ \cite{Adan},

\item\label{intro:AdaBelief}
the \AdaBelief\ \cite{ZhuangAdaBelief2020},

\item\label{intro:AMSGrad}
the AMSGrad \cite{ReddiKale2019}, and

\item\label{intro:Yogi}
the Yogi \cite{Yogi_NEURIPS2018}
\end{enumerate}
optimizers. The most popular of such optimization methods is presumably the \Adam\ optimizer (see \cref{intro:Adam} above) proposed in 2014 by Kingma \& Ba~\cite{KingmaBaAdam}. Despite the popularity and the success of \Adam, it remained an open question to prove (or disprove) strong convergence of every bounded \Adam\ optimization trajectory in the training of \DNNs.

In this work we partially solve this problem by establishing -- in the training of fully-connected feedforward \DNNs\ with analytic activations -- convergence with convergence rates not only for every bounded optimization trajectory of \Adam\ but for each of the above mentioned optimizers \ref{intro:momentum}--\ref{intro:Yogi}. Specifically, it is a key contribution of the main result of this work, see \cref{cor:gen_convergence_rate_Loja_l_rate} in \cref{subsec:convergence_UGD_KL_Lipschitz} below, to develop a \emph{unified framework}
\begin{enumerate}[label=(\alph*)]
\item which includes \emph{each} of the above mentioned optimization methods \ref{intro:momentum}--\ref{intro:Yogi} as a special case and

\item under which every bounded trajectory of the optimizer converges with convergence rates to a critical point for general \KL\ objective functions with \emph{locally Lipschitz continuous} gradients.
\end{enumerate}
In \cref{def:Loja} in \cref{subsec:KL_functions} below we recall the notion of a \KL\ function (cf., \eg, \cite[Definition~1]{MR3785672}, \cite[Definition~9]{pmlr-v129-barakat20a}, and \cite[Definition~9.1.2]{JentzenBookDeepLearning2023}). \KL\ objective functions cover, \eg, the training of \DNNs\ with analytic activation functions as a special case; see \cref{subsec:dnns} and \cref{thm:Adam_for_dnns_2} below for details.

To briefly illustrate the contribution of this work in this introductory section, we present in \cref{thm_intro_strong_convergence_1} below in the following subsection a special case of our general unified convergence analysis in \cref{cor:gen_convergence_rate_Loja_l_rate}.

\subsection{Main result: Unified error analysis for gradient descent (GD) optimizers}

In \cref{thm_intro_strong_convergence_1} the natural number $\fd \in \N = \{ 1, 2, 3, \dots \}$ represents the dimensionality of the considered optimization problem and the function $\cL \colon \R^{\fd} \to \R$ is the objective function of the optimization problem, that is, $\cL \colon \R^{\fd} \to \R$ is the function that we intend to minimize in \cref{thm_intro_strong_convergence_1}.

\cfclear
\begin{savenotes}
\begin{samepage}
\begin{tcolorbox}[colback=white!95!gray,
                  colframe=black,
                  boxrule=0.5pt,
                  sharp corners,
                  enhanced,
                 ]
\begin{athm}{theorem}{thm_intro_strong_convergence_1}
Let $\fd \in \N$, $\alpha \in [0, 1)$, $\lrexpo \in (\nicefrac{3}{4}, 1]$, let $\cL \colon \allowbreak \R^{\fd} \allowbreak \to \R$ be a\cfadd{def:Loja} \KL[Kurdyka-{\L}ojasiewicz] function, assume that $\nabla \cL$ is locally Lipschitz continuous, let $\bbA \colon \N \to \R^{\fd \times \fd}$, $\Theta \colon \N_0 \allowbreak \to \allowbreak \R^\fd$, $\mu \colon \N \to \R^{\fd}$, and $p \colon \N \to \R$ be bounded, let $\m \allowbreak \colon \allowbreak \N_0 \allowbreak \to \R^\fd$ satisfy for all $n \in \N$ that
\begin{gather}\llabel{m}
\textstyle \textcolor{magenta}{\m_n = \alpha \m_{n - 1} + (1 - \alpha) \bigl[(\nabla \cL)(\Theta_{n - 1}) + n^{- \lrexpo} \mu_n\bigr]} \\ \llabel{Theta}
\textstyle \andqq \textcolor{magenta}{\Theta_n = \Theta_{n - 1} - n^{- \lrexpo} \bbA_n \bigl[p_n \m_n + (1 - p_n) (\nabla \cL)(\Theta_{n - 1})\bigr]},
\end{gather}
and assume for all $n \in \N$ that $\bbA_n - \bbI_{\fd}$ is symmetric positive semi-definite \cfout. Then there exist \textcolor{magenta}{$\vartheta \allowbreak \in \allowbreak \R^{\fd}$}, \textcolor{magenta}{$\rho \allowbreak \in (0, \allowbreak \infty)$} which satisfy\footnote{\Nobs that for all $m \in \N$, $x = (x_1, \dots, x_m) \in \R^m$ it holds that $\norm{x} = (\sum_{i = 1}^m \abs{x_i}^2)^{\nicefrac{1}{2}}$ (standard norm).} for all \textcolor{magenta}{$n \in \N$} that
\begin{equation}\llabel{result}
\textstyle \textcolor{magenta}{\norm{(\nabla \cL) (\vartheta)} + \abs{\cL(\Theta_n) - \cL(\vartheta)} + \norm{\Theta_n - \vartheta}^{\rho} \le \rho \bigl[\sum_{j = 1}^{n} j^{- \lrexpo}\bigr]^{- 1}} \dott
\end{equation}
\end{athm}
\end{tcolorbox}
\end{samepage}
\end{savenotes}

\cref{thm_intro_strong_convergence_1} is a direct consequence of our more general convergence result in \cref{cor:gen_convergence_rate_Loja_l_rate}, in which the momentum decay factor $\alpha \in [0, 1)$ in \cref{thm_intro_strong_convergence_1.m} may depend on $n \in \N$.

The process $(\Theta_n)_{n \in \N_0}$ in \cref{thm_intro_strong_convergence_1.Theta} in \cref{thm_intro_strong_convergence_1} is the considered optimization process. The process $(\m_n)_{n \in \N}$ in \cref{thm_intro_strong_convergence_1.m} is an abstract momentum process and the process $(\mu_n)_{n \in \N}$ in \cref{thm_intro_strong_convergence_1.m} allows us to handle perturbations from the abstract momentum process. The matrix valued process $\bbA = (\bbA_n)_{n \in \N} \colon \N \to \R^{\fd \times \fd}$ in \cref{thm_intro_strong_convergence_1.Theta} in \cref{thm_intro_strong_convergence_1} allows us to cover the learning rate as well as other adaptivity features of the considered optimizer, such as adaptive control of the learning rates as in \RMSprop\ and \Adam\ (cf.\ \cref{thm:Adam_for_dnns_2.Theta} below for details). We label the perspective in \cref{thm_intro_strong_convergence_1.m,thm_intro_strong_convergence_1.Theta} to handle a large class of optimizers within a unified framework as \UGD\ approach and we refer to \cref{subsec:UGD} and \cref{section:UGD_applications} below for details on this \UGD\ approach.

\cref{thm_intro_strong_convergence_1} proves that the considered optimization process $(\Theta_n)_{n \in \N_0}$ converges with a suitable speed of convergence to a \emph{critical point} $\vartheta \in \R^{\fd}$ of the objective function $\cL \colon \R^{\fd} \to \R$ (a critical point of $\cL$ is a point in $\R^{\fd}$ at which the gradient of the objective function vanishes). The term $[\sum_{j = 1}^n j^{- \lrexpo}]^{- 1}$ on the right hand side of \cref{thm_intro_strong_convergence_1.result} ensures that the optimization process converges with a suitable (polynomial) speed of convergence; see, \eg, \cref{lemma:l_rate_specific} in \cref{subsec:convergence_UGD_KL_Lipschitz} below for details. In particular, we observe that $\lim_{n \to \infty} ([\sum_{j = 1}^n j^{- \lrexpo}]^{- 1}) = 0$.

A key contribution of \cref{thm_intro_strong_convergence_1} (and the more general results in this work, respectively) is that the setup in \cref{thm_intro_strong_convergence_1.m}--\cref{thm_intro_strong_convergence_1.Theta} in \cref{thm_intro_strong_convergence_1} is so general that it covers the \Adam\ optimizer and \emph{each} of the other above mentioned optimization methods \ref{intro:momentum}--\ref{intro:Yogi} as special cases; see \cref{section:UGD_applications} below for details.

In \cref{thm_intro_strong_convergence_1} and our more general convergence results we employ the assumption that the optimization process is bounded. We refer, \eg, to \cite{DereichDoJentzen_uniformbound_Adam} (cf., \eg, also \cite{DereichGraeberJentzenRiekert_asymtoticstability_Adam}) for sufficient conditions and example optimization problems that ensure that the \Adam\ optimization process is bounded and thus fulfills this requirement.

The objective function $\cL \colon \R^{\fd} \to \R$ in \cref{thm_intro_strong_convergence_1} (the function that we intend to minimize in \cref{thm_intro_strong_convergence_1}) is only assumed to be a \emph{\KL\ function} (cf.\ \cref{def:Loja}) and to have a \emph{locally} Lipschitz continuous first derivative. These assumptions, in particular, cover the training of \DNNs\ with analytic activation functions (see \cref{subsec:dnns} below for details).

To sketch the generality of the setup in \cref{thm_intro_strong_convergence_1}, we now illustrate \cref{thm_intro_strong_convergence_1} in the special situation where \Adam\ is applied for the training of \DNNs\ with analytic activations. This is precisely the subject of \cref{thm:Adam_for_dnns_2} presented in the next subsection.

\subsection{Application of the main result: Adam training deep neural networks (DNNs)}

In \cref{thm:Adam_for_dnns_2} the natural number $M \in \N$ is the number of input-output data pairs, the real vectors $\fx_1, \fx_2, \dots, \fx_M \in \R^{\ell_0}$ represent the input data, and the real vectors $\fy_1, \allowbreak \fy_2, \allowbreak \dots, \allowbreak \fy_M \allowbreak \in \R^{\ell_L}$ represent the output data of the considered learning problem. We first again present \cref{thm:Adam_for_dnns_2} in a completely self-contained way with all mathematical details and thereafter provide explaining sentences.

\cfclear
\begin{savenotes}
\begin{samepage}
\begin{tcolorbox}[colback=white!95!gray,
                  colframe=black,
                  boxrule=0.5pt,
                  sharp corners,
                  enhanced,
                 ]
\begin{athm}{cor}{thm:Adam_for_dnns_2}[\resname{\Adam\ training \DNNs}]
Let $\fd, L, M \in \N$, $\ell_0, \allowbreak \ell_1, \allowbreak \dots, \allowbreak \ell_L \allowbreak \in \N$, $\fx_1, \allowbreak \fx_2, \allowbreak \dots, \allowbreak \fx_M \allowbreak \in \allowbreak \R^{\ell_0}$, $\fy_1, \allowbreak \fy_2, \allowbreak \dots, \allowbreak \fy_M \allowbreak \in \allowbreak \R^{\ell_L}$ satisfy $\fd = \allowbreak \sum_{i = 1}^{L} \ell_i (\ell_{i - 1} + 1)$, let $\act \colon \R \to \R$ be analytic, for every $\theta = \allowbreak (\theta_1, \allowbreak \dots, \allowbreak \theta_{\fd}) \allowbreak \in \R^{\fd}$ let $\cN^{k, \theta} = \allowbreak (\cN^{k, \theta}_1, \allowbreak \dots, \allowbreak \cN^{k, \theta}_{\ell_k}) \colon \allowbreak \R^{\ell_0} \allowbreak \to \R^{\ell_k}$, $k \in \allowbreak \{0, \allowbreak 1, \allowbreak \dots, \allowbreak L\}$, satisfy for all $k \in \allowbreak \{0, \allowbreak 1, \allowbreak \dots, \allowbreak L - 1\}$, $x = (x_1, \allowbreak \dots, \allowbreak x_{\ell_0}) \allowbreak \in \allowbreak \R^{\ell_0}$, $i \in \allowbreak \{1, \allowbreak 2, \allowbreak \dots, \allowbreak \ell_{k + 1}\}$ that
\begin{multline}
\textstyle \textcolor{magenta}{\cN^{k + 1, \theta}_i(x) = \theta_{\ell_{k + 1} \ell_{k} + i + \sum_{h = 1}^{k} \ell_h (\ell_{h - 1} + 1)}} \\
\textstyle \textcolor{magenta}{+ \sum_{j = 1}^{\ell_{k}} \theta_{(i - 1) \ell_{k} + j + \sum_{h = 1}^{k} \ell_h (\ell_{h - 1} + 1)} \bigl[x_j \mathbbm{1}_{\{0\}}(k) + \act(\cN^{k, \theta}_j(x)) \mathbbm{1}_{\N}(k)\bigr]},
\end{multline}
let $\cL \colon \R^{\fd} \to \R$ satisfy for all $\theta \in \R^{\fd}$ that $\textcolor{magenta}{\cL(\theta) = \frac{1}{M} \sum_{m = 1}^{M} \norm{\cN^{L, \theta}(\fx_m) - \fy_m}^2}$, let $\beta_1, \allowbreak \beta_2 \in [0, \allowbreak 1)$, $\eps \allowbreak \in (0, \allowbreak \infty)$, $\lrexpo \allowbreak \in (\nicefrac{3}{4}, \allowbreak 1]$, let $\gamma \colon \N \allowbreak \to (0, \allowbreak \infty)$ and $\Theta\Index{k}{} \allowbreak = (\Theta\Index{k}[1]{}, \allowbreak \dots, \allowbreak \Theta\Index{k}[\fd]{}) \allowbreak \colon \allowbreak \N_0 \allowbreak \to \allowbreak \R^\fd$, $k \in \N_0$, satisfy for all $k \in \allowbreak \{1, \allowbreak 2\}$, $n \in \N$, $i \in \allowbreak \{1, \allowbreak 2, \allowbreak \dots, \allowbreak \fd\}$ that
\begin{gather}\llabel{m_M}
\textstyle \textcolor{magenta}{\Theta\Index{k}[i]{n} = \beta_k \Theta\Index{k}[i]{n - 1} + (1 - \beta_k) \bigl[(\nabla \cL)^i(\Theta\Index{0}{n - 1})\bigr]^k} \\ \llabel{Theta}
\textstyle \andqq \textcolor{magenta}{\Theta\Index{0}[i]{n} = \Theta\Index{0}[i]{n - 1} - \gamma_n \bigl[\frac{\Theta\Index{1}[i]{n}}{1 - (\beta_1)^n}\bigr] \Bigl[\eps + \bigl[\frac{\abs{\Theta\Index{2}[i]{n}}}{1 - (\beta_2)^n}\bigr]^{\nicefrac{1}{2}}\Bigr]^{-1}},
\end{gather}
and assume $\sup_{n \in \N} (\norm{\Theta\Index{0}{n}} + \sum_{k = -1}^{1} (\gamma_n n^{\lrexpo})^k) < \infty$ \cfout. Then there exist $\textcolor{magenta}{\vartheta \allowbreak \in \allowbreak \R^{\fd}}$, $\textcolor{magenta}{\rho \in (0, \allowbreak \infty)}$ which satisfy for all $\textcolor{magenta}{n \in \N}$ that
\begin{equation}
\textstyle \textcolor{magenta}{\norm{(\nabla \cL) (\vartheta)} + \abs{\cL(\Theta\Index{0}{n}) - \cL(\vartheta)} + \norm{\Theta\Index{0}{n} - \vartheta}^{\rho} \le \rho \bigl[\sum_{j = 1}^{n} j^{- \lrexpo}\bigr]^{-1}} \dott
\end{equation}
\end{athm}
\end{tcolorbox}
\end{samepage}
\end{savenotes}

\cref{thm:Adam_for_dnns_2} is proved as \cref{cor:Adam_for_dnns_2} in \cref{subsec:Adam} below. \cref{cor:Adam_for_dnns_2} follows from \cref{cor:gen_convergence_rate_Loja_l_rate} below and \cref{thm_intro_strong_convergence_1} above, respectively.

In \cref{thm:Adam_for_dnns_2} we do not consider the general setup in \cref{thm_intro_strong_convergence_1.m}--\cref{thm_intro_strong_convergence_1.Theta} that covers a large class of optimizers but instead we exclusively restrict ourselves to the \Adam\ optimizer in \cref{thm:Adam_for_dnns_2.m_M}--\cref{thm:Adam_for_dnns_2.Theta}, which is currently presumably the most popular optimizer for the training of \AI\ systems. In \cref{thm:Adam_for_dnns_2} we study the training of \DNNs\ by means of \Adam\ and for every $\theta \in \R^{\fd}$ we have that $\cN^{L, \theta} \colon \R^{\ell_0} \to \R^{\ell_L}$ is the realization function of the fully-connected feedforward \DNN\ with the \DNN\ parameter vector $\theta$, with the activation function $\act \colon \R \to \R$, and with the architecture consisting of $L + 1$ layers ($L - 1$ hidden layers)
\begin{itemize}
\item with $\ell_0$ neurons on the input layer (on the $1^{\text{st}}$ layer),

\item with $\ell_1$ neurons on the $1^{\text{st}}$ hidden layer (on the $2^{\text{nd}}$ layer),

\item with $\ell_2$ neurons on the $2^{\text{nd}}$ hidden layer (on the $3^{\text{rd}}$ layer),

\dots
\item with $\ell_{L - 1}$ neurons on the $(L - 1)^{\text{th}}$ hidden layer (on the $L^{\text{th}}$ layer), and

\item with $\ell_L$ neurons on the output layer (on the $(L + 1)^{\text{th}}$ layer).
\end{itemize}
The natural number $\fd \in \N$ satisfying $\fd = \allowbreak \smallsum_{i = 1}^{L} \ell_i (\ell_{i - 1} + 1)$ in \cref{thm:Adam_for_dnns_2} represents the overall number of parameters of the considered \DNNs. The function $\cL \colon \R^{\fd} \to \R$ in \cref{thm:Adam_for_dnns_2} specifies the objective function of the considered minimization problem and the \Adam\ optimization process is specified in \cref{thm:Adam_for_dnns_2.m_M}--\cref{thm:Adam_for_dnns_2.Theta} in \cref{thm:Adam_for_dnns_2}.

\subsection{Literature overview}
\label{subsec:literature}

In this subsection we briefly review selected further results from the literature that theoretically analyze \GD\ optimization algorithms, where we particularly focus on convergence results for accelerated and adaptive \GD\ optimization methods.

\subsubsection*{Adaptive and accelerated \GD\ optimization methods beyond \Adam}

Since \cite{ReddiKale2019} has demonstrated that \Adam\ can fail to converge even on simple one-dimensional convex objectives, a proliferation of \Adam\ variants has emerged. \emph{Weak convergence} -- convergence of function values at the iterates -- has been established for convex objective functions in methods such as \AdaBelief\ ~\cite{ZhuangAdaBelief2020}, FastAdaBelief~\cite{MR4637063}, and SAdam~\cite{WangSAdam2020}. Relaxing convexity, analyses of variants like AdamL~\cite{XiaAdamL2023} achieve weak convergence by assuming the \PL\ condition. In broader non-convex regimes, weak convergence proofs are achieved assuming (lower) boundedness of the objective function and/or Lipschitz continuous gradients, as demonstrated in the analyses of \Adan\ ~\cite{Adan}, Yogi~\cite{Yogi_NEURIPS2018}, AdaBound~\cite{LIU2022300,LuoAdaBoundAMSBound2019}, AMSBound~\cite{LuoAdaBoundAMSBound2019}, SignSGD~\cite{SignSGD18}, generalized SignSGD~\cite{Crawshaw2022}, AdaNorm~\cite{AdaNorm22,FawTziotis2022}, AdaGrad~\cite{HongLin2024,WangZhang2023,Defossez2022}, and \SGD\ with momentum and difference (\SGD\ (MD)) \cite{YuanNonAdaptiveOPT2022}. Beyond weak guarantees, \emph{strong convergence} -- convergence of the iterates themselves -- has been established in highly restricted settings, such as for quadratic objective functions in, \eg, AdaLoss~\cite{WuAdaLoss2021}. Moreover, other variants of \Adam\ have been introduced based purely on numerical simulations and empirical performance (see, \eg, RAdam~\cite{LiuJiang2019} and EAdam~\cite{YuanEAdam2020}).

Beyond adaptive optimizers, foundational algorithms such as \SGD\ and its momentum-based extensions -- including the \HB\ and \NAG\ methods -- have been extensively analyzed. However, for general smooth objectives, the established theoretical guarantees predominantly remain weak. In non-convex settings, weak convergence and associated rates have been established for \SGD\ \cite{ZhangFang2020} and momentum \GD\ \cite{YangMomentum2016,NEURIPS2020_d3f5d4de,pmlr-v178-liu22d} under bounded variance and standard or generalized Lipschitz continuous gradients. Similarly, under global convexity, analyses of the \HB\ method provide weak convergence guarantees in both stochastic and deterministic regimes \cite{SunHBconvergence2018, 7330562}.

Strong convergence for these foundational momentum algorithms is largely confined to highly restricted structural settings. For instance, (global) linear strong convergence for momentum methods is established, but for convex quadratic objectives, such as least-squares regression problems \cite{LoizouSHB2017, MR4174637}. Comparable strong guarantees and lower error bounds for standard \SGD\ also rely on quadratic landscapes \cite{JENTZEN2020101438}, or demand rigid structural constraints when applied to broader convex functions \cite{pmlr-v134-sebbouh21a}.

To bridge the gap to strong convergence in more general non-convex landscapes, recent analyses have increasingly leveraged geometric properties such as the \PL\ and \KL\ conditions. Strong convergence of deterministic momentum \GD\ to a local minimizer of semialgebraic objective functions (which inherently satisfy the \KL\ property) is established in \cite{MR4663539} with constant hyperparameters. However, this relies on restrictive assumptions regarding the Hessian of the objective function -- specifically imposing the strict saddle property by requiring a strictly negative eigenvalue at every critical point that is not a local minimum point. Moreover, strong convergence of \SGD\ has been proven for analytic functions satisfying the \KL\ property \cite{MR3315611}, while momentum \SGD\ has been investigated under both \PL\ \cite{GessKassingMomentum2023} and \KL\ \cite{DereichKassing,QiuMaMilzarek_momentum_KL_24} geometries. Furthermore, \cite{MR4298987} established strong convergence for a generalized class of stochastic momentum methods which encompasses \Adam. However, these strong convergence results enforce coupling between the momentum decay parameters $(\alpha_n)_{n \in \N} \subseteq [0, 1]$ and non-summable learning rates $(\gamma_n)_{n \in \N} \subseteq (0, \infty)$, \eg, setting $\alpha_n = 1 - r_n \gamma_n$ for some $(r_n)_{n \in \N} \subseteq (0, \infty)$ in \cite{GessKassingMomentum2023,DereichKassing,MR4298987} and $\alpha_n = \lambda \gamma_{n - 1} (\gamma_n)^{- 1}$ for arbitrary $\lambda \in [0, 1)$ in \cite{QiuMaMilzarek_momentum_KL_24}. Such coupled restrictions limit their generality, in contrast to our deterministic framework, which accommodates arbitrary hyperparameter decoupling.

\subsubsection*{Convergence analyses for the \Adam\ optimizer}

Under standard assumptions of (lower) boundedness of the objective function and/or Lipschitz continuous gradients, weak convergence -- typically formulated as bounds for the gradient or the objective function evaluated at the iterates -- has been widely established for \Adam\ \cite{ChenLiuSunHong2018,ZouShen2019,Defossez2022,ZhangChen2022,HeAdamConvergence2023} and unified \Adam-type frameworks \cite{Jiang2023UAdam}. Weak convergence bounds for \Adam\ have been derived under generalized $(L_0, L_\rho)$-smoothness conditions \cite{Li2023}, and frameworks encompassing the broader \Adam\ family have been analyzed in non-smooth landscapes using subdifferentials \cite{MR4723898}.

Strong convergence for adaptive gradient methods remains rare and is confined to highly restricted structural settings. For example, local strong convergence for \Adam\ has been established, but only within local regions where the objective function behaves as strongly convex (exhibiting a positive definite Hessian) \cite{DereichJentzenRiekert_sharpconvergencerates_Adam,BockAdamLocalConvergence2019}. Scaling to global guarantees, recent unconditional error analyses \cite{DereichDoJentzen_uniformbound_Adam} have secured strong convergence for \Adam, but require the objective function to be globally strongly convex. Similarly, while \cite{GodichonBaggioni2023} establishes strong convergence for \SGD\ with adaptive scaling matrix $(\bbA_n)_{n \in \N}$ -- similar to our formulation in \cref{thm_intro_strong_convergence_1.Theta} in \cref{thm_intro_strong_convergence_1} when $p_n = 0$ or $\abs{p_n - 1} + \alpha + \norm{\mu_n} = 0$ for all $n \in \N$ -- using smallest eigenvalue control $\liminf_{n \to \infty} \lambda_{\min} (\bbA_n) > 0$, their analysis is restricted to quasi-strongly convex objective functions.

Seeking broader generality beyond strict convexity, Barakat \& Bianchi investigated \Adam, initially deriving weak convergence for an objective function with a Lipschitz continuous gradient \cite{pmlr-v129-barakat20a}, and subsequently proving strong convergence under the \KL\ property to a critical point via an associated \ODE\ limit \cite{Barakat2021}. While decaying learning rates are a fundamental necessity for achieving stochastic convergence -- as demonstrated by non-convergence results for non-decaying learning rate schedules even for quadratic optimization \cite{DereichGraeberJentzen_Adam_nonconvergence} -- achieving strong convergence via continuous-time \ODE\ limits \cite{Barakat2021}, and the aforementioned weak convergence results \cite{Jiang2023UAdam,Li2023,MR4723898}, additionally require the first-order momentum decay parameters $(\alpha_n)_{n \in \N} \subseteq [0, 1]$ to converge to $1$ or to be coupled with non-summable learning rates $(\gamma_n)_{n \in \N} \subseteq (0, \infty)$, setting, \eg, $\alpha_n = 1 - r_n \gamma_n$ for some $(r_n)_{n \in \N} \subseteq (0, \infty)$. Consequently, an analysis that secures strong convergence for \Adam\ in general non-convex \KL\ landscapes -- without, \eg, chaining the momentum hyperparameters to the learning rate -- remained an open problem that our unified framework resolved.

\subsection{Structure of this article}

The remainder of this work is organized as follows. In \cref{section:weak_convergence} we establish under suitable assumptions within a unified framework, which we refer to as \UGD\ approach, convergence of the evaluation of the objective function at the considered optimization process to a critical level of the objective function (to the evaluation of the objective function at a critical point of the objective function); see \cref{prop:gen_convergence_ohne_Loja} for the main result of \cref{section:weak_convergence}. In \cref{section:UGD_strong_convergence} we employ this weak convergence result for the \UGD\ method from \cref{prop:gen_convergence_ohne_Loja} to establish for every bounded optimization trajectory of the \UGD\ method strong convergence to a critical point for general \KL\ objective functions with locally Lipschitz continuous gradients; see \cref{theorem:gen_convergence_Loja_no_delta}, \cref{cor:gen_convergence_rate_Loja_l_rate}, \cref{cor:gen_convergence_rate_Loja_l_rate_TEMP}, and \cref{cor:gen_convergence_rate_Loja_gen_l_rate_Lipschitz_alpha}. Moreover, in \cref{subsec:dnns} we illustrate the conclusion of the strong convergence result in \cref{cor:gen_convergence_rate_Loja_gen_l_rate_Lipschitz_alpha} in the situation of the training of \DNNs. Finally, in \cref{section:UGD_applications} we apply \cref{cor:gen_convergence_rate_Loja_gen_l_rate_Lipschitz_alpha} to several concrete optimization methods such as each of the optimizers in \ref{intro:momentum}--\ref{intro:Yogi} above.

\section{Weak convergence for GD op\-ti\-mi\-za\-tion methods}
\label{section:weak_convergence}

In this section we study each of the optimization methods \ref{intro:momentum}--\ref{intro:Yogi} from \cref{sec:introduction} within a \emph{unified framework} -- which we refer to as \UGD\ method -- presented in \cref{setting:UGD} in \cref{subsec:UGD} below. In particular, in \cref{prop:gen_convergence_ohne_Loja} below (which is the main result of this section) we establish within this framework under suitable assumptions convergence of the evaluation of the objective function at the considered optimization process to a critical level of the objective function (to the evaluation of the objective function at a critical point of the objective function).

This kind of weak convergence result for \GD\ methods is crucially used in the proof of our main convergence result in \cref{cor:gen_convergence_rate_Loja_l_rate} in \cref{section:UGD_strong_convergence} below and thus also in our proof of \cref{thm_intro_strong_convergence_1} in the introduction. The arguments in our proof of \cref{prop:gen_convergence_ohne_Loja} are in parts inspired by the arguments, \eg, in the proof of \cite[Theorem~3.1]{GodichonBaggioni2023}.

\subsection{Unified gradient descent (UGD) method}
\label{subsec:UGD}

\cfclear
\begin{savenotes}
\begin{samepage}
\begin{tcolorbox}[colback=white!95!gray,
                  colframe=black,
                  boxrule=0.5pt,
                  sharp corners,
                  enhanced,
                 ]
\begin{setting}[\resname{\UGD\ method}]\label{setting:UGD}
Let $\fd \in \N$, $\delta \in (0, \infty)$, $\cL \in C^1(\R^{\fd}, \R)$ and let $\alpha \colon \N \allowbreak \to [0, \allowbreak 1]$, $p \colon \N \allowbreak \to [0, \allowbreak 1]$, $\gamma \colon \N \to (0, \infty)$, $\bbA \colon \N \to \R^{\fd \times \fd}$, $\m \colon \N_0 \to \R^{\fd}$, $\mu \colon \N \to \R^{\fd}$, and $\Theta \colon \N_0 \to \R^{\fd}$ satisfy for all $n \in \N$ that
\begin{gather}
\textstyle \m_n = \alpha_n \m_{n - 1} + (1 - \alpha_n) (\nabla \cL)(\Theta_{n - 1}) + (\gamma_n)^{\delta} \mu_n \\
\textstyle \andqq \Theta_n = \Theta_{n - 1} - \gamma_n \bbA_n \bigl[p_n \m_n + (1 - p_n) (\nabla \cL)(\Theta_{n - 1})\bigr] \dott \label{setting:UGD:Theta}
\end{gather}
\end{setting}
\end{tcolorbox}
\end{samepage}
\end{savenotes}

As in many other works in scientific literature which study gradient based optimization methods analytically (cf., \eg, \cite[Lemmas~2.13 and 2.14]{DereichDoJentzenPhilippeAdamSymmetry}, \cite[(10) in Theorem~2.5]{DereichJentzenAdamRates}, \cite[Lemmas~3.1 and 3.2]{JENTZEN2020101438}, and \cite[Assumption~2.2]{Polyak_accelerating}), we also impose in this work suitable assumptions on the step-sizes $\gamma_n \in (0, \infty)$, $n \in \N$, (such as summability and non-summability of powers of the step-sizes as well as conditions on the scaled increments $(\gamma_n)^{- 2} (\gamma_n - \gamma_{n + 1}) \in [0, \infty)$, $n \in \N$, of the step-sizes) which ensure that the step-sizes converge sufficiently quickly to zero but also that the step-sizes do not converge too quickly to zero (see, \eg, \cref{thm:univ_scheme_convergence_Loja.Theta,thm:univ_scheme_convergence_Loja.bound1} in \cref{thm:univ_scheme_convergence_Loja}, \cref{theorem:gen_convergence_Loja.bound_2} and below \cref{theorem:gen_convergence_Loja.bound_2} in \cref{theorem:gen_convergence_Loja}, and \cref{theorem:gen_convergence_Loja_no_delta.bound_2} and below \cref{theorem:gen_convergence_Loja_no_delta.bound_2} in \cref{theorem:gen_convergence_Loja_no_delta}).

We also point out that in some of later results for the \UGD\ method we will choose the step-size sequence $(\gamma_n)_{n \in \N}$ to coincide with the sequence $n^{- \lrexpo}$, $n \in \N$, for some $\lrexpo \in (\nicefrac{3}{4}, 1]$ (see above \cref{cor:gen_convergence_rate_Loja_l_rate.m_new} in the proof of \cref{cor:gen_convergence_rate_Loja_l_rate}).

However, even after having made the choice that for all $n \in \N$ it holds that $\gamma_n = n^{- \lrexpo}$, the framework is still general enough to handle optimization methods with general learning rates (that do not necessarily agree with the sequence $(n^{ - \lrexpo})_{n \in \N})$ by employing the matrix valued sequence $\bbA_n \in \R^{\fd \times \fd}$, $n \in \N$, in \cref{setting:UGD:Theta} in \cref{setting:UGD}. For details on this issue we refer to \cref{cor:gen_convergence_rate_Loja_gen_l_rate_Lipschitz_alpha.Theta} and the proof of \cref{cor:gen_convergence_rate_Loja_gen_l_rate_Lipschitz_alpha}.

\subsection{Weak convergence of general gradient method}
\label{subsec:weak_convergence_general_GD}

In this subsection we establish in \cref{prop:general_convergence_ohne_Loja_2} below a weak convergence result, convergence of the evaluation of the objective function $\cL \colon \R^{\fd} \to \R$, for a general abstract class of optimization processes $(\Theta_n)_{n \in \N_0}$ (see \cref{prop:general_convergence_ohne_Loja_2.boundedness} in \cref{prop:general_convergence_ohne_Loja_2} for details). In \cref{prop:general_convergence_ohne_Loja_2} we also establish for every $n \in \N$ an upper bound for the value of $\cL(\Theta_n)$ in terms of $\cL(\Theta_{n - 1})$. We employ \cref{prop:general_convergence_ohne_Loja_2} in our proof of \cref{prop:gen_convergence_ohne_Loja} (the main result of this section) which establishes under suitable assumptions weak convergence for the \UGD\ method (see \cref{setting:UGD}).

In our proof of \cref{prop:general_convergence_ohne_Loja_2} we employ the well-known Taylor approximation estimate for continuously differentiable functions with locally H\"{o}lder continuous derivatives in \cref{lemma:Hoelder_cont_der}. Only for completeness we include here a detailed proof for \cref{lemma:Hoelder_cont_der}.

\cfclear
\begin{savenotes}
\begin{samepage}
\begin{tcolorbox}[colback=white!95!gray,
                  colframe=black,
                  boxrule=0.5pt,
                  sharp corners,
                  enhanced,
                 ]
\begin{athm}{lemma}{lemma:Hoelder_cont_der}[\resname{Taylor remainder bound}]
Let $\fd \in \N$, $\delta \in (0, \infty)$ and let $\cL \in C^1(\R^{\fd}, \R)$ have a locally $\delta$-H\"{o}lder continuous derivative. Then for every $B \in \R$ there exists $D \in \R$ such that for all $u, w \in \{x \in \R^{\fd} \colon \norm{x} \le B\}$ it holds that
\begin{equation}
\textstyle \abs{\cL(u) - \cL(w) - \spro{u - w, (\nabla \cL)(w)}} \le D \norm{u - w}^{1 + \delta}.
\end{equation}
\end{athm}
\end{tcolorbox}
\end{samepage}
\end{savenotes}
\begin{aproof}
Throughout this proof let $B \in [0, \infty)$.
\startnewargseq
\argument{the assumption that $\cL$ has a locally $\delta$-H\"{o}lder continuous derivative}{that there exists $L \in \R$ which satisfies for all $u, \allowbreak w \in \allowbreak \{x \in \allowbreak \R^{\fd} \colon \allowbreak \norm{x} \le B\}$ that
\begin{equation}\llabel{eqn1}
\textstyle \norm{(\nabla \cL)(u) - (\nabla \cL)(w)} \le L \norm{u - w}^{\delta}\dott
\end{equation}
}
\argument{\lref{eqn1}; the triangle inequality; H\"{o}lder's inequality; the fact that for all $u, \allowbreak w \in \allowbreak \{x \in \allowbreak \R^{\fd} \colon \allowbreak \norm{x} \le B\}$, $t \in [0, 1]$ it holds that $\norm{w + t (u - w)} \le B$}{ for all $u, \allowbreak w \in \allowbreak \{x \in \allowbreak \R^{\fd} \colon \allowbreak \norm{x} \le B\}$ that
\begin{equation}\llabel{eqn2}
\begin{split}
& \textstyle \abs{\cL(u) - \cL(w) - \spro{u - w, (\nabla \cL)(w)}} = \Abs{[\cL(w + t (u - w))]_{t = 0}^{t = 1} - \spro{u - w, (\nabla \cL)(w)}} \\
& \textstyle = \Abs{\int_0^1 \spro{u - w, (\nabla \cL)(w + t (u - w))} \d t - \spro{u - w, (\nabla \cL)(w)}} \\
& \textstyle = \Abs{\int_0^1 \spro{u - w, (\nabla \cL)(w + t (u - w)) - (\nabla \cL)(w)} \d t} \\
& \textstyle \le \int_0^1 \norm{u - w} \norm{(\nabla \cL)(w + t (u - w)) - (\nabla \cL)(w)} \d t \le L \norm{u - w}^{1 + \delta} \int_0^1 t^{\delta} \d t \\
& \textstyle = L (1 + \delta)^{-1} \norm{u - w}^{1 + \delta}\dott
\end{split}
\end{equation}
}
\argument{\lref{eqn2}}{for all $u, w \in \{x \in \R^{\fd} \colon \norm{x} \le B\}$ that
\begin{equation}
\textstyle \abs{\cL(u) - \cL(w) - \spro{u - w, (\nabla \cL)(w)}} \le L (1 + \delta)^{-1} \norm{u - w}^{1 + \delta}\dott
\end{equation}}
\end{aproof}

\cfclear
\begin{savenotes}
\begin{samepage}
\begin{tcolorbox}[colback=white!95!gray,
                  colframe=black,
                  boxrule=0.5pt,
                  sharp corners,
                  enhanced,
                 ]
\begin{athm}{prop}{prop:general_convergence_ohne_Loja_1}
Let $\fd \in \N$, $\delta \in (0, \infty)$, let $\cL \in C^1(\R^{\fd}, \R)$ have a locally $\delta$-H\"{o}lder continuous derivative, let $\gamma \colon \N \to (0, \infty)$, $\sigma \colon \N \to [0, \infty)$, $\Theta \allowbreak \colon \allowbreak \N_0 \allowbreak \to \allowbreak \R^\fd$, $\Gamma \colon \N \to \R^{\fd}$, and $\mu \colon \allowbreak \N \allowbreak \to \allowbreak \R^{\fd}$ satisfy for all $n \in \N$ that
\begin{equation}\llabel{eqn:prop:general_convergence_ohne_Loja_1}
\textstyle \Theta_n = \Theta_{n - 1} - \gamma_n \bigl(\Gamma_n + \sigma_n \mu_n\bigr) \ifnocf,
\end{equation}
and assume $\#\{n \in \N \colon \spro{\Gamma_n, (\nabla \cL)(\Theta_{n - 1})} < 0\} < \infty$ and
\begin{equation}\label{eqn:prop:general_convergence_ohne_Loja_1:boundedness}
\textstyle \bigl[\sum_{n = 1}^{\infty} \gamma_n [\sigma_n + (\gamma_n)^{\delta}]\bigr] + \limsup_{n \to \infty} \bigl[\norm{\Theta_n} + \norm{\Gamma_n} + \norm{\mu_n}\bigr] < \infty\dott
\end{equation}
\cfout[.]Then there exist $\fC \in \R$, $\vartheta \in \R^{\fd}$ such that
\begin{enumerate}[label=(\roman*)]
\item
\label{item1:prop:general_convergence_ohne_Loja_1} it holds for all $n \in \N$ that
\begin{equation}\llabel{result}
\begin{split}
& \textstyle \cL(\Theta_n) - \cL(\Theta_{n - 1}) + \gamma_n \spro{\Gamma_n, (\nabla \cL)(\Theta_{n - 1})} \\
& \textstyle \le \fC \min\bigl\{\gamma_n [\sigma_n + (\gamma_n)^{\delta}], \gamma_n \sigma_n \norm{(\nabla \cL)(\Theta_{n - 1})} + (\gamma_n)^{1 + \delta} \norm{\Gamma_n + \sigma_n \mu_n}^{1 + \delta}\bigr\} 
\end{split}
\end{equation}
and
\item
\label{item2:prop:general_convergence_ohne_Loja_1} it holds that $\limsup_{n \to \infty} \abs{\cL(\Theta_n) - \cL(\vartheta)} = 0$.
\end{enumerate}
\end{athm}
\end{tcolorbox}
\end{samepage}
\end{savenotes}
\begin{aproof}
\startnewargseq
\argument{Lemma~\ref{lemma:Hoelder_cont_der} (applied with $\fd \with \fd$, $\delta \with \delta$, $\cL \with \cL$ in the notation of Lemma~\ref{lemma:Hoelder_cont_der}); \cref{eqn:prop:general_convergence_ohne_Loja_1:boundedness}; the assumption that $\cL$ has a locally $\delta$-H\"{o}lder continuous derivative}{that there exist $B \in \R$, $D \in [B, \infty)$, $N \in \N$ which satisfy for all $n \in \allowbreak \N \allowbreak \cap \allowbreak [N, \allowbreak \infty)$, $u, w \in \allowbreak \{x \in \allowbreak \R^{\fd} \colon \allowbreak \norm{x} \allowbreak \le B\}$ that
\begin{equation}\llabel{eqn1}
\textstyle 1 + \bigl[\sum_{j = 1}^{\infty} \gamma_j [\sigma_j + (\gamma_j)^{\delta}]\bigr] + \sup_{j \in \N} \bigl[\norm{\Theta_{j - 1}} + \norm{(\nabla \cL)(\Theta_{j - 1})} + \norm{\Gamma_j} + \norm{\mu_j}\bigr] \le B,
\end{equation}
\begin{equation}\llabel{eqn2}
\textstyle \spro{\Gamma_n, (\nabla \cL)(\Theta_{n - 1})} \ge 0, \qandq \abs{\cL(u) - \cL(w) - \spro{u - w, (\nabla \cL)(w)}} \le D \norm{u - w}^{1 + \delta}.
\end{equation}
}
Let $\fC \in \R$ and $(\bbL_{n})_{n \in \N_0} \subseteq \R$ satisfy for all $n \in \N_0$ that
\begin{gather}\llabel{fC}
\textstyle \fC \ge \max\bigl\{D^2 + 2^{\delta} D^{2 + 2 \delta}, \max_{j \in \{1, 2, \dots, N\}}[(\gamma_j)^{-1}[\sigma_j + (\gamma_j)^{\delta}]^{-1} (\cL(\Theta_j) - \cL(\Theta_{j - 1}))]\bigr\} \\
\llabel{bbL}
\textstyle \andqq \bbL_n = \cL(\Theta_n) + \fC \sum_{j = n + 1}^{\infty} \gamma_j [\sigma_j + (\gamma_j)^{\delta}].
\end{gather}
\startnewargseq
\argument{\lref{eqn1}; \lref{eqn2}}{for all $n \in \N$ that
\begin{equation}\llabel{L_Theta_bound}
\textstyle \abs{\cL(\Theta_n) - \cL(\Theta_{n - 1}) - \spro{\Theta_n - \Theta_{n - 1}, (\nabla \cL)(\Theta_{n - 1})}} \le D \norm{\Theta_n - \Theta_{n - 1}}^{1 + \delta}\dott
\end{equation}
}
\argument{the triangle inequality; H\"{o}lder's inequality; \lref{eqn:prop:general_convergence_ohne_Loja_1}; \lref{L_Theta_bound}; \lref{eqn1}}[verbs=s]{for all $n \in \N$ that
\begin{equation}\llabel{pre_res0}
\begin{split}
& \textstyle \cL(\Theta_n) - \cL(\Theta_{n - 1}) \le \spro{\Theta_n - \Theta_{n - 1}, (\nabla \cL)(\Theta_{n - 1})} + D \norm{\Theta_n - \Theta_{n - 1}}^{1 + \delta} \\
& \textstyle = - \gamma_n \spro{\Gamma_n + \sigma_n \mu_n, (\nabla \cL)(\Theta_{n - 1})} + D (\gamma_n)^{1 + \delta} \norm{\Gamma_n + \sigma_n \mu_n}^{1 + \delta} \\
& \textstyle = - \gamma_n \spro{\Gamma_n, (\nabla \cL)(\Theta_{n - 1})} - \gamma_n \sigma_n \spro{\mu_n, (\nabla \cL)(\Theta_{n - 1})} + D (\gamma_n)^{1 + \delta} \norm{\Gamma_n + \sigma_n \mu_n}^{1 + \delta} \\
& \textstyle \le - \gamma_n \spro{\Gamma_n, (\nabla \cL)(\Theta_{n - 1})} + \gamma_n \sigma_n \norm{\mu_n} \norm{(\nabla \cL)(\Theta_{n - 1})} + D (\gamma_n)^{1 + \delta} \norm{\Gamma_n + \sigma_n \mu_n}^{1 + \delta} \\
& \textstyle \le - \gamma_n \spro{\Gamma_n, (\nabla \cL)(\Theta_{n - 1})} + D \bigl[\gamma_n \sigma_n \norm{(\nabla \cL)(\Theta_{n - 1})} + (\gamma_n)^{1 + \delta} \norm{\Gamma_n + \sigma_n \mu_n}^{1 + \delta}\bigr].
\end{split}
\end{equation}
}
\argument{the triangle inequality; H\"{o}lder's inequality; \lref{eqn1}; \lref{eqn2}; \lref{fC}; \lref{pre_res0}}{that for all $n \in \N$ it holds that
\begin{equation}\llabel{pre_res}
\begin{split}
& \textstyle \cL(\Theta_n) - \cL(\Theta_{n - 1}) + \gamma_n \spro{\Gamma_n, (\nabla \cL)(\Theta_{n - 1})} \\
& \textstyle \le D \bigl[\gamma_n \sigma_n \norm{(\nabla \cL)(\Theta_{n - 1})} + (\gamma_n)^{1 + \delta} \norm{\Gamma_n + \sigma_n \mu_n}^{1 + \delta}\bigr] \\
& \textstyle \le D^2 \gamma_n \sigma_n + 2^{\delta} D (\gamma_n)^{1 + \delta} (\norm{\Gamma_n}^{1 + \delta} + (\sigma_n)^{1 + \delta} \norm{\mu_n}^{1 + \delta}) \\
& \textstyle \le D^2 \gamma_n \sigma_n + 2^{\delta} D^{2 + \delta} (\gamma_n)^{1 + \delta} (1 + (\sigma_n)^{1 + \delta}) \\
& \textstyle = \gamma_n \sigma_n (D^2 + 2^{\delta} D^{2 + \delta} (\gamma_n \sigma_n)^{\delta}) + 2^{\delta} D^{2 + \delta} (\gamma_n)^{1 + \delta} \\
& \textstyle \le \gamma_n \sigma_n (D^2 + 2^{\delta} D^{2 + 2 \delta}) + 2^{\delta} D^{2 + \delta} (\gamma_n)^{1 + \delta} \\
& \textstyle \le (D^2 + 2^{\delta} D^{2 + 2 \delta}) \gamma_n [\sigma_n + (\gamma_n)^{\delta}] \le \fC \gamma_n [\sigma_n + (\gamma_n)^{\delta}].
\end{split}
\end{equation}
}
\argument{\lref{pre_res0}; \lref{pre_res}}[verbs=e]{\cref{item1:prop:general_convergence_ohne_Loja_1}. \llabel{empty}}
\argument{\lref{eqn2}; \lref{fC}; \lref{pre_res}}{for all $n \in \N$ that
\begin{equation}\llabel{res1}
\textstyle \cL(\Theta_n) - \cL(\Theta_{n - 1}) \le \fC \gamma_n [\sigma_n + (\gamma_n)^{\delta}].
\end{equation}
}
\argument{\lref{eqn1}; \lref{eqn2}; \lref{bbL}; \lref{res1}}{for all $n \in \N$ that $\sup_{k \in \N} \abs{\bbL_{k - 1}} \allowbreak < \allowbreak \infty$ and
\begin{equation}\llabel{pre_res2}
\begin{split}
& \textstyle \bbL_n - \bbL_{n - 1} = \bigl(\cL(\Theta_n) + \fC \sum_{j = n + 1}^{\infty} \gamma_j [\sigma_j + (\gamma_j)^{\delta}]\bigr) - \bigl(\cL(\Theta_{n - 1}) + \fC \sum_{j = n}^{\infty} \gamma_j [\sigma_j + (\gamma_j)^{\delta}]\bigr) \\
& = \cL(\Theta_n) - \cL(\Theta_{n - 1}) - \fC \gamma_n [\sigma_n + (\gamma_n)^{\delta}] \le 0.
\end{split}
\end{equation}
}
\argument{\lref{pre_res2}; Monotone convergence theorem}{that there exists $\scrL \in \R$ which satisfy 
\begin{equation}\llabel{scrL}
\textstyle \limsup_{n \to \infty} \abs{\bbL_n - \scrL} = 0.
\end{equation}
}
\argument{\lref{eqn1}; the fact that $\cL$ is continuous}{that there exist an increasing sequence $(n_k)_{k \in \N} \subseteq \N$ and $\vartheta \in \R^{\fd}$ which satisfy
\begin{equation}\llabel{subseq}
\textstyle \limsup_{k \to \infty} \bigl[\abs{\cL(\Theta_{n_k}) - \cL(\vartheta)} + \sum_{j = k}^{\infty} \gamma_k [\sigma_k + (\gamma_k)^{\delta}]\bigr] = 0.
\end{equation}
}
\argument{\lref{bbL}; \lref{scrL}; \lref{subseq}}{that
\begin{equation}\llabel{res2}
\textstyle \scrL = \cL(\vartheta) \qqandqq \limsup_{n \to \infty} \abs{\cL(\Theta_n) - \cL(\vartheta)} = 0.
\end{equation}
}
\argument{\lref{res2}}[verbs=e]{\cref{item2:prop:general_convergence_ohne_Loja_1}. }
\end{aproof}

\cfclear
\begin{savenotes}
\begin{samepage}
\begin{tcolorbox}[colback=white!95!gray,
                  colframe=black,
                  boxrule=0.5pt,
                  sharp corners,
                  enhanced,
                 ]
\begin{athm}{prop}{prop:general_convergence_ohne_Loja_2}
Let $\fd \in \N$, $c, \kappa, \delta \in (0, \infty)$, let $\cL \in C^1(\R^{\fd}, \R)$ have a locally $\delta$-H\"{o}lder continuous derivative, let $\Theta \allowbreak \colon \allowbreak \N_0 \allowbreak \to \allowbreak \R^\fd$, $\Gamma \colon \N \to \R^{\fd}$, and $\mu \colon \allowbreak \N \allowbreak \to \allowbreak \R^{\fd}$ be bounded, let $\gamma \colon \N \to [0, \infty)$ and $\sigma \colon \N \to [0, \infty)$ satisfy for all $n \in \N$ that
\begin{equation}\llabel{boundedness}
\textstyle \Theta_n = \Theta_{n - 1} - \gamma_n (\Gamma_n + \sigma_n \mu_n) \qqandqq \sum_{k = 1}^{\infty} \gamma_k (\sigma_k + (\gamma_k)^{\delta}) < \infty = \sum_{k = 1}^{\infty} \gamma_k,
\end{equation}
and assume for all $n \in \N \cap [c, \infty)$ that $\norm{(\nabla \cL)(\Theta_n)}^{\kappa} \le c \spro{\Gamma_{n + 1}, (\nabla \cL)(\Theta_n)}$. Then there exist $\fc \in \R$, $\vartheta \in (\nabla \cL)^{-1} (\{0\})$ which satisfy for all $n \in \N$ that
\begin{equation}\llabel{result}
\textstyle \cL(\Theta_n) - \cL(\Theta_{n - 1}) \le \fc \gamma_n [\sigma_n + (\gamma_n)^{\delta}] \qqandqq \limsup_{k \to \infty} \abs{\cL(\Theta_k) - \cL(\vartheta)} = 0.
\end{equation}
\end{athm}
\end{tcolorbox}
\end{samepage}
\end{savenotes}
\begin{aproof}
\startnewargseq
\argument{the assumption that for all $n \in \N \cap [c, \infty)$ it holds that $c \spro{\Gamma_{n + 1}, (\nabla \cL)(\Theta_n)} \ge \norm{(\nabla \cL)(\Theta_n)}^{\kappa}$}{that there exist $D \in [c, \infty)$, $N \in \N_0$ which satisfy for all $n \in \N \cap [N + 1, \infty)$ that
\begin{gather}\llabel{N_and_C}
\textstyle \spro{\Gamma_n, (\nabla \cL)(\Theta_{n - 1})} \ge \frac{1}{D} \norm{(\nabla \cL)(\Theta_{n - 1})}^{\kappa}, \quad \#\{j \in \N \colon \spro{\Gamma_j, (\nabla \cL)(\Theta_{j - 1})} < 0\} \le N, \\
\llabel{C}
\textstyle \andqq D \ge \max_{j \in \{1, 2, \dots, N\}}[(\gamma_j)^{-1}[\sigma_j + (\gamma_j)^{\delta}]^{-1} (\cL(\Theta_j) - \cL(\Theta_{j - 1}))].
\end{gather}
}
\argument{\lref{C}; \cref{prop:general_convergence_ohne_Loja_1} (applied with $\fd \with \fd$, $\delta \with \delta$, $\cL \with \cL$, $\gamma \with \gamma$, $\sigma \with \sigma$, $\Theta \with \Theta$, $\Gamma \with \Gamma$, $\mu \with \mu$ in the notation of \cref{prop:general_convergence_ohne_Loja_1})}{that there exist $\fC \in [D, \infty)$, $\psi \in \R^{\fd}$ which satisfy for all $n \in \N$ that $\limsup_{k \to \infty} \abs{\cL(\Theta_k) - \cL(\psi)} = 0$ and
\begin{equation}\llabel{fC}
\textstyle \cL(\Theta_n) - \cL(\Theta_{n - 1}) + \gamma_n \spro{\Gamma_n, (\nabla \cL)(\Theta_{n - 1})} \le \fC \gamma_n [\sigma_n + (\gamma_n)^{\delta}].
\end{equation}
}
\argument{\lref{N_and_C}; \lref{C}; \lref{fC}}{for all $n \in \N$ that
\begin{equation}\llabel{res0}
\textstyle \cL(\Theta_n) - \cL(\Theta_{n - 1}) \le \fC \gamma_n [\sigma_n + (\gamma_n)^{\delta}] \qandq \limsup_{k \to \infty} \abs{\cL(\Theta_k) - \cL(\psi)} = 0.
\end{equation}
}
\argument{\lref{N_and_C}; \lref{fC}}{for all $n \in \N \cap [N + 1, \infty)$ that
\begin{equation}\llabel{C_and_fC}
\textstyle \cL(\Theta_n) - \cL(\Theta_{n - 1}) + \frac{\gamma_n}{D} \norm{(\nabla \cL)(\Theta_{n - 1})}^{\kappa} \le \fC \gamma_n [\sigma_n + (\gamma_n)^{\delta}].
\end{equation}
}
\argument{\lref{C_and_fC}}{for all $n \in \N \cap [N + 1, \infty)$ that
\begin{equation}\llabel{convergence0}
\begin{split}
& \textstyle \frac{1}{D} \sum_{j = N + 1}^{n} \gamma_j \norm{(\nabla \cL)(\Theta_{j - 1})}^{\kappa} \\
& \textstyle = \cL(\Theta_N) - \cL(\Theta_n) + \sum_{j = N + 1}^{n} \bigl[\cL(\Theta_j) - \cL(\Theta_{j - 1}) + \frac{1}{D} \gamma_j \norm{(\nabla \cL)(\Theta_{j - 1})}^{\kappa}\bigr] \\
& \textstyle \le \cL(\Theta_N) - \cL(\Theta_n) + \fC \sum_{j = N + 1}^{n} \gamma_j [\sigma_j + (\gamma_j)^{\delta}].
\end{split}
\end{equation}
}
\argument{\lref{boundedness}; \lref{res0}; \lref{convergence0}}{that
\begin{equation}\llabel{weak_convergence0}
\textstyle \sum_{j = N + 1}^{\infty} \gamma_j \norm{(\nabla \cL)(\Theta_{j - 1})}^{\kappa} < \infty = \sum_{j = N + 1}^{\infty} \gamma_j.
\end{equation}
}
\argument{\lref{weak_convergence0}}{that $\liminf_{n \to \infty} \norm{(\nabla \cL)(\Theta_{n - 1})} = 0$. }
\argument{\lref{boundedness}; \lref{res0}; \lref{weak_convergence0}; the fact that $\cL$ and $\nabla \cL$ are continous}[verbs=d]{that there exists $\vartheta \in (\nabla \cL)^{-1}(\{0\}) \cap \cL^{-1}(\{\cL(\psi)\})$ which satisfies $\limsup_{k \to \infty} \abs{\cL(\Theta_k) - \cL(\vartheta)} = 0$. }
\end{aproof}

\subsection{Convergence rates for auxiliary momentum processes}
\label{subsec:auxiliary_momentum_approximation}

In this subsection we establish in \cref{lemma:momentum_bound} below convergence with a rate of convergence of a suitable auxiliary momentum process $(\m_n)_{n \in \N_0}$ (see \cref{eqn:lemma:momentum_bound:1} in \cref{lemma:momentum_bound} for details) to evaluations of the gradient of the objective function (in \cref{lemma:momentum_bound} denoted by the function $\cG \colon \allowbreak \R^{\fd} \allowbreak \to \R^{\fd}$) at the underlying optimization process $(\Theta_n)_{n \in \N_0}$. We employ \cref{lemma:momentum_bound} in our proof of \cref{prop:gen_convergence_ohne_Loja} (the main result of this section) which establishes under suitable assumptions weak convergence for the \UGD\ method (see \cref{setting:UGD} above).

\cfclear
\begin{savenotes}
\begin{samepage}
\begin{tcolorbox}[colback=white!95!gray,
                  colframe=black,
                  boxrule=0.5pt,
                  sharp corners,
                  enhanced,
                 ]
\begin{athm}{lemma}{lemma:momentum}
Let $\fd \in \N$ and let $\alpha \colon \N \to \R$, $g \colon \N_0 \to \R^{\fd}$, $m \colon \allowbreak \N_0 \allowbreak \to \R^\fd$, and $\mu \colon \N_0 \to \R^{\fd}$ satisfy for all $n \in \N$ that
\begin{equation}\label{eqn:lemma:momentum:1}
\textstyle m_n = \alpha_n m_{n - 1} + (1 - \alpha_n) g_{n - 1} + \mu_n.
\end{equation}
Then it holds for all $n, \ell \in \N$ with $\ell \le n$ that
\begin{equation}\label{eqn:lemma:momentum:2}
\textstyle m_n - g_n = \bigl[\prod_{j = \ell}^n \alpha_j\bigr] [m_{\ell - 1} - g_{\ell - 1}] - \sum_{j = \ell}^n \bigl(\bigl[\prod_{i = j + 1}^n \alpha_i\bigr] [(g_j - g_{j - 1}) - \mu_j]\bigr).
\end{equation}
\end{athm}
\end{tcolorbox}
\end{samepage}
\end{savenotes}
\begin{aproof}
Throughout this proof let $\kappa \colon \N_0 \to \R^{\fd}$ satisfy for all $n \in \N_0$ that
\begin{equation}\label{eqn:lemma:momentum:kappa}
\textstyle \kappa_n = m_n - g_n.
\end{equation}
\argument{\cref{eqn:lemma:momentum:1}; \cref{eqn:lemma:momentum:kappa}}{that for all $n \in \N$ it holds that
\begin{equation}\llabel{kappa_eq}
\begin{split}
\textstyle \kappa_n & \textstyle = m_n - g_n = (\alpha_n m_{n - 1} + (1 - \alpha_n) g_{n - 1} + \mu_n) - g_n \\
& \textstyle = \alpha_n (m_{n - 1} - g_{n - 1}) - (g_n - g_{n - 1}) + \mu_n = \alpha_n \kappa_{n - 1} - (g_n - g_{n - 1}) + \mu_n.
\end{split}
\end{equation}
}
\argument{\lref{kappa_eq}}{for all $n, \ell \in \N$ with $\ell \le n$ that
\begin{equation}\llabel{pre_res}
\begin{split}
& \textstyle \kappa_n - \bigl[\prod_{j = \ell}^n \alpha_j\bigr] \kappa_{\ell - 1} = \sum_{j = \ell}^n \bigl(\bigl[\prod_{i = j + 1}^n \alpha_i\bigr] \kappa_j - \bigl[\prod_{i = j}^n \alpha_i\bigr] \kappa_{j - 1}\bigr) \\
& \textstyle = \sum_{j = \ell}^n \bigl(\bigl[\prod_{i = j + 1}^n \alpha_i\bigr] (\kappa_j - \alpha_j \kappa_{j - 1})\bigr) = - \sum_{j = \ell}^n \bigl(\bigl[\prod_{i = j + 1}^n \alpha_i\bigr] [(g_j - g_{j - 1}) - \mu_j]\bigr).
\end{split}
\end{equation}
}
\argument{\cref{eqn:lemma:momentum:kappa}; \lref{pre_res}}[verbs=e]{\cref{eqn:lemma:momentum:2}. }
\end{aproof}

\cfclear
\begin{savenotes}
\begin{samepage}
\begin{tcolorbox}[colback=white!95!gray,
                  colframe=black,
                  boxrule=0.5pt,
                  sharp corners,
                  enhanced,
                 ]
\begin{athm}{lemma}{lemma:bound_geometric_sum}
Let $\alpha \colon \N \to [0, \infty)$ and $\gamma \colon \N \to (0, \infty)$ satisfy
\begin{equation}\label{eqn:lemma:bound_geometric_sum:gamma}
\textstyle \limsup_{n \to \infty} [(\gamma_{n + 1})^{-1} \alpha_{n + 1} \gamma_n] < 1.
\end{equation}
Then
\begin{equation}\label{eqn:lemma:bound_geometric_sum:result}
\textstyle \sup_{n \in \N} \bigl((\gamma_n)^{- 1} \sum_{j = 1}^{n} \bigl(\bigl[\prod_{i = j + 1}^n \alpha_i\bigr]\gamma_j\bigr)\bigr) < \infty \dott
\end{equation}
\end{athm}
\end{tcolorbox}
\end{samepage}
\end{savenotes}
\begin{aproof}
Throughout this proof let $A \colon \N \to \R$ satisfy for all $n \in \N$ that
\begin{equation}\label{eqn:lemma:bound_geometric_sum:A_n}
\textstyle A_n = \sum_{j = 1}^{n} \bigl(\bigl[\prod_{i = j + 1}^n \alpha_i\bigr]\gamma_j\bigr).
\end{equation}
\argument{\cref{eqn:lemma:bound_geometric_sum:A_n}}{that for all $n \in \N$ it holds that
\begin{equation}\label{eqn:lemma:bound_geometric_sum:A_n_2}
\begin{split}
\textstyle A_{n + 1} & \textstyle = \sum_{j = 1}^{n + 1} \bigl(\bigl[\prod_{i = j + 1}^{n + 1} \alpha_i\bigr]\gamma_j\bigr) = \gamma_{n + 1} + \sum_{j = 1}^{n} \bigl(\bigl[\prod_{i = j + 1}^{n + 1} \alpha_i\bigr]\gamma_j\bigr) \\
& \textstyle = \gamma_{n + 1} + \alpha_{n + 1} \sum_{j = 1}^{n} \bigl(\bigl[\prod_{i = j + 1}^{n} \alpha_i\bigr]\gamma_j\bigr) = \gamma_{n + 1} + \alpha_{n + 1} A_n.
\end{split}
\end{equation}
}
\argument{\cref{eqn:lemma:bound_geometric_sum:gamma}}{that there exist $N \in \N$, $\delta \in (0, \infty)$ which satisfy
\begin{equation}\label{eqn:lemma:bound_geometric_sum:delta}
\textstyle \Forall n \in \N \cap [N, \infty) \colon (\gamma_{n + 1})^{-1} \alpha_{n + 1} \gamma_n \le 1 - \delta.
\end{equation}
}
Let $c \in \R$ satisfy
\begin{equation}\label{eqn:lemma:bound_geometric_sum:c}
\textstyle c \ge \max\bigl\{\delta^{-1}, \max_{n \in \{1, 2, \dots, N\}} [(\gamma_n)^{-1} A_n]\bigr\}.
\end{equation}
\argument{\cref{eqn:lemma:bound_geometric_sum:delta}; \cref{eqn:lemma:bound_geometric_sum:c}}{for all $n \in \N \cap [N, \infty)$ that
\begin{equation}\llabel{c2}
\begin{split}
& \textstyle \gamma_{n + 1} (1 - c) + \alpha_{n + 1} \gamma_n c = \gamma_{n + 1} + c(\alpha_{n + 1} \gamma_n - \gamma_{n + 1}) = \gamma_{n + 1} + c \gamma_{n + 1} \bigl(\frac{\alpha_{n + 1}\gamma_n}{\gamma_{n + 1}} - 1\bigr) \\
& \textstyle \le \gamma_{n + 1} - c \gamma_{n + 1} \delta = \gamma_{n + 1} \delta \bigl(\delta^{-1} - c\bigr) \le 0.
\end{split}
\end{equation}
}
\argument{\cref{eqn:lemma:bound_geometric_sum:A_n_2}; \lref{c2}}{that for all $n \in \N \allowbreak \cap \allowbreak [N, \allowbreak \infty)$ it holds that
\begin{equation}\llabel{A_monotonicity}
\begin{split}
& \textstyle A_{n + 1} - c \gamma_{n + 1} = (\gamma_{n + 1} + \alpha_{n + 1} A_n) - c \gamma_{n + 1} = \gamma_{n + 1} (1 - c) + \alpha_{n + 1} A_n \\
& \textstyle = (\gamma_{n + 1} (1 - c) + \alpha_{n + 1} \gamma_n c) + \alpha_{n + 1} (A_n - c \gamma_n) \le \alpha_{n + 1} (A_n - c \gamma_n).
\end{split}
\end{equation}
}
\argument{\lref{A_monotonicity}}{for all $n \in \N \cap [N, \infty)$ that
\begin{equation}\llabel{A_monotonicity2}
\begin{split}
& \textstyle A_n - c \gamma_n - \bigl[\prod_{j = N + 1}^n \alpha_j\bigr] (A_N - c \gamma_N) \\
& \textstyle = \sum_{j = N}^{n - 1} \bigl(\bigl[\prod_{i = j + 2}^n \alpha_i\bigr] (A_{j + 1} - c \gamma_{j + 1}) - \bigl[\prod_{i = j + 1}^n \alpha_i\bigr] (A_j - c \gamma_j)\bigr) \\
& \textstyle = \sum_{j = N}^{n - 1} \bigl(\bigl[\prod_{i = j + 2}^n \alpha_i\bigr] [A_{j + 1} - c \gamma_{j + 1} - \alpha_{j + 1} (A_j - c \gamma_j)]\bigr) \le 0.
\end{split}
\end{equation}
}
\argument{\cref{eqn:lemma:bound_geometric_sum:c}; \lref{A_monotonicity2}}{that for all $n \in \N \cap [N, \infty)$ it holds that
\begin{equation}\llabel{pre_res}
\textstyle A_n - c \gamma_n \le \bigl[\prod_{j = N + 1}^n \alpha_j\bigr] (A_N - c \gamma_N) \le 0.
\end{equation}
}
\argument{\cref{eqn:lemma:bound_geometric_sum:A_n}; \cref{eqn:lemma:bound_geometric_sum:c}; \lref{pre_res}}{that for all $n \in \N$ it holds that
\begin{equation}\llabel{res}
\textstyle  \sum_{j = 1}^{n} \bigl(\bigl[\prod_{i = j + 1}^n \alpha_i\bigr]\gamma_j\bigr) = A_n \le c \gamma_n.
\end{equation}
}
\argument{\lref{res}}[verbs=e]{\cref{eqn:lemma:bound_geometric_sum:result}. }
\end{aproof}

\cfclear
\begin{savenotes}
\begin{samepage}
\begin{tcolorbox}[colback=white!95!gray,
                  colframe=black,
                  boxrule=0.5pt,
                  sharp corners,
                  enhanced,
                 ]
\begin{athm}{prop}{lemma:momentum_bound}
Let $\fd \in \N$, $\delta \in (0, \infty)$, let $\cG \colon \allowbreak \R^{\fd} \allowbreak \to \R^{\fd}$ be locally $\delta$-H\"{o}lder continuous, let $\Theta \colon \allowbreak \N_0 \allowbreak \to \allowbreak \R^{\fd}$ and $\mu \colon \allowbreak \N_0 \allowbreak \to \allowbreak \R^{\fd}$ be bounded, let $\alpha \colon \N \to [0, \infty)$, $\gamma \colon \N \to (0, \infty)$, and $\m \colon \allowbreak \N_0 \allowbreak \to \allowbreak \R^\fd$ satisfy for all $n \in \N$ that
\begin{equation}\label{eqn:lemma:momentum_bound:1}
\textstyle \m_n = \alpha_n \m_{n - 1} + (1 - \alpha_n) \cG(\Theta_{n - 1}) + (\gamma_n)^{\delta} \mu_n,
\end{equation}
and assume for all $k \in \{0, 1\}$ that
\begin{equation}\llabel{lemma:momentum_bound:gamma}
\textstyle \limsup_{n \to \infty} \bigl[(\gamma_n)^{- 1} (\alpha_n)^{k / \delta} (\gamma_{n - 1})^k \norm{\Theta_n - \Theta_{n - 1}}^{1 - k}\bigr] < \sup_{\eps \in (0, 1)} [k + \eps]^{- 1} \dott
\end{equation}
Then
\begin{equation}\llabel{eqn:lemma:momentum_bound:result}
\textstyle \sup_{n \in \N} \bigl[(\gamma_n)^{- \delta} (\sum_{k = 0}^{1} \norm{\m_n - \cG(\Theta_{n - k})})\bigr] < \infty.
\end{equation}
\end{athm}
\end{tcolorbox}
\end{samepage}
\end{savenotes}
\begin{aproof}
Throughout this proof for every $n \in \N$ let $g_{n - 1} = \cG(\Theta_{n - 1})$.
\argument{\lref{lemma:momentum_bound:gamma}; the assumption that $\cG$ is locally $\delta$-H\"{o}lder continuous}{that there exist $L, B \in \R$ which satisfy for all $n \in \N$, $u, \allowbreak w \in \allowbreak \{x \in \allowbreak \R^{\fd} \colon \allowbreak \norm{x} \le B\}$ that
\begin{equation}\label{eqn:lemma:momentum_bound:bounds}
\textstyle \norm{\Theta_{n - 1}} + \norm{\mu_n} + (\gamma_n)^{-1} \norm{\Theta_n - \Theta_{n - 1}} \le B \qandq \norm{\cG(u) - \cG(w)} \le L \norm{u - w}^{\delta}.
\end{equation}
}
\argument{\lref{lemma:momentum_bound:gamma}; \cref{lemma:bound_geometric_sum} (applied for every $n \in \N$ with $\alpha_n \with \alpha_n$, $\gamma_n \with (\gamma_n)^{\delta}$ in the notation of \cref{lemma:bound_geometric_sum})}{that there exists $\fC \in \R$ which satisfies for all $n \in \N$ that
\begin{equation}\label{eqn:lemma:momentum_bound:fC}
\textstyle \sum_{j = 1}^{n} \bigl(\bigl[\prod_{i = j + 1}^n \alpha_i\bigr](\gamma_j)^{\delta}\bigr) \le \fC (\gamma_n)^{\delta}.
\end{equation}
}
\argument{\cref{eqn:lemma:momentum_bound:1}; \cref{lemma:momentum} (applied for every $n \in \N$ with $\fd \with \fd$, $\alpha_n \with \alpha_n$, $g_{n - 1} \with g_{n - 1}$, $m_n \with \m_n$, $\mu_n \with (\gamma_n)^{\delta} \mu_n$ in the notation of \cref{lemma:momentum})}[verbs=s]{for all $n \in \N$ that
\begin{equation}\llabel{eqn:lemma:momentum_bound:m}
\textstyle \m_n - g_n = \bigl[\prod_{j = 1}^n \alpha_j\bigr] [\m_0 - g_0] - \sum_{j = 1}^n \bigl(\bigl[\prod_{i = j + 1}^n \alpha_i\bigr] [(g_j - g_{j - 1}) - (\gamma_j)^{\delta} \mu_j]\bigr).
\end{equation}
}
\argument{\cref{eqn:lemma:momentum_bound:bounds}; \cref{eqn:lemma:momentum_bound:fC}; \lref{eqn:lemma:momentum_bound:m}; the triangle inequality}{for all $n \in \N$ that
\begin{equation}\llabel{pre_res}
\begin{split}
& \textstyle \norm{\m_n - g_n} \\
& \textstyle = \Norm{\bigl[\prod_{j = 1}^n \alpha_j\bigr] [\m_0 - g_0] - \sum_{j = 1}^n \bigl(\bigl[\prod_{i = j + 1}^n \alpha_i\bigr] [(g_j - g_{j - 1}) - (\gamma_j)^{\delta} \mu_j]\bigr)} \\
& \textstyle \le \bigl[\prod_{j = 1}^n \alpha_j\bigr] \norm{\m_0 - g_0} + \sum_{j = 1}^n \bigl(\bigl[\prod_{i = j + 1}^n \alpha_i\bigr] (\norm{g_j - g_{j - 1}} + (\gamma_j)^{\delta} \norm{\mu_j})\bigr) \\
& \textstyle \le \bigl[\prod_{j = 1}^n \alpha_j\bigr] \norm{\m_0 - g_0} + \sum_{j = 1}^n \bigl(\bigl[\prod_{i = j + 1}^n \alpha_i\bigr] (L \norm{\Theta_j - \Theta_{j - 1}}^{\delta} + (\gamma_n)^{\delta} \norm{\mu_j})\bigr) \\
& \textstyle \le \bigl[\prod_{j = 1}^n \alpha_j\bigr] \norm{\m_0 - g_0} + \sum_{j = 1}^n \bigl(\bigl[\prod_{i = j + 1}^n \alpha_i\bigr] (L B^{\delta} (\gamma_j)^{\delta} + (\gamma_n)^{\delta} B)\bigr) \\
& \textstyle = \bigl[\prod_{j = 1}^n \alpha_j\bigr] \norm{\m_0 - g_0} + (L B^{\delta} + B) \sum_{j = 1}^n \bigl(\bigl[\prod_{i = j + 1}^n \alpha_i\bigr] (\gamma_j)^{\delta}\bigr) \\
& \textstyle \le \bigl[\prod_{j = 1}^n \alpha_j\bigr] \norm{\m_0 - g_0} + (L B^{\delta} + B) \fC (\gamma_n)^{\delta}.
\end{split}
\end{equation}
}
\argument{\lref{lemma:momentum_bound:gamma}}{that there exists $N \in \N$ which satisfies for all $n \in \N \cap [N, \infty)$ that
\begin{equation}\llabel{telescopic_prod}
\textstyle \gamma_N (\gamma_n)^{-1} \bigl[\prod_{j = N + 1}^n (\alpha_j)^{1 / \delta}\bigr] = \prod_{j = N}^{n - 1} \bigl[(\gamma_{j + 1})^{-1} (\alpha_{j + 1})^{1 / \delta} \gamma_j\bigr] \le 1 \dott
\end{equation}
}
\argument{\lref{telescopic_prod}}{that
\begin{equation}\llabel{telescopic_prod_bound}
\textstyle \limsup_{n \to \infty} (\gamma_n)^{- \delta} \bigl[\prod_{j = 1}^n \alpha_j \bigr] < \infty \dott
\end{equation}
}
\argument{\lref{pre_res}; \lref{telescopic_prod_bound}}{that
\begin{equation}\llabel{res}
\textstyle \limsup_{n \to \infty} \bigl[(\gamma_n)^{- \delta} \norm{\m_n - g_n}\bigr] < \infty.
\end{equation}
}
\argument{\cref{eqn:lemma:momentum_bound:bounds}; \lref{res}; the triangle inequality}{that
\begin{equation}\llabel{res1}
\begin{split}
& \textstyle \sup\limits_{n \in \N} \bigl(\frac{1}{(\gamma_n)^{\delta}} \norm{\m_n - g_{n - 1}}\bigr) \le \sup\limits_{n \in \N} \frac{1}{(\gamma_n)^{\delta}} (\norm{\m_n - g_n} + \norm{g_n - g_{n - 1}}) \\
& \textstyle \le \sup\limits_{n \in \N} \frac{1}{(\gamma_n)^{\delta}} \bigl(\norm{\m_n - g_n} + L \norm{\Theta_n - \Theta_{n - 1}}^{\delta}\bigr) \le \sup\limits_{n \in \N} \bigl(\frac{1}{(\gamma_n)^{\delta}} \norm{\m_n - g_n} + LB^{\delta}\bigr) < \infty.
\end{split}
\end{equation}
}
\argument{\lref{res}; \lref{res1}}[verbs=e]{\lref{eqn:lemma:momentum_bound:result}. }
\end{aproof}

\subsection{Weak convergence of UGD}
\label{subsec:UGD_weak_convergence}

In this subsection we combine \cref{prop:general_convergence_ohne_Loja_1} and \cref{prop:general_convergence_ohne_Loja_2} from \cref{subsec:weak_convergence_general_GD} with \cref{lemma:momentum_bound} from \cref{subsec:auxiliary_momentum_approximation} to establish in \cref{prop:gen_convergence_ohne_Loja} weak convergence of the \UGD\ method from \cref{setting:UGD}. In our proof of \cref{prop:gen_convergence_ohne_Loja} we also employ the elementary upper bound for geometrically weighted averages (momentum-type sequences) in \cref{lemma:momentum_priori_bound_2}. Only for completeness we include here a detailed proof for \cref{lemma:momentum_priori_bound_2}.

\cfclear
\begin{savenotes}
\begin{samepage}
\begin{tcolorbox}[colback=white!95!gray,
                  colframe=black,
                  boxrule=0.5pt,
                  sharp corners,
                  enhanced,
                 ]
\begin{athm}{lemma}{lemma:momentum_priori_bound_2}
Let $\alpha \colon \N \to [0, 1]$, $A \colon \N_0 \to \R$, and $x \colon \N_0 \to \R$ satisfy for all $n \in \N$ that $x_n \le \alpha_n x_{n - 1} + (1 - \alpha_n) A_{n - 1}$. Then it holds for all $n, m \in \N_0$ with $m < n$ that
\begin{equation}\llabel{res}
\textstyle x_n \le \max\{x_n, \sup_{v \in \N_0} A_v\} \le \max\{x_m, \sup_{v \in \N_0} A_v\}.
\end{equation}
\end{athm}
\end{tcolorbox}
\end{samepage}
\end{savenotes}
\begin{aproof}
\argument{the fact that for all $n \in \N$ it holds that $0 \le \alpha_n \le 1$; the fact that for all $n \in \N$ it holds that $x_n \le \alpha_n x_{n - 1} + (1 - \alpha_n) A_{n - 1}$}[verbs=d]{that for all $n \in \N$ it holds that
\begin{equation}\llabel{induction}
\begin{split}
\textstyle x_n & \textstyle \le \alpha_n x_{n - 1} + (1 - \alpha_n) A_{n - 1} \le \alpha_n \max\{x_{n - 1}, A_{n - 1}\} + (1 - \alpha_n) \max\{x_{n - 1}, A_{n - 1}\} \\
& \textstyle = \max\{x_{n - 1}, A_{n - 1}\} \le \max\{x_{n - 1}, \sup_{v \in \N_0} A_v\}.
\end{split}
\end{equation}
}
\argument{\lref{induction}}{that for all $n \in \N$ it holds that
\begin{equation}\llabel{pre_res}
\begin{split}
\textstyle \max\{x_n, \sup_{v \in \N_0} A_v\} & \textstyle \le \max\{\max\{x_{n - 1}, \sup_{v \in \N_0} A_v\}, \sup_{v \in \N_0} A_v\} \\
& \textstyle = \max\{x_{n - 1}, \sup_{v \in \N_0} A_v\}.
\end{split}
\end{equation}
}
\argument{\lref{pre_res}}[verbs=e]{\lref{res}. }
\end{aproof}

\cfclear
\begin{savenotes}
\begin{samepage}
\begin{tcolorbox}[colback=white!95!gray,
                  colframe=black,
                  boxrule=0.5pt,
                  sharp corners,
                  enhanced,
                 ]
\begin{athm}{cor}{prop:gen_convergence_ohne_Loja}[\resname{Weak convergence for the \UGD\ method}]
Let $\fd \in \N$, $\delta \in (0, \infty)$, let $\cL \in C^1(\R^{\fd}, \R)$ have a locally $\delta$-H\"{o}lder continuous derivative, let $\bbA \colon \N \to \R^{\fd \times \fd}$, $\Theta \colon \allowbreak \N_0 \allowbreak \to \allowbreak \R^\fd$, $\mu \colon \N \to \R^{\fd}$, and $p \colon \N \allowbreak \to \R$ be bounded, let $\alpha \colon \N \to [0, 1]$, $\gamma \colon \N \to (0, \infty)$, and $\m \colon \allowbreak \N_0 \allowbreak \to \R^\fd$ satisfy for all $n \allowbreak \in \N$ that
\begin{gather}
\textstyle \m_n = \alpha_n \m_{n - 1} + (1 - \alpha_n) [(\nabla \cL)(\Theta_{n - 1}) + (\gamma_n)^{\delta} \mu_n] \llabel{m} \\
\textstyle \andqq \Theta_n = \Theta_{n - 1} - \gamma_n \bbA_n \bigl[p_n \m_n + (1 - p_n) (\nabla \cL)(\Theta_{n - 1})\bigr] \ifnocf, \llabel{Theta}
\end{gather}
assume for all $n \in \N$ that $\bbA_n - \bbI_{\fd}$ is symmetric positive semi-definite, and assume
\begin{equation}\llabel{boundedness}
\textstyle \sum_{n = 1}^{\infty} (\gamma_n)^{1 + \delta} < \infty \qqandqq \limsup_{n \to \infty} \bigl[(\gamma_{n + 1})^{-1} (\alpha_{n + 1})^{1 / \delta} \gamma_n\bigr] < 1\dott
\end{equation}
\cfout[.]Then there exist $\fC \in \R$, $\vartheta \in \R^{\fd}$ which satisfy for all $n \in \N$ that
\begin{gather}
\textstyle (\gamma_n)^{-1} \norm{\Theta_n - \Theta_{n - 1}} + (\gamma_n)^{- \delta} \norm{\bbA_n (\m_n - (\nabla \cL)(\Theta_{n - 1}))} \le \fC, \llabel{result_1} \\
\textstyle \bigl[1 - \mathbbm{1}_{\R}\bigl(\sum_{k = 1}^{\infty} \gamma_k\bigr)\bigr] \norm{(\nabla \cL)(\vartheta)} + \limsup_{k \to \infty} \abs{\cL(\Theta_k) - \cL(\vartheta)} = 0, \llabel{result} \\
\textstyle \andqq \cL(\Theta_n) \le \cL(\Theta_{n - 1}) \allowbreak + \fC (\gamma_n)^{1 + \delta} \dott \llabel{result_compl}
\end{gather}
\end{athm}
\end{tcolorbox}
\end{samepage}
\end{savenotes}
\begin{aproof}
Throughout this proof for every $n \in \N$ let
\begin{equation}\llabel{Gamma_L_fm}
\textstyle \Gamma_n = \bbA_n (\nabla \cL)(\Theta_{n - 1}), \qquad \fm_n = (\gamma_n)^{-\delta} \bbA_n p_n (\m_n - (\nabla \cL)(\Theta_{n - 1})),
\end{equation}
and $K_n = \{m \in \N \cap [n, \infty) \colon (\nabla \cL)(\Theta_{m - 1}) \neq 0\}$.
\argument{\lref{Theta}; \lref{Gamma_L_fm}}{that for all $n \in \N$ it holds that
\begin{equation}\llabel{Theta_new}
\textstyle \Theta_n = \Theta_{n - 1} - \gamma_n (\Gamma_n + (\gamma_n)^{\delta} \fm_n).
\end{equation}
}
\argument{\cref{lemma:momentum_priori_bound_2} (applied for every $n \in \N$ with $\alpha_n \with \alpha_n$, $x_{n - 1} \with \norm{\m_{n - 1}}$, $A_{n - 1} \with \norm{(\nabla \cL)(\Theta_{n - 1}) + (\gamma_n)^{\delta} \mu_n}$, in the notation of \cref{lemma:momentum_priori_bound_2}); \lref{m}; \lref{boundedness}; the fact that $\nabla \cL$ is $\delta$-H\"{o}lder continuous; the assumption that $\Theta$ and $\mu$ is bounded; the triangle inequality}{that $\sup_{n \in \N_0} \norm{(\nabla \cL)(\Theta_n)} < \infty$ and
\begin{equation}\llabel{boundedness_m}
\textstyle \sup_{n \in \N_0} \norm{\m_n} \le \max\{\norm{\m_0}, \sup_{n \in \N} \norm{(\nabla \cL)(\Theta_{n - 1}) + (\gamma_n)^{\delta} \mu_n}\} < \infty.
\end{equation}
}
\argument{\lref{Theta}; \lref{boundedness}; \lref{boundedness_m}; H\"{o}lder's inequality; the assumption that $\bbA$ and $p$ are bounded; the triangle inequality}[verbs=d]{that $\sup_{n \in \N} \gamma_n < \infty$ and
\begin{equation}\llabel{boundedness_gamma_Theta}
\textstyle \sup_{n \in \N} [(\gamma_n)^{-1} \norm{\Theta_n - \Theta_{n - 1}}] = \sup_{n \in \N} \norm{\bbA_n [p_n \m_n + (1 - p_n) (\nabla \cL)(\Theta_{n - 1})]} < \infty.
\end{equation}
}
\argument{\cref{lemma:momentum_bound} (applied for every $n \in \N$ with $\fd \with \fd$, $\delta \with \delta$, $\zeta \with \zeta$, $\cG \with \nabla \cL$, $\alpha_n \with \alpha_n$, $\gamma_n \with \gamma_n$, $\Theta_n \with \Theta_n$, $\m_n \with \m_n$, $\mu_n \with (1 - \alpha_n) \mu_n$ in the notation of \cref{lemma:momentum_bound}); \lref{boundedness_gamma_Theta}; \lref{boundedness}; \lref{Gamma_L_fm}; the assumption that $\bbA$ and $p$ are bounded; H\"{o}lder's inequality}{that
\begin{gather}
\textstyle \sup_{n \in \N} (\gamma_n)^{- \delta} \norm{\bbA_n (\m_n - (\nabla \cL)(\Theta_{n - 1}))} < \infty \\
\textstyle \andqq \sup_{n \in \N} \norm{\fm_n} < \infty \dott \llabel{fm_bound}
\end{gather}
}
\argument{\lref{boundedness_gamma_Theta}; \lref{fm_bound}}[verbs=e]{\lref{result_1}. }
\argument{\lref{Gamma_L_fm}; the fact that for every $n \in \N$ it holds that $\bbA_n - \bbI_{\fd}$ is symmetric positive semi-definite; the fact that $\nabla \cL$ is $\delta$-H\"{o}lder continuous; the assumption that $\bbA$ and $\Theta$ are bounded; H\"{o}lder's inequality}{that for all $n \in \allowbreak \N$ it holds that
\begin{equation}\llabel{Gamma_bound1}
\textstyle \spro{\Gamma_n, (\nabla \cL)(\Theta_{n - 1})} = \spro{\bbA_n (\nabla \cL)(\Theta_{n - 1}), (\nabla \cL)(\Theta_{n - 1})} \ge \norm{(\nabla \cL)(\Theta_{n - 1})}^2
\end{equation}
and
\begin{equation}\llabel{Gamma_bound2}
\textstyle \sup_{k \in \N} \norm{\Gamma_k} = \sup_{k \in \N} \norm{\bbA_k (\nabla \cL)(\Theta_{k - 1})} \le \sup_{k \in \N} (\norm{\bbA_k} \norm{(\nabla \cL)(\Theta_{k - 1})}) < \infty\dott
\end{equation}
}
\argument{\lref{Gamma_bound2}}{that
\begin{equation}\llabel{L_bound}
\textstyle \liminf_{n \to \infty} \Bigl[\inf_{m \in K_n} \frac{\spro{\Gamma_m, (\nabla \cL)(\Theta_{m - 1})}}{\norm{(\nabla \cL)(\Theta_{m - 1})}^{2}}\Bigr] > 0
\end{equation}
and $\#\{n \in \N \colon \spro{\Gamma_n, (\nabla \cL)(\Theta_{n - 1})} < 0\} < \infty$ \dott
}
\argument{\cref{prop:general_convergence_ohne_Loja_1} (applied for every $n \in \N$ with $\fd \with \fd$, $\delta \with \delta$, $\cL \with \cL$, $\gamma_n \with \gamma_n$, $\sigma_n \with (\gamma_n)^{\delta}$, $\Theta_n \with \Theta_n$, $\Gamma_n \with \Gamma_n$, $\mu_n \with \fm_n$ in the notation of \cref{prop:general_convergence_ohne_Loja_1}); \cref{prop:general_convergence_ohne_Loja_2} (applied for every $n \in \N$ with $\fd \with \fd$, $\delta \with \delta$, $\kappa \with 2$, $\cL \with \cL$, $\gamma_n \with \gamma_n$, $\sigma_n \with (\gamma_n)^{\delta}$, $\Theta_n \with \Theta_n$, $\Gamma_n \with \Gamma_n$, $\mu_n \with \fm_n$ in the notation of \cref{prop:general_convergence_ohne_Loja_2}); \lref{boundedness}; \lref{Theta_new}; \lref{fm_bound}; \lref{Gamma_bound2}; \lref{L_bound}}[verbs=e]{\lref{result} and \lref{result_compl}. }
\end{aproof}

\cfclear
\begin{savenotes}
\begin{samepage}
\begin{tcolorbox}[colback=white!95!gray,
                  colframe=black,
                  boxrule=0.5pt,
                  sharp corners,
                  enhanced,
                 ]
\begin{athm}{cor}{cor_intro_weak_convergence}[\resname{Weak convergence for the \UGD\ method}]
Let $\fd \in \N$, $\lrexpo \allowbreak \in (\nicefrac{1}{2}, \allowbreak 1]$, let $\cL \allowbreak \in C^1(\R^{\fd}, \allowbreak \R)$ have a locally Lipschitz continuous derivative, let $\alpha \colon \N \to [0, 1]$, $p \colon \N \to \R$, $\bbA \colon \N \to \R^{\fd \times \fd}$, $\Theta \allowbreak \colon \allowbreak \N_0 \allowbreak \to \allowbreak \R^\fd$, $\m \allowbreak \colon \allowbreak \N_0 \allowbreak \to \R^\fd$, and $\mu \colon \N \to \R^{\fd}$ satisfy for all $n \in \N$ that
\begin{gather}\llabel{m}
\textstyle \m_n = \alpha_n \m_{n - 1} + (1 - \alpha_n) \bigl[(\nabla \cL)(\Theta_{n - 1}) + \mu_n\bigr], \\
\llabel{Theta}
\textstyle \Theta_n = \Theta_{n - 1} - \bbA_n \bigl[p_n \m_n + (1 - p_n) (\nabla \cL)(\Theta_{n - 1})\bigr],
\end{gather}
and $\limsup_{k \to \infty} \alpha_k < 1 \le 1 + \sup_{k \in \N}(\norm{\Theta_k} + \abs{p_k} + k^{\lrexpo}(\norm{\bbA_k} + \norm{\mu_k})) < \infty$, and assume for all $n \in \N$ that $n^{\lrexpo} \bbA_n - \bbI_{\fd}$ is symmetric positive semi-definite \cfout. Then there exists $\vartheta \allowbreak \in \allowbreak \R^{\fd}$ such that
\begin{equation}\llabel{result}
\textstyle (\nabla \cL) (\vartheta) = 0 \qqandqq \limsup_{n \to \infty} \abs{\cL(\Theta_n) - \cL(\vartheta)} = 0\dott
\end{equation}
\end{athm}
\end{tcolorbox}
\end{samepage}
\end{savenotes}
\begin{aproof}
\argument{\cref{prop:gen_convergence_ohne_Loja} (applied for every $n \in \N$ with $\fd \with \fd$, $\delta \with \allowbreak 1$, $\cL \with \cL$, $\bbA_n \with n^{\lrexpo} \bbA_n$, $\Theta_{n - 1} \with \Theta_{n - 1}$, $\mu_n \with n^{\lrexpo} \mu_n$, $p_n \with p_n$, $\alpha_n \with \alpha_n$, $\gamma_n \with n^{- \lrexpo}$, $\m_{n - 1} \with \m_{n - 1}$ in the notation of \cref{prop:gen_convergence_ohne_Loja})}[verbs=e]{\lref{result}. }
\end{aproof}

\cfclear
\begin{savenotes}
\begin{samepage}
\begin{tcolorbox}[colback=white!95!gray,
                  colframe=black,
                  boxrule=0.5pt,
                  sharp corners,
                  enhanced,
                 ]
\begin{athm}{cor}{cor:gen_convergence_ohne_Loja}[\resname{Weak convergence for the standard \GD\ method}]
Let $\fd \in \N$, $\delta \in (0, \allowbreak \infty)$, let $\cL \allowbreak \in C^1(\R^{\fd}, \allowbreak \R)$ have a locally $\delta$-H\"{o}lder continuous derivative, let $\gamma \colon \N \to [0, \infty)$ and $\Theta \allowbreak \colon \allowbreak \N_0 \allowbreak \to \allowbreak \R^\fd$ satisfy for all $n \in \N$ that
\begin{equation}\llabel{Theta}
\textstyle \Theta_n = \Theta_{n - 1} - \gamma_n (\nabla \cL)(\Theta_{n - 1}) \ifnocf,
\end{equation}
and assume
\begin{equation}\llabel{boundedness}
\textstyle \bigl[\sum_{n = 1}^{\infty} (\gamma_n)^{1 + \delta}\bigr] + \limsup_{n \to \infty} \norm{\Theta_n} < \infty\ifnocf.
\end{equation}
\cfout[.]Then there exists $\vartheta \in \R^{\fd}$ such that
\begin{equation}\llabel{result}
\textstyle \bigl[1 - \mathbbm{1}_{\R}\bigl(\sum_{n = 1}^{\infty} \gamma_n\bigr)\bigr] \norm{(\nabla \cL)(\vartheta)} + \limsup_{n \to \infty} \abs{\cL(\Theta_n) - \cL(\vartheta)} = 0.
\end{equation}
\end{athm}
\end{tcolorbox}
\end{samepage}
\end{savenotes}
\begin{aproof}
\argument{\cref{prop:gen_convergence_ohne_Loja} (applied for every $n \in \N$ with $\fd \with \fd$, $\delta \with \delta$, $\cL \with \cL$, $\bbA_n \with \bbI_{\fd}$, $\Theta_n \with \Theta_n$, $p_n \with 0$, $\gamma_n \with \gamma_n$ in the notation of \cref{prop:gen_convergence_ohne_Loja}); \lref{Theta}; \lref{boundedness}}[verbs=e]{\lref{result}. }
\end{aproof}

We have established \cref{cor:gen_convergence_ohne_Loja} for the standard \GD\ method here through a direct application of the more general weak convergence result for the \UGD\ method in \cref{prop:gen_convergence_ohne_Loja} above. The conclusion of \cref{cor:gen_convergence_ohne_Loja} for the standard \GD\ method is, of course, well-known in the scientific literature. For example, Theorem~1.3 in \cite{DereichKassing} directly implies \cref{cor:gen_convergence_ohne_Loja}.

\section{Strong convergence for GD optimization methods}
\label{section:UGD_strong_convergence}

In this section we establish in \cref{cor:gen_convergence_rate_Loja_l_rate} below the main result of this work, in which we develop a unified framework which includes each of the optimization methods \ref{intro:momentum}--\ref{intro:Yogi} from \cref{sec:introduction} as a special case and under which every bounded trajectory of the optimizer converges with convergence rates to a critical point for general \KL\ objective functions with locally Lipschitz continuous gradients.

In \cref{cor:gen_convergence_rate_Loja_gen_l_rate_Lipschitz_alpha} we reformulate and specialize \cref{cor:gen_convergence_rate_Loja_l_rate} to the situation where the momentum decay factors do not depend on the number of time steps and in \cref{section:UGD_applications} we apply \cref{cor:gen_convergence_rate_Loja_gen_l_rate_Lipschitz_alpha} to each of the optimization methods \ref{intro:momentum}--\ref{intro:Yogi} from \cref{sec:introduction}. The arguments in our proof of \cref{thm:univ_scheme_convergence_Loja} are in parts inspired by the arguments, \eg, in the proofs of \cite[Theorem~3.2]{MR2197994} and \cite[Proposition~8.1]{MR4716376}.

\subsection{Kurdyka-{\L}ojasiewicz (KL) functions}
\label{subsec:KL_functions}

In the following notion, \cref{def:Loja}, we recall the notion of a \KL\ function (cf., \eg, \cite[Definition~1]{MR3785672}, \cite[Definition~9]{pmlr-v129-barakat20a}, and \cite[Definition~9.1.2]{JentzenBookDeepLearning2023}).

\begin{definition}[\resname{\KL\ functions}]\label{def:Loja}
Let $\fd \in \N$ and let $\cL \colon \R^{\fd} \to \R$ be a function. Then we say that $\cL$ is a \KL\ function if and only if we have that
\begin{enumerate}[label=(\roman*)]
\item it holds that $\cL \in C^1(\R^{\fd}, \R)$ and
\item it holds for all $\theta \in \R^{\fd}$ that there exist $\eps, \scrc \in (0, \infty)$, $\alpha \in (0, 1)$ such that for all $\vartheta \in \{v \in \R^{\fd} \colon \norm{v - \theta} < \eps\}$ it holds that
\begin{equation}
\textstyle \abs{\cL(\vartheta) - \cL(\theta)}^{\alpha} \le \scrc \norm{(\nabla \cL)(\vartheta)}.
\end{equation}
\end{enumerate}
\end{definition}

Regarding \cref{def:Loja} we briefly recall the fundamental result of {\L}ojasiewicz~\cite[Proposition~1]{Lojasiewicz} that for every $\fd \in \N$ we have that every analytic function from $\R^{\fd}$ to $\R$ is also a \KL\ function (cf., \eg, also \cite[Proposition~6.8]{PMIHES_1988_67_5_0}, \cite[Theorem~3.1]{MR2274510}, and \cite[Corollary~9.10.2]{JentzenBookDeepLearning2023}).

\subsection{UGD applied to KL functions with locally Hölder continuous gradients}

In this subsection we first establish in \cref{thm:univ_scheme_convergence_Loja} strong convergence, different type of error estimates, and certain decent-type properties for a general class of abstract optimization processes (including the \UGD\ framework from \cref{setting:UGD} and the optimization methods \ref{intro:momentum}--\ref{intro:Yogi} from \cref{sec:introduction}, respectively, as special cases). 

Thereafter, we combine \cref{thm:univ_scheme_convergence_Loja} with the elementary upper bound for the partial sum of the step-sizes $(\gamma_n)_{n \in \N} \subseteq (0, \infty)$ in \cref{lr_for_rate} and the elementary upper bound for the real recursion in \cref{rate_lemma} to establish in \cref{theorem:gen_convergence_Loja} strong convergence with a rate of convergence for the \UGD\ method from \cref{setting:UGD}.

The general class of abstract optimization processes in \cref{thm:univ_scheme_convergence_Loja.Theta} and \cref{thm:univ_scheme_convergence_Loja.bound2} in \cref{thm:univ_scheme_convergence_Loja} can be regarded as a slightly modified variant of the deterministic version of the class of adaptive algorithms in \cite[(1.1.1) in Chapter~1]{MR1082341}.

\cfclear
\begin{savenotes}
\begin{samepage}
\begin{tcolorbox}[colback=white!95!gray,
                  colframe=black,
                  boxrule=0.5pt,
                  sharp corners,
                  enhanced,
                 ]
\begin{athm}{prop}{thm:univ_scheme_convergence_Loja}
Let $\fd \in \N$, $\delta \in (0, \infty)$, $\tau \in (0, 1]$, $\kappa \in [\frac{2 + \delta}{1 + \delta}, \infty)$, let $\cL \colon \R^{\fd} \to \R$ be a\cfadd{def:Loja} \KL\ function, assume that $\nabla \cL$ is locally $\delta$-H\"{o}lder continuous, let $\Theta \allowbreak \colon \allowbreak \N_0 \allowbreak \to \allowbreak \R^\fd$, $\Gamma \colon \N \to \R^{\fd}$, and $\mu \colon \allowbreak \N \allowbreak \to \allowbreak \R^{\fd}$ be bounded, let $\gamma \colon \allowbreak \N \allowbreak \to \allowbreak (0, \allowbreak \infty)$ and $\sigma \colon \N \allowbreak \to \allowbreak [0, \allowbreak \infty)$ satisfy for all $m \in \N$ that
\begin{gather}\llabel{Theta}
\textstyle \Theta_m = \Theta_{m - 1} - \gamma_m (\Gamma_m + \sigma_m \mu_m), \qquad \sum_{n = 1}^{\infty} \gamma_n [\sigma_n + (\gamma_n)^{\delta}]^{\tau} < \infty, \\
\textstyle \andqq \liminf_{n \to \infty} \Bigl[\frac{\gamma_n [\sigma_n + (\gamma_n)^{\delta}]^{\frac{\kappa \tau}{\kappa - 1}} - \gamma_{n + 1} [\sigma_{n + 1} + (\gamma_{n + 1})^{\delta}]^{\frac{\kappa \tau}{\kappa - 1}}}{\gamma_n [\sigma_n + (\gamma_n)^{\delta}]^\frac{\kappa}{\kappa - 1}}\Bigr] > 0 = \limsup_{n \to \infty} \sigma_n, \llabel{bound1}
\end{gather}
and let $c \in (0, \infty)$ satisfy for all $n \in \N \cap [c, \infty)$ that
\begin{equation}\llabel{bound2}
\textstyle \norm{(\nabla \cL)(\Theta_n)}^{\kappa} \le c \spro{\Gamma_{n + 1}, (\nabla \cL)(\Theta_n)} \qqandqq \norm{\Gamma_{n + 1}} \le c \norm{(\nabla \cL)(\Theta_n)}^{\kappa - 1} \dott
\end{equation}
\cfout[.]Then there exist $\scrC, \eps \in (0, \infty)$, $\rho \in (0, 1)$, $\vartheta \in \R^{\fd}$ which satisfy for all $n \in \N$, $\theta \in \R^{\fd}$ with $\norm{\theta - \vartheta} < \eps$ that
\begin{gather}\llabel{Lojaaa}
\textstyle \abs{\cL(\theta) - \cL(\vartheta)}^{\rho} \le \scrC \norm{(\nabla \cL)(\theta)}, \\
\llabel{descending}
\textstyle \cL(\Theta_n) + \scrC^{-1} \gamma_n \norm{(\nabla \cL)(\Theta_{n - 1})}^{\kappa} \le \cL(\Theta_{n - 1}) + \scrC \gamma_n [\sigma_n + (\gamma_n)^{\delta}]^{\frac{\kappa}{\kappa - 1}}, \\
\llabel{rate_prep}
\textstyle \scrC^{-1} (\cL(\vartheta) - \cL(\Theta_{n - 1})) \le \sum_{j = n}^{\infty} \gamma_j [\sigma_j + (\gamma_j)^{\delta}]^{\frac{\kappa}{\kappa - 1}} \le \scrC \gamma_n [\sigma_n + (\gamma_n)^{\delta}]^{\frac{\kappa \tau}{\kappa - 1}}, \\
\llabel{rate}
\textstyle \frac{1}{\scrC} \norm{\Theta_{n - 1} - \vartheta}\! \le \!\Bigl[\cL(\Theta_{n - 1}) - \cL(\vartheta) + \scrC \!\sum\limits_{j = n}^{\infty}\! \gamma_j [\sigma_j + (\gamma_j)^{\delta}]^{\frac{\kappa}{\kappa - 1}}\Bigr]^{1 - \rho}\! + \sum\limits_{j = n}^{\infty}\! \gamma_j [\sigma_j + (\gamma_j)^{\delta}]^{\tau}, \\
\llabel{result}
\textstyle \andqq \bigl[1 - \mathbbm{1}_{\R}\bigl(\sum_{k = 1}^{\infty} \gamma_k\bigr)\bigr] \norm{(\nabla \cL)(\vartheta)} + \limsup_{k \to \infty} \norm{\Theta_k - \vartheta} = 0 \dott
\end{gather}
\end{athm}
\end{tcolorbox}
\end{samepage}
\end{savenotes}
\begin{aproof}
Throughout this proof let $\cG \colon \R^{\fd} \to \R^{\fd}$ satisfy $\cG = \nabla \cL$.
\argument{\lref{Theta}; \lref{bound1}; \lref{bound2}}{that 
\begin{equation}\llabel{spro_neg}
\textstyle \limsup_{n \to \infty} \gamma_n = 0 = \limsup_{n \to \infty} \sigma_n \qqandqq \#\{n \in \N \colon \spro{\Gamma_n, \cG(\Theta_{n - 1})} < 0\} < \infty.
\end{equation}
}
\argument{\lref{spro_neg}}{that
\begin{equation}\llabel{gamma_sigma_lim}
\textstyle \limsup_{n \to \infty} [\sigma_n + (\gamma_n)^{\delta}] = 0.
\end{equation}
}
\argument{\lref{Theta}; \lref{gamma_sigma_lim}; the fact that $0 < \tau \le 1$}{that
\begin{equation}\llabel{gamma_sigma_sum}
\textstyle \bigl(\sum_{n = 1}^{\infty} \gamma_n [\sigma_n + (\gamma_n)^{\delta}]^{\frac{\kappa}{\kappa - 1}}\bigr) + \bigl(\sum_{n = 1}^{\infty} \gamma_n [\sigma_n + (\gamma_n)^{\delta}]\bigr) < \infty.
\end{equation}
}
\argument{\cref{prop:general_convergence_ohne_Loja_1} (applied with $\fd \with \fd$, $\delta \with \delta$, $\cL \with \cL$, $\gamma \with \gamma$, $\sigma \with \sigma$, $\Theta \with \Theta$, $\Gamma \with \Gamma$, $\mu \with \mu$ in the notation of \cref{prop:general_convergence_ohne_Loja_1}); \lref{Theta}; \lref{spro_neg}; \lref{gamma_sigma_sum}}{that there exist $\fC \in (0, \infty)$, $\psi \in \R^{\fd}$ which satisfy for all $n \in \N$ that
\begin{equation}\llabel{key_ineq}
\textstyle \cL(\Theta_n) - \cL(\Theta_{n - 1}) + \gamma_n \spro{\Gamma_n, \cG(\Theta_{n - 1})} \le \fC \gamma_n \bigl[\sigma_n \norm{\cG(\Theta_{n - 1})} + (\gamma_n)^{\delta} \norm{\Gamma_n + \sigma_n \mu_n}^{1 + \delta}\bigr],
\end{equation}
\begin{equation}\llabel{limit_cL}
\textstyle \andqq \limsup_{k \to \infty} \abs{\cL(\Theta_k) - \cL(\psi)} = 0\dott
\end{equation}
}
\argument{the assumption that $\Theta$ is bounded; \lref{limit_cL}; the fact that $\cL$ is differentiable}{that there exist an increasing sequence $(n_k)_{k \in \N} \subseteq \N$ and $\vartheta \in \allowbreak \R^{\fd}$  which satisfy
\begin{equation}\llabel{vartheta}
\textstyle \cL(\vartheta) = \cL(\psi) \qqandqq \limsup_{k \to \infty} \norm{\Theta_{n_k} - \vartheta} = 0.
\end{equation}
}
\argument{\lref{bound1}; \lref{bound2}; \lref{gamma_sigma_lim}; the assumption that $\cL$ is a \KL\ function; the assumption that $\cG$ is H\"{o}lder continuous; the fact that $\kappa \ge \frac{2 + \delta}{1 + \delta}$}{that there exist $N \in \N$, $B, C, D, \eps \in (0, \infty)$, $\rho \in (0, 1)$ which satisfy for all $n \in \allowbreak \N \allowbreak \cap [N, \allowbreak \infty)$, $\theta \in \allowbreak \{u \in \allowbreak \R^{\fd} \colon \allowbreak \norm{u - \vartheta} < \eps\}$ that
\begin{gather}\llabel{Loja}
\textstyle \max\bigl\{1, \fC, \norm{\cG(\Theta_{n - N})}^{\kappa - 1 - \frac{1}{1 + \delta}}, \norm{\mu_{n - N + 1}}\bigr\} \le C, \qquad \abs{\cL(\theta) - \cL(\vartheta)}^{\rho} \le C \norm{\cG(\theta)}, \\
\llabel{N_C_1}
\textstyle \norm{\cG(\Theta_{n - 1})}^{\kappa} \le C \spro{\Gamma_n, \cG(\Theta_{n - 1})}, \qquad \norm{\Gamma_n} \le C \norm{\cG(\Theta_{n - 1})}^{\kappa - 1}, \\
\llabel{N_C_2}
\textstyle \gamma_n [\sigma_n + (\gamma_n)^{\delta}]^{\frac{\kappa \tau}{\kappa - 1}} - \gamma_{n + 1} [\sigma_{n + 1} + (\gamma_{n + 1})^{\delta}]^{\frac{\kappa \tau}{\kappa - 1}} \ge C^{-1} \gamma_n [\sigma_n + (\gamma_n)^{\delta}]^\frac{\kappa}{\kappa - 1}, \\
\llabel{B}
\textstyle 2^{3 + \delta} C^{4 + 2 \delta} \le B, \quad [\sigma_n + (\gamma_n)^{\delta}]^{\tau} \le B^{-1}, \qandq \max\{3 B C, 2^{\rho + 1} C^3 \frac{1}{1 - \rho}\} \le D\dott
\end{gather}
}
Let $\scrN_1, \scrN_2, \scrN_3 \subseteq \N \cap [N, \infty)$ satisfy
\begin{gather}\llabel{scrN_1}
\textstyle \Forall n \in \scrN_1 \colon \norm{\cG(\Theta_{n - 1})}^{\kappa - 1} \le B [\sigma_n + (\gamma_n)^{\delta}], \\
\llabel{scrN_2}
\textstyle \Forall n \in \scrN_2 \colon B [\sigma_n + (\gamma_n)^{\delta}] < \norm{\cG(\Theta_{n - 1})}^{\kappa - 1} \le B [\sigma_n + (\gamma_n)^{\delta}]^{\tau}, \\
\llabel{scrN_3}
\textstyle \Forall n \in \scrN_3 \colon B [\sigma_n + (\gamma_n)^{\delta}]^{\tau} < \norm{\cG(\Theta_{n - 1})}^{\kappa - 1},
\end{gather}
and $\scrN_1 \cup \scrN_2 \cup \scrN_3 = \N \cap [N, \infty)$.
\startnewargseq
\argument{\lref{scrN_1}; \lref{scrN_2}; \lref{scrN_3}}{that
\begin{equation}\llabel{scrNs_cap}
\textstyle (\scrN_1 \cap \scrN_2) \cup (\scrN_2 \cap \scrN_3) \cup (\scrN_3 \cap \scrN_1) = \emptyset.
\end{equation}
}
\argument{\lref{N_C_2}}{for all $n, m \in \allowbreak \N \allowbreak \cap [N, \allowbreak \infty)$ with $n \le m$ that
\begin{equation}\llabel{telescopic_sum}
\begin{split}
& \textstyle \gamma_n [\sigma_n + (\gamma_n)^{\delta}]^{\frac{\kappa \tau}{\kappa - 1}} - \gamma_m [\sigma_m + (\gamma_m)^{\delta}]^{\frac{\kappa \tau}{\kappa - 1}} \\
& \textstyle = \sum_{j = n}^{m - 1} \bigl(\gamma_j [\sigma_j + (\gamma_j)^{\delta}]^{\frac{\kappa \tau}{\kappa - 1}} - \gamma_{j + 1} [\sigma_{j + 1} + (\gamma_{j + 1})^{\delta}]^{\frac{\kappa \tau}{\kappa - 1}}\bigr) \\
& \textstyle \ge C^{-1} \sum_{j = n}^{m - 1} \gamma_j [\sigma_j + (\gamma_j)^{\delta}]^\frac{\kappa}{\kappa - 1}\dott
\end{split}
\end{equation}
}
\argument{\lref{gamma_sigma_lim}; \lref{telescopic_sum}}{that for all $n \in \N \cap [N, \infty)$ it holds that
\begin{equation}\llabel{Cesaro_ineq}
\begin{split}
& \textstyle \gamma_n [\sigma_n + (\gamma_n)^{\delta}]^{\frac{\kappa \tau}{\kappa - 1}} = \lim_{m \to \infty} \bigl(\gamma_n [\sigma_n + (\gamma_n)^{\delta}]^{\frac{\kappa \tau}{\kappa - 1}} - \gamma_m [\sigma_m + (\gamma_m)^{\delta}]^{\frac{\kappa \tau}{\kappa - 1}}\bigr) \\
& \textstyle \ge C^{-1} \lim_{m \to \infty} \Bigl(\sum_{j = n}^{m - 1} \gamma_j [\sigma_j + (\gamma_j)^{\delta}]^\frac{\kappa}{\kappa - 1}\Bigr) = C^{-1} \sum_{j = n}^{\infty} \gamma_j [\sigma_j + (\gamma_j)^{\delta}]^\frac{\kappa}{\kappa - 1}\dott
\end{split}
\end{equation}
}
\argument{\lref{Loja}; \lref{N_C_1}; the triangle inequality}{that for all $n \in \allowbreak \N \cap \allowbreak [N, \allowbreak \infty)$ it holds that
\begin{equation}\llabel{temp1}
\textstyle \norm{\Gamma_n + \sigma_n \mu_n}^{1 + \delta} \le (\norm{\Gamma_n} + \sigma_n \norm{\mu_n})^{1 + \delta} \le (C \norm{\cG(\Theta_{n - 1})}^{\kappa - 1} + C \sigma_n)^{1 + \delta}.
\end{equation}
}
\argument{\lref{Loja}; \lref{temp1}}{that for all $n \in \allowbreak \N \cap \allowbreak [N, \allowbreak \infty)$ it holds that
\begin{equation}\llabel{step_to_bbL_lower_bound_1}
\begin{split}
& \textstyle \sigma_n \norm{\cG(\Theta_{n - 1})} + (\gamma_n)^{\delta} \norm{\Gamma_n + \sigma_n \mu_n}^{1 + \delta} \\
& \textstyle \le \sigma_n \norm{\cG(\Theta_{n - 1})} + C^{1 + \delta} (\gamma_n)^{\delta} (\norm{\cG(\Theta_{n - 1})}^{\kappa - 1} + \sigma_n)^{1 + \delta} \\
& \textstyle \le C^{1 + \delta} \bigl[\sigma_n \norm{\cG(\Theta_{n - 1})} + (\gamma_n)^{\delta} (\norm{\cG(\Theta_{n - 1})}^{\kappa - 1} + \sigma_n)^{1 + \delta}\bigr].
\end{split} 
\end{equation}
}
\argument{\lref{key_ineq}; \lref{N_C_1}; \lref{step_to_bbL_lower_bound_1}}{for all $n \in  \N \cap [N, \infty)$ that
\begin{equation}\llabel{key_ineq_2}
\begin{split}
& \textstyle \cL(\Theta_n) - \cL(\Theta_{n - 1}) + C^{-1} \gamma_n \norm{\cG(\Theta_{n - 1})}^{\kappa} \le \cL(\Theta_n) - \cL(\Theta_{n - 1}) + \gamma_n \spro{\Gamma_n, \cG(\Theta_{n - 1})} \\
& \textstyle \le \fC \gamma_n \bigl[\sigma_n \norm{\cG(\Theta_{n - 1})} + (\gamma_n)^{\delta} \norm{\Gamma_n + \sigma_n \mu_n}^{1 + \delta}\bigr] \\
& \textstyle \le \fC C^{1 + \delta} \gamma_n \bigl[\sigma_n \norm{\cG(\Theta_{n - 1})} + (\gamma_n)^{\delta} (\norm{\cG(\Theta_{n - 1})}^{\kappa - 1} + \sigma_n)^{1 + \delta}\bigr]\dott
\end{split}
\end{equation}
}
\argument{\lref{Loja}; \lref{N_C_1}; \lref{B}; \lref{scrN_1}; the fact that $\kappa \ge 1 + \frac{1}{1 + \delta}$; the fact that $0 < \tau \le 1$}{that for all $n \in \scrN_1$ it holds that
\begin{equation}\llabel{step_to_bbL_lower_bound_scrN_1}
\begin{split}
& \textstyle \sigma_n \norm{\cG(\Theta_{n - 1})} + (\gamma_n)^{\delta} (\norm{\cG(\Theta_{n - 1})}^{\kappa - 1} + \sigma_n)^{1 + \delta} \\
& \textstyle \le \sigma_n B^{\frac{1}{\kappa - 1}} [\sigma_n + (\gamma_n)^{\delta}]^{\frac{1}{\kappa - 1}} + (\gamma_n)^{\delta} (B [\sigma_n + (\gamma_n)^{\delta}] + \sigma_n)^{1 + \delta} \\
& \textstyle \le \sigma_n 2^{1 + \delta} B^{\frac{1}{\kappa - 1}} [\sigma_n + (\gamma_n)^{\delta}]^{\frac{1}{\kappa - 1}} + (2B)^{1 + \delta} (\gamma_n)^{\delta} [\sigma_n + (\gamma_n)^{\delta}]^{1 + \delta} \\
& \textstyle = 2^{1 + \delta} B^{\frac{1}{\kappa - 1}} [\sigma_n + (\gamma_n)^{\delta}]^{\frac{1}{\kappa - 1}} \bigl[\sigma_n + (\gamma_n)^{\delta} (B [\sigma_n + (\gamma_n)^{\delta}])^{1 + \delta - \frac{1}{\kappa - 1}}\bigr] \\
& \textstyle = 2^{1 + \delta} B^{\frac{1}{\kappa - 1}} [\sigma_n + (\gamma_n)^{\delta}]^{\frac{1}{\kappa - 1}} \bigl[\sigma_n + (\gamma_n)^{\delta} (B [\sigma_n + (\gamma_n)^{\delta}])^{\frac{1 + \delta}{\kappa - 1} (\kappa - 1 - \frac{1}{1 + \delta})}\bigr] \\
& \textstyle \le 2^{1 + \delta} B^{\frac{1}{\kappa - 1}} [\sigma_n + (\gamma_n)^{\delta}]^{\frac{1}{\kappa - 1}} [\sigma_n + (\gamma_n)^{\delta}] = 2^{1 + \delta} B^{\frac{1}{\kappa - 1}} [\sigma_n + (\gamma_n)^{\delta}]^{\frac{\kappa}{\kappa - 1}}.
\end{split}
\end{equation}
}
\argument{\lref{Loja}; \lref{B}; \lref{scrN_2}; \lref{scrN_3}; the fact that $0 < \tau \le 1$}{that for all $n \in \scrN_2 \cup \scrN_3$ it holds that
\begin{equation}\llabel{step_to_bbL_lower_bound_scrN_2_3}
\begin{split}
& \textstyle \sigma_n \norm{\cG(\Theta_{n - 1})} + (\gamma_n)^{\delta} (\norm{\cG(\Theta_{n - 1})}^{\kappa - 1} + \sigma_n)^{1 + \delta} \\
& \textstyle \le \sigma_n \norm{\cG(\Theta_{n - 1})} + (\gamma_n)^{\delta} (\norm{\cG(\Theta_{n - 1})}^{\kappa - 1} + B[\sigma_n + (\gamma_n)^{\delta}])^{1 + \delta} \\
& \textstyle \le \sigma_n \norm{\cG(\Theta_{n - 1})} + 2^{1 + \delta} (\gamma_n)^{\delta} \norm{\cG(\Theta_{n - 1})}^{(\kappa - 1)(1 + \delta)} \\
& \textstyle = \sigma_n \norm{\cG(\Theta_{n - 1})} + 2^{1 + \delta} (\gamma_n)^{\delta} \norm{\cG(\Theta_{n - 1})} \norm{\cG(\Theta_{n - 1})}^{(1 + \delta)(\kappa - 1 - \frac{1}{1 + \delta})} \\
& \textstyle \le \sigma_n \norm{\cG(\Theta_{n - 1})} + 2^{1 + \delta} (\gamma_n)^{\delta} \norm{\cG(\Theta_{n - 1})} C^{1 + \delta} \\
& \textstyle \le (2 C)^{1 + \delta} \norm{\cG(\Theta_{n - 1})} [\sigma_n + (\gamma_n)^{\delta}] \le (2 C)^{1 + \delta} \norm{\cG(\Theta_{n - 1})} B^{-1} \norm{\cG(\Theta_{n - 1})}^{\kappa - 1} \\
& \textstyle = (2 C)^{1 + \delta} B^{-1} \norm{\cG(\Theta_{n - 1})}^{\kappa}.
\end{split}
\end{equation}
}
\argument{\lref{Loja}; \lref{B}; \lref{key_ineq_2}; \lref{step_to_bbL_lower_bound_scrN_1}}{for all $n \in \scrN_1$ that
\begin{equation}\llabel{key_ineq_scrN_1}
\begin{split}
& \textstyle \cL(\Theta_n) - \cL(\Theta_{n - 1}) + (2 C)^{-1} \gamma_n \norm{\cG(\Theta_{n - 1})}^{\kappa} \\
& \textstyle \le \cL(\Theta_n) - \cL(\Theta_{n - 1}) + C^{-1} \gamma_n \norm{\cG(\Theta_{n - 1})}^{\kappa} \\
& \textstyle \le \fC C^{1 + \delta} \gamma_n \bigl[\sigma_n \norm{\cG(\Theta_{n - 1})} + (\gamma_n)^{\delta} (\norm{\cG(\Theta_{n - 1})}^{\kappa - 1} + \sigma_n)^{1 + \delta}\bigr] \\
& \textstyle \le \fC (2 C)^{1 + \delta} B^{\frac{1}{\kappa - 1}} \gamma_n [\sigma_n + (\gamma_n)^{\delta}]^{\frac{\kappa}{\kappa - 1}} \le 2^{1 + \delta} C^{2 + \delta} B^{\frac{1}{\kappa - 1}} \gamma_n [\sigma_n + (\gamma_n)^{\delta}]^{\frac{\kappa}{\kappa - 1}} \\
& \textstyle = 2^{1 + \delta} C^{2 + \delta} B^{\frac{\kappa}{\kappa - 1}} B^{-1} \gamma_n [\sigma_n + (\gamma_n)^{\delta}]^{\frac{\kappa}{\kappa - 1}} \\
& \textstyle \le 2^{1 + \delta} C^{2 + \delta} B^{\frac{\kappa}{\kappa - 1}} (2^{3 + \delta} C^{4 + 2 \delta})^{-1} \gamma_n [\sigma_n + (\gamma_n)^{\delta}]^{\frac{\kappa}{\kappa - 1}} \\
& \textstyle = 2^{-2} C^{-2 - \delta} B^{\frac{\kappa}{\kappa - 1}} \gamma_n [\sigma_n + (\gamma_n)^{\delta}]^{\frac{\kappa}{\kappa - 1}} \le (2 C)^{-2} B^{\frac{\kappa}{\kappa - 1}} \gamma_n [\sigma_n + (\gamma_n)^{\delta}]^{\frac{\kappa}{\kappa - 1}}\dott
\end{split}
\end{equation}
}
\argument{\lref{Loja}; \lref{B}; \lref{key_ineq_2}; \lref{step_to_bbL_lower_bound_scrN_2_3}}{for all $n \in \scrN_2 \cup \scrN_3$ that
\begin{equation}\llabel{key_ineq_scrN_2_3}
\begin{split}
& \textstyle \cL(\Theta_n) - \cL(\Theta_{n - 1}) + C^{-1} \gamma_n \norm{\cG(\Theta_{n - 1})}^{\kappa} \\
& \textstyle \le \fC C^{1 + \delta} \gamma_n \bigl[\sigma_n \norm{\cG(\Theta_{n - 1})} + (\gamma_n)^{\delta} (\norm{\cG(\Theta_{n - 1})}^{\kappa - 1} + \sigma_n)^{1 + \delta}\bigr] \\
& \textstyle \le \fC 2^{1 + \delta} C^{2 + 2 \delta} B^{-1} \gamma_n \norm{\cG(\Theta_{n - 1})}^{\kappa} \le 2^{1 + \delta} C^{3 + 2 \delta} B^{-1} \gamma_n \norm{\cG(\Theta_{n - 1})}^{\kappa} \\
& \textstyle \le 2^{1 + \delta} C^{3 + 2 \delta} (2^{3 + \delta} C^{4 + 2 \delta})^{-1} \gamma_n \norm{\cG(\Theta_{n - 1})}^{\kappa} = 2^{-2} C^{-1} \gamma_n \norm{\cG(\Theta_{n - 1})}^{\kappa} \\
& \textstyle \le (2 C)^{-1} \gamma_n \norm{\cG(\Theta_{n - 1})}^{\kappa} + (2 C)^{-2} B^{\frac{\kappa}{\kappa - 1}} \gamma_n [\sigma_n + (\gamma_n)^{\delta}]^{\frac{\kappa}{\kappa - 1}}\dott
\end{split}
\end{equation}
}
\argument{\lref{key_ineq_scrN_1}; \lref{key_ineq_scrN_2_3}; the fact that $\scrN_1 \cup \scrN_2 \cup \scrN_3 = \N \cap [N, \infty)$}{that for all $n \in \allowbreak \N \cap \allowbreak [N, \allowbreak \infty)$ it holds that
\begin{equation}\llabel{result_item:C}
\textstyle \cL(\Theta_n) - \cL(\Theta_{n - 1}) + (2 C)^{-1} \gamma_n \norm{\cG(\Theta_{n - 1})}^{\kappa} \le (2 C)^{-2} B^{\frac{\kappa}{\kappa - 1}} \gamma_n [\sigma_n + (\gamma_n)^{\delta}]^{\frac{\kappa}{\kappa - 1}} \dott
\end{equation}
}
Let $\bbL \colon \N_0 \to \R$ satisfy for all $n \in \N_0$ that
\begin{equation}\llabel{bbL}
\textstyle \bbL_n = \cL(\Theta_n) - \cL(\vartheta) + (2 C)^{-2} B^{\frac{\kappa}{\kappa - 1}} \sum_{j = n + 1}^{\infty} \gamma_j [\sigma_j + (\gamma_j)^{\delta}]^{\frac{\kappa}{\kappa - 1}}\dott
\end{equation}
\startnewargseq
\argument{\lref{result_item:C}; \lref{bbL}}{that for all $n \in \N \cap [N, \infty)$ it holds that
\begin{equation}\llabel{bbL_diff}
\begin{split}
& \textstyle \bbL_n - \bbL_{n - 1} \\
& \textstyle = \cL(\Theta_n) - \cL(\vartheta) + (2 C)^{-2} B^{\frac{\kappa}{\kappa - 1}} \sum_{j = n + 1}^{\infty} \gamma_j [\sigma_j + (\gamma_j)^{\delta}]^{\frac{\kappa}{\kappa - 1}} \\
& \textstyle \quad - \cL(\Theta_{n - 1}) + \cL(\vartheta) - (2 C)^{-2} B^{\frac{\kappa}{\kappa - 1}} \sum_{j = n}^{\infty} \gamma_j [\sigma_j + (\gamma_j)^{\delta}]^{\frac{\kappa}{\kappa - 1}} \\
& \textstyle = \cL(\Theta_n) - \cL(\Theta_{n - 1}) - (2 C)^{-2} B^{\frac{\kappa}{\kappa - 1}} \gamma_n [\sigma_n + (\gamma_n)^{\delta}]^{\frac{\kappa}{\kappa - 1}} \\
& \textstyle \le - (2 C)^{-1} \gamma_n \norm{\cG(\Theta_{n - 1})}^{\kappa}.
\end{split}
\end{equation}
}
\argument{\lref{gamma_sigma_sum}; \lref{limit_cL}; \lref{vartheta}; \lref{bbL}; \lref{bbL_diff}}{for all $n \in \N \cap [N, \infty)$ that
\begin{equation}\llabel{bbL_ineq}
\textstyle \limsup_{k \to \infty} \abs{\bbL_k} = 0 \le (2 C)^{-1} \gamma_n \norm{\cG(\Theta_{n - 1})}^{\kappa} \le \bbL_{n - 1} - \bbL_n.
\end{equation}
}
\argument{\lref{bbL_ineq}}{that for all $n, m \in \N \cap [N, \infty)$ with $n \le m$ it holds that
\begin{equation}\llabel{bbL_telescopic}
\textstyle (2 C)^{-1} \sum_{j = n}^{m - 1} \gamma_j \norm{\cG(\Theta_{j - 1})}^{\kappa} \le \sum_{j = n}^{m - 1} (\bbL_{j - 1} - \bbL_j) = \bbL_{n - 1} - \bbL_{m - 1}.
\end{equation}
}
\argument{\lref{bbL_ineq}; \lref{bbL_telescopic}}{that for all $n \in \N \cap [N, \infty)$ it holds that
\begin{equation}\llabel{bbL_lower_bound}
\textstyle (2 C)^{-1} \sum_{j = n}^{\infty} \gamma_j \norm{\cG(\Theta_{j - 1})}^{\kappa} \le \lim_{m \to \infty} (\bbL_{n - 1} - \bbL_{m - 1}) = \bbL_{n - 1}.
\end{equation}
}
\argument{\lref{Theta}; \lref{Loja}; \lref{N_C_1}; the triangle inequality}[verbs=d]{that for all $n \in \N \cap [N, \infty)$ it holds that
\begin{equation}\llabel{Theta_ineq}
\begin{split}
& \textstyle \norm{\Theta_n - \Theta_{n - 1}} = \gamma_n \norm{\Gamma_n + \sigma_n \mu_n} \le \gamma_n (\norm{\Gamma_n} + \sigma_n \norm{\mu_n}) \\
& \textstyle \le C \gamma_n \norm{\cG(\Theta_{n - 1})}^{\kappa - 1} + \gamma_n \sigma_n C = C \gamma_n (\norm{\cG(\Theta_{n - 1})}^{\kappa - 1} + \sigma_n).
\end{split}
\end{equation}
}
\argument{\lref{Loja}; \lref{B}; \lref{scrN_1}; \lref{scrN_2}; \lref{Theta_ineq}; the fact that $0 < \tau \le 1$}{for all $n \in \scrN_1 \cup \scrN_2$ that
\begin{equation}\llabel{Theta_scrN_1_2}
\begin{split}
& \textstyle \norm{\Theta_n - \Theta_{n - 1}} \le C \gamma_n (\norm{\cG(\Theta_{n - 1})}^{\kappa - 1} + \sigma_n) \\
& \textstyle \le C \gamma_n (B [\sigma_n + (\gamma_n)^{\delta}]^{\tau} + \sigma_n) \le 2 B C \gamma_n [\sigma_n + (\gamma_n)^{\delta}]^{\tau}.
\end{split}
\end{equation}
}
\argument{\lref{bbL_ineq}; \lref{Theta_ineq}}{for all $n \in \allowbreak \N \cap \allowbreak [N, \allowbreak \infty)$ that
\begin{equation}\llabel{Theta_bbL_bound}
\begin{split}
& \textstyle \bbL_{n - 1} - \bbL_n \ge (2 C)^{-1} \gamma_n \norm{\cG(\Theta_{n - 1})}^{\kappa} = 2^{-1} C^{-2} \norm{\cG(\Theta_{n - 1})} (C \gamma_n \norm{\cG(\Theta_{n - 1})}^{\kappa - 1}) \\
& \ge 2^{-1} C^{-2} \norm{\cG(\Theta_{n - 1})} (\norm{\Theta_n - \Theta_{n - 1}} - \gamma_n \sigma_n C).
\end{split}
\end{equation}
}
\argument{\lref{Loja}; \lref{Theta_bbL_bound}}{for all $n \in \allowbreak \N \cap \allowbreak [N, \allowbreak \infty)$ with $\norm{\Theta_{n - 1} - \vartheta} < \eps$ that
\begin{equation}\llabel{Theta_scrN_3_0}
\begin{split}
& \textstyle \bbL_{n - 1} - \bbL_n \ge 2^{-1} C^{-2} \norm{\cG(\Theta_{n - 1})} (\norm{\Theta_n - \Theta_{n - 1}} - \gamma_n \sigma_n C) \\
& \textstyle \ge 2^{-1} C^{-3} \abs{\cL(\Theta_{n - 1}) - \cL(\vartheta)}^{\rho} (\norm{\Theta_n - \Theta_{n - 1}} - \gamma_n \sigma_n C).
\end{split}
\end{equation}
}
\argument{\lref{scrN_3}; \lref{Cesaro_ineq}; \lref{bbL}; \lref{bbL_lower_bound}; the fact that for all $n \in \N$ it holds that $\gamma_n > 0$}{for all $n \in \scrN_3$ that
\begin{equation}\llabel{step_to_bbL_lower_bound_scrN_3}
\begin{split}
& \textstyle \bbL_{n - 1} \ge (2 C)^{-1} \sum_{j = n}^{\infty} \gamma_j \norm{\cG(\Theta_{j - 1})}^{\kappa} \ge (2 C)^{-1} \gamma_n \norm{\cG(\Theta_{n - 1})}^{\kappa} \\
& \textstyle \ge (2 C)^{-1} B^{\frac{\kappa}{\kappa - 1}} \gamma_n [\sigma_n + (\gamma_n)^{\delta}]^{\frac{\kappa \tau}{\kappa - 1}} \ge (2 C)^{-1} B^{\frac{\kappa}{\kappa - 1}} C^{-1} \sum_{j = n}^{\infty} \gamma_j [\sigma_j + (\gamma_j)^{\delta}]^\frac{\kappa}{\kappa - 1} \\
& \textstyle = 2 (2 C)^{-2} B^{\frac{\kappa}{\kappa - 1}} \sum_{j = n}^{\infty} \gamma_j [\sigma_j + (\gamma_j)^{\delta}]^\frac{\kappa}{\kappa - 1} = 2 (\bbL_{n - 1} - \cL(\Theta_{n - 1}) + \cL(\vartheta))
\end{split}
\end{equation}
and
\begin{equation}\llabel{bbL_scrN_3_positive}
\textstyle 2C \bbL_{n - 1} \ge \sum_{j = n}^{\infty} \gamma_j \norm{\cG(\Theta_{j - 1})}^{\kappa} \ge \gamma_n \norm{\cG(\Theta_{n - 1})}^{\kappa} \ge B^{\frac{\kappa}{\kappa - 1}} \gamma_n [\sigma_n + (\gamma_n)^{\delta}]^{\frac{\kappa \tau}{\kappa - 1}} > 0\dott
\end{equation}
}
\argument{\lref{bbL_scrN_3_positive}}{that for all $n \in \scrN_3$ it holds that
\begin{equation}\llabel{cL_lower_bound_bbL}
\textstyle \cL(\Theta_{n - 1}) - \cL(\vartheta) \ge 2^{-1} \bbL_{n - 1} > 0\dott
\end{equation}
}
\argument{\lref{Theta_scrN_3_0}; \lref{cL_lower_bound_bbL}}{for all $n \in \scrN_3$ with $\norm{\Theta_{n - 1} - \vartheta} < \eps$ that
\begin{equation}\llabel{Theta_scrN_3_1}
\begin{split}
& \textstyle \bbL_{n - 1} - \bbL_n \ge 2^{-1} C^{-3} \abs{\cL(\Theta_{n - 1}) - \cL(\vartheta)}^{\rho} (\norm{\Theta_n - \Theta_{n - 1}} - \gamma_n \sigma_n C) \\
& \textstyle \ge C^{-3} \frac{1}{2^{\rho + 1}} (\bbL_{n - 1})^{\rho} (\norm{\Theta_n - \Theta_{n - 1}} - \gamma_n \sigma_n C).
\end{split}
\end{equation}
}
\argument{\lref{bbL_scrN_3_positive}; \lref{Theta_scrN_3_1}; the fact that for all $n \in \N \cap [N, \infty)$ it holds that $\bbL_{n - 1} \ge \allowbreak \bbL_{n} \allowbreak \ge \allowbreak 0$}{that for all $n \in \scrN_3$ with $\norm{\Theta_{n - 1} - \vartheta} < \eps$ it holds that
\begin{equation}\llabel{Theta_scrN_3}
\begin{split}
& \textstyle \norm{\Theta_n - \Theta_{n - 1}} \le 2^{\rho + 1} C^3 \frac{\bbL_{n - 1} - \bbL_n}{(\bbL_{n - 1})^{\rho}} + \gamma_n \sigma_n C = \gamma_n \sigma_n C + 2^{\rho + 1} C^3 \int_{\bbL_n}^{\bbL_{n - 1}} (\bbL_{n - 1})^{- \rho} \d u \\
& \textstyle \le \gamma_n \sigma_n C + 2^{\rho + 1} C^3 \int_{\bbL_n}^{\bbL_{n - 1}} u^{- \rho} \d u = \gamma_n \sigma_n C + 2^{\rho + 1} C^3 \frac{1}{1 - \rho} ([\bbL_{n - 1}]^{1 - \rho} - [\bbL_n]^{1 - \rho}).
\end{split}
\end{equation}
}
\argument{\lref{Loja}; \lref{B}; \lref{Theta_scrN_1_2}; \lref{Theta_scrN_3}; the fact that $\scrN_1 \cup \scrN_2 \cup \scrN_3 = \N \cap [N, \infty)$; the fact that $0 < \tau \le 1$}{that for all $n \in \N \cap [N, \infty)$ with $\norm{\Theta_{n - 1} - \vartheta} < \eps$ it holds that
\begin{equation}\llabel{Theta_upper_bound_gen}
\begin{split}
& \textstyle \norm{\Theta_n - \Theta_{n - 1}} \le 2 B C \gamma_n [\sigma_n + (\gamma_n)^{\delta}]^{\tau} + \gamma_n \sigma_n C + 2^{\rho + 1} C^3 \frac{1}{1 - \rho} ([\bbL_{n - 1}]^{1 - \rho} - [\bbL_n]^{1 - \rho}) \\
& \textstyle \le 3 B C \gamma_n [\sigma_n + (\gamma_n)^{\delta}]^{\tau} + 2^{\rho + 1} C^3 \frac{1}{1 - \rho} ([\bbL_{n - 1}]^{1 - \rho} - [\bbL_n]^{1 - \rho}) \\
& \textstyle \le D \bigl(\gamma_n [\sigma_n + (\gamma_n)^{\delta}]^{\tau} + [\bbL_{n - 1}]^{1 - \rho} - [\bbL_n]^{1 - \rho}\bigr).
\end{split}
\end{equation}
}
\argument{\lref{Theta}; \lref{vartheta}; \lref{bbL_ineq}; the fact that $0 < \rho < 1$}{that there exists $M \in \N \cap [N, \infty)$ which satisfy for all $n \in \N \cap [M, \infty)$ that
\begin{equation}\llabel{M}
\textstyle \norm{\Theta_M - \vartheta} \le \frac{\eps}{4} \qqandqq [\bbL_{n - 1}]^{1 - \rho} + \sum_{j = n}^{\infty} \gamma_j [\sigma_j + (\gamma_j)^{\delta}]^{\tau} \le \frac{\eps}{4 D}.
\end{equation}
}
\argument{\lref{Theta_upper_bound_gen}; \lref{M}; the triangle inequality; the fact that for all $n \in \N \cap [N, \infty)$ it holds that $\bbL_{n - 1} \ge 0$}{for all $s \in \N_0$ with $\Forall t \in \{M, M + 1, \dots, M + s\} \colon \norm{\Theta_t - \vartheta} \le \eps / 2$ that
\begin{equation}\llabel{Theta_induction_step}
\begin{split}
& \textstyle \norm{\Theta_{M + s + 1} - \vartheta} \le \norm{\Theta_M - \vartheta} + \sum_{j = M + 1}^{M + s + 1} \norm{\Theta_{j} - \Theta_{j - 1}} \\
& \textstyle \le \frac{\eps}{4} + D \sum_{j = M + 1}^{M + s + 1} \bigl(\gamma_j [\sigma_j + (\gamma_j)^{\delta}]^{\tau} + [\bbL_{j - 1}]^{1 - \rho} - [\bbL_j]^{1 - \rho}\bigr) \\
& \textstyle = \frac{\eps}{4} + D \bigl([\bbL_{M}]^{1 - \rho} - [\bbL_{M + s + 1}]^{1 - \rho} + \sum_{j = M + 1}^{M + s + 1} \gamma_j [\sigma_j + (\gamma_j)^{\delta}]^{\tau}\bigr) \\
& \textstyle \le \frac{\eps}{4} + D \bigl([\bbL_{M}]^{1 - \rho} + \sum_{j = M + 1}^{\infty} \gamma_j [\sigma_j + (\gamma_j)^{\delta}]^{\tau}\bigr) \le \frac{\eps}{4} + \frac{D \eps}{4 D} = \frac{\eps}{2}.
\end{split}
\end{equation}
}
\argument{\lref{M}; \lref{Theta_induction_step}; the induction}{that
\begin{equation}\llabel{Theta_eps_M}
\textstyle \Forall n \in \N \cap[M, \infty) \colon \norm{\Theta_n - \vartheta} \le \frac{\eps}{2} < \eps.
\end{equation}
}
\argument{\lref{Theta}; \lref{bbL_ineq}; \lref{Theta_upper_bound_gen}; \lref{Theta_eps_M}}{that for all $n \in \N \cap [M + 1, \infty)$ it holds that
\begin{equation}\llabel{Theta_convergence}
\begin{split}
\textstyle \sum_{j = n}^{\infty} \norm{\Theta_j - \Theta_{j - 1}} & \textstyle \le D \sum_{j = n}^{\infty} \bigl(\gamma_j [\sigma_j + (\gamma_j)^{\delta}]^{\tau} + [\bbL_{j - 1}]^{1 - \rho} - [\bbL_j]^{1 - \rho}\bigr) \\
& \textstyle = D [\bbL_{n - 1}]^{1 - \rho} + D \sum_{j = n}^{\infty} \gamma_j [\sigma_j + (\gamma_j)^{\delta}]^{\tau} < \infty.
\end{split}
\end{equation}
}
\argument{\lref{vartheta}; \lref{Theta_convergence}}{that
\begin{equation}\llabel{Theta_limit}
\textstyle \limsup_{n \to \infty} \norm{\Theta_n - \vartheta} = 0.
\end{equation}
}
\argument{\lref{bbL_lower_bound}; \lref{Theta_limit}; the fact that $\cG$ is H\"{o}lder continuous}{that
\begin{equation}\llabel{vartheta_critical_point}
\textstyle \bigl[1 - \mathbbm{1}_{\R}\bigl(\sum_{n = 1}^{\infty} \gamma_n\bigr)\bigr] \norm{\cG(\vartheta)} = 0.
\end{equation}
}
Let $\scrC \in (0, \infty)$ satisfy 
\begin{gather}\llabel{scrC1}
\textstyle \scrC \ge \max\{D, 2C, (2 C)^{-2} B^{\frac{\kappa}{\kappa - 1}}\}, \\ \llabel{scrC2_0}
\textstyle \scrC \ge \max_{n \in \{1, 2, \dots, M\}} \bigl[(\gamma_n)^{-1} [\sigma_n + (\gamma_n)^{\delta}]^{\frac{- \kappa \tau}{\kappa - 1}}  \sum_{j = n}^{\infty} \gamma_j [\sigma_j + (\gamma_j)^{\delta}]^\frac{\kappa}{\kappa - 1}\bigr] \\
\llabel{scrC2}
\textstyle \scrC \ge \max_{n \in \{1, 2, \dots, M\}} \bigl[(\gamma_n)^{-1} [\sigma_n + (\gamma_n)^{\delta}]^{\frac{- \kappa}{\kappa - 1}} (\cL(\Theta_n) - \cL(\Theta_{n - 1}) + \gamma_n \norm{\cG(\Theta_{n - 1})}^{\kappa})\bigr], \\
\llabel{scrC3}
\textstyle \scrC \ge \max_{n \in \{1, 2, \dots, M\}} \bigl[\bigl(\sum_{j = n}^{\infty} \gamma_j [\sigma_j + (\gamma_j)^{\delta}]^{\frac{\kappa}{\kappa - 1}}\bigr)^{-1} (\cL(\vartheta) - \cL(\Theta_{n - 1}))\bigr], \\
\llabel{scrC4}
\textstyle \andqq \scrC \ge \max_{n \in \{1, 2, \dots, M\}} \bigl[\bigl(\sum_{j = n}^{\infty} \gamma_j [\sigma_j + (\gamma_j)^{\delta}]^{\tau}\bigr)^{-1} \norm{\Theta_{n - 1} - \vartheta}\bigr]\dott
\end{gather}
\startnewargseq
\argument{\lref{Loja}; \lref{scrC1}}{that for all $\theta \in \R^{\fd}$ with $\norm{\theta - \vartheta} < \eps$ it holds that
\begin{equation}\llabel{proof:item:Loja}
\textstyle \abs{\cL(\theta) - \cL(\vartheta)}^{\rho} \le \scrC \norm{\cG(\theta)}.
\end{equation}
}
\argument{\lref{Theta_limit}; \lref{vartheta_critical_point}}{that
\begin{equation}\llabel{proof:item:vartheta}
\textstyle \bigl[1 - \mathbbm{1}_{\R}\bigl(\sum_{n = 1}^{\infty} \gamma_n\bigr)\bigr] \norm{(\nabla \cL)(\vartheta)} + \limsup_{n \to \infty} \norm{\Theta_n - \vartheta} = 0.
\end{equation}
}
\argument{\lref{Loja}; \lref{scrC2}}{that for all $n \in \{1, 2, \dots, M\}$ it holds that
\begin{equation}\llabel{proof:item:C:before_M}
\textstyle \scrC^{-1} \gamma_n \norm{\cG(\Theta_{n - 1})}^{\kappa} \le \gamma_n \norm{\cG(\Theta_{n - 1})}^{\kappa} \le - \cL(\Theta_n) + \cL(\Theta_{n - 1}) + \scrC \gamma_n [\sigma_n + (\gamma_n)^{\delta}]^{\frac{\kappa}{\kappa - 1}}\dott
\end{equation}
}
\argument{\lref{result_item:C}; \lref{scrC1}; the fact that $M \ge N$}{that for all $n \in \allowbreak \N \allowbreak \cap \allowbreak (M, \allowbreak \infty)$ it holds that
\begin{equation}\llabel{proof:item:C:after_M}
\begin{split}
& \textstyle \cL(\Theta_n) - \cL(\Theta_{n - 1}) + \frac{1}{\scrC} \gamma_n \norm{\cG(\Theta_{n - 1})}^{\kappa} \le \cL(\Theta_n) - \cL(\Theta_{n - 1}) + \frac{1}{2 C} \gamma_n \norm{\cG(\Theta_{n - 1})}^{\kappa} \\
& \textstyle \le (2 C)^{-2} B^{\frac{\kappa}{\kappa - 1}} \gamma_n [\sigma_n + (\gamma_n)^{\delta}]^{\frac{\kappa}{\kappa - 1}} \le \scrC \gamma_n [\sigma_n + (\gamma_n)^{\delta}]^{\frac{\kappa}{\kappa - 1}}\dott
\end{split}
\end{equation}
}
\argument{\lref{proof:item:C:before_M}; \lref{proof:item:C:after_M}}{that for all $n \in \N$ it holds that
\begin{equation}\llabel{proof:item:C}
\textstyle \cL(\Theta_n) - \cL(\Theta_{n - 1}) + \scrC^{-1} \gamma_n \norm{\cG(\Theta_{n - 1})}^{\kappa} \le \scrC \gamma_n [\sigma_n + (\gamma_n)^{\delta}]^{\frac{\kappa}{\kappa - 1}}\dott
\end{equation}
}
\argument{\lref{bbL}; \lref{scrC1}; \lref{scrC3}; the fact that $M \ge N$; the fact that for all $n \in \allowbreak \N \cap \allowbreak [N, \allowbreak \infty)$ it holds that $\bbL_{n - 1} \ge \allowbreak \bbL_n \allowbreak \ge 0$}{that for all $n \in \N$ it holds that
\begin{equation}\llabel{proof:Theta_bound_bbL}
\textstyle \cL(\Theta_{n - 1}) - \cL(\vartheta) + \scrC \!\sum_{j = n}^{\infty} \gamma_j [\sigma_j + (\gamma_j)^{\delta}]^{\frac{\kappa}{\kappa - 1}} \ge [\mathbbm{1}_{(M, \infty)}(n)] \bbL_{n - 1} \ge 0\dott
\end{equation}
}
\argument{\lref{scrC4}}{that for all $n \in \{1, 2, \dots, M\}$ it holds that
\begin{equation}\llabel{proof:item:rate:before_M}
\textstyle \norm{\Theta_{n - 1} - \vartheta} \le \scrC \sum_{j = n}^{\infty} \gamma_j [\sigma_j + (\gamma_j)^{\delta}]^{\tau}.
\end{equation}
}
\argument{\lref{Theta_convergence}; \lref{Theta_limit}; \lref{scrC1}; \lref{proof:Theta_bound_bbL}; the triangle inequality}{for all $n \in \allowbreak \N \cap \allowbreak (M, \allowbreak \infty)$ that
\begin{equation}\llabel{proof:item:rate:after_M}
\begin{split}
& \textstyle \norm{\Theta_{n - 1} - \vartheta} \le \sum_{j = n}^{\infty} \norm{\Theta_j - \Theta_{j - 1}} \le D [\bbL_{n - 1}]^{1 - \rho} + D \sum_{j = n}^{\infty} \gamma_j [\sigma_j + (\gamma_j)^{\delta}]^{\tau} \\
& \textstyle \le \scrC [\bbL_{n - 1}]^{1 - \rho} + \scrC \sum_{j = n}^{\infty} \gamma_j [\sigma_j + (\gamma_j)^{\delta}]^{\tau} \\
& \textstyle \le \scrC \bigl[\cL(\Theta_{n - 1}) - \cL(\vartheta) + \scrC \!\sum_{j = n}^{\infty} \gamma_j [\sigma_j + (\gamma_j)^{\delta}]^{\frac{\kappa}{\kappa - 1}}\bigr]^{1 - \rho} + \scrC \sum_{j = n}^{\infty} \gamma_j [\sigma_j + (\gamma_j)^{\delta}]^{\tau}\dott
\end{split}
\end{equation}
}
\argument{\lref{proof:Theta_bound_bbL}; \lref{proof:item:rate:before_M}; \lref{proof:item:rate:after_M}}{that for all $n \in \N$ it holds that
\begin{multline}\llabel{proof:item:rate}
\textstyle \norm{\Theta_{n - 1} - \vartheta} \\
\textstyle \le \scrC \bigl[\cL(\Theta_{n - 1}) - \cL(\vartheta) + \scrC \!\sum_{j = n}^{\infty} \gamma_j [\sigma_j + (\gamma_j)^{\delta}]^{\frac{\kappa}{\kappa - 1}}\bigr]^{1 - \rho} + \scrC \!\sum_{j = n}^{\infty} \gamma_j [\sigma_j + (\gamma_j)^{\delta}]^{\tau}\dott
\end{multline}
}
\argument{\lref{Cesaro_ineq}; \lref{scrC1}; \lref{scrC2_0}; \lref{proof:item:Loja}; \lref{proof:item:vartheta}; \lref{proof:item:C}; \lref{proof:Theta_bound_bbL}; \lref{proof:item:rate}}[verbs=e]{\lref{Lojaaa,descending,rate_prep,rate,result}. }
\end{aproof}

\cfclear
\begin{savenotes}
\begin{samepage}
\begin{tcolorbox}[colback=white!95!gray,
                  colframe=black,
                  boxrule=0.5pt,
                  sharp corners,
                  enhanced,
                 ]
\begin{athm}{lemma}{lr_for_rate}
Let $\delta \in (0, 2]$, $\tau \in [0, \infty)$ and let $\gamma \colon \N \to (0, \infty)$ satisfy
\begin{equation}\llabel{gamma}
\textstyle \liminf_{n \to \infty} \bigl[\frac{(\gamma_n)^{1 + \tau} - (\gamma_{n + 1})^{1 + \tau}}{(\gamma_n)^{1 + \delta}}\bigr] > 0 \dott
\end{equation}
Then
\begin{equation}\llabel{result}
\textstyle \sup_{n \in \N} \bigl[(\gamma_n)^{1 + \tau} \sum_{j = 1}^{n} \gamma_j\bigr] < \infty \dott
\end{equation}
\end{athm}
\end{tcolorbox}
\end{samepage}
\end{savenotes}
\begin{aproof}
\argument{\lref{gamma}}{that there exists $N \in \N$ which satisfies for all $n \in \N \cap [N, \infty)$ that
\begin{equation}\llabel{gamma_monotone}
\textstyle \gamma_n \ge \gamma_{n + 1} \qqandqq \sup_{k \in \N} \bigl[(\gamma_k)^{2 - \delta + 2 \tau}\bigr] < \infty.
\end{equation}
}
\argument{\lref{gamma_monotone}}{that
\begin{equation}\llabel{gamma_bounded}
\textstyle \liminf\limits_{n \to \infty} \Bigl[\frac{(\gamma_n)^{1 + \delta}}{(\gamma_n)^{1 + \tau} (\gamma_{n + 1})^{2 + \tau}}\Bigr] \ge \liminf\limits_{n \to \infty} \Bigl[\frac{(\gamma_n)^{1 + \delta}}{(\gamma_n)^{3 + 2 \tau}}\Bigr] = \liminf\limits_{n \to \infty} \Bigl[\frac{1}{(\gamma_n)^{2 - \delta + 2 \tau}}\Bigr] > 0 \dott
\end{equation}
}
\argument{\lref{gamma}; \lref{gamma_bounded}}{that
\begin{equation}\llabel{telescope_lower}
\textstyle 0 < \liminf\limits_{n \to \infty} \Bigl[\frac{(\gamma_n)^{1 + \tau} - (\gamma_{n + 1})^{1 + \tau}}{(\gamma_n)^{1 + \delta}} \cdot \frac{(\gamma_n)^{1 + \delta}}{(\gamma_n)^{1 + \tau} (\gamma_{n + 1})^{2 + \tau}}\Bigr]\! = \liminf\limits_{n \to \infty} \Bigl[\frac{(\gamma_n)^{1 + \tau} - (\gamma_{n + 1})^{1 + \tau}}{\gamma_{n + 1} (\gamma_n)^{1 + \tau} (\gamma_{n + 1})^{1 + \tau}}\Bigr]\! \dott
\end{equation}
}
\argument{\lref{telescope_lower}}{that there exist $\eps \in (0, \infty)$, $N \in \N$ which satisfy for all $n \in \N \cap [N, \infty)$ that
\begin{equation}\llabel{telescope}
\textstyle \eps \gamma_{n + 1} \le \frac{(\gamma_n)^{1 + \tau} - (\gamma_{n + 1})^{1 + \tau}}{(\gamma_n)^{1 + \tau} (\gamma_{n + 1})^{1 + \tau}} = \frac{1}{(\gamma_{n + 1})^{1 +\tau}} - \frac{1}{(\gamma_n)^{1 + \tau}} \dott
\end{equation}
}
\argument{\lref{telescope}}{for all $n \in \N \cap [N, \infty)$ that
\begin{equation}\llabel{res}
\textstyle \eps \sum_{j = N}^{n} \gamma_{j + 1} \le \sum_{j = N}^{n} \Bigl[\frac{1}{(\gamma_{j + 1})^{1 + \tau}} - \frac{1}{(\gamma_j)^{1 + \tau}}\Bigr] = \frac{1}{(\gamma_{n + 1})^{1 + \tau}} - \frac{1}{(\gamma_N)^{1 + \tau}} < \frac{1}{(\gamma_{n + 1})^{1 + \tau}} \dott
\end{equation}
}
\argument{\lref{gamma_monotone}; \lref{res}}[verbs=e]{\lref{result}. }
\end{aproof}

\cfclear
\begin{savenotes}
\begin{samepage}
\begin{tcolorbox}[colback=white!95!gray,
                  colframe=black,
                  boxrule=0.5pt,
                  sharp corners,
                  enhanced,
                 ]
\begin{athm}{lemma}{rate_lemma}
Let $a \in (0, \infty)$, $b, \tau \in [0, \infty)$, let $\gamma \colon \N \to (0, \infty)$ satisfy
\begin{equation}\llabel{gamma}
\textstyle \sup_{n \in \N} \bigl[(\gamma_n)^{1 + \tau} \sum_{j = 1}^{n} \gamma_j\bigr] < \infty,
\end{equation}
and let $e \colon \N_0 \to [0, \infty)$ satisfy for all $n \in \N$ that $e_n (1 + a \gamma_n e_n) \le e_{n - 1} + b (\gamma_n)^{3 + 2 \tau}$. Then
\begin{equation}\llabel{result}
\textstyle \sup_{n \in \N} \bigl[e_n \sum_{j = 1}^{n} \gamma_j\bigr] < \infty.
\end{equation}
\end{athm}
\end{tcolorbox}
\end{samepage}
\end{savenotes}
\begin{aproof}
\argument{\lref{gamma}}{that there exist $C, \fC \in \R$ which satisfy for all $n \in \N$ that
\begin{equation}\llabel{C_fC}
\textstyle \gamma_n + (\gamma_n)^{1 + \tau} \sum_{j = 1}^{n} \gamma_j \le C \qqandqq \max\bigl\{b C^2, \frac{2 \gamma_1 + C}{a \gamma_1}, e_1 \gamma_1\bigr\} \le \fC \dott
\end{equation}
}
\startnewargseq
\argument{\lref{C_fC}}{that for all $n \in \N \cap [2, \infty)$ it holds that
\begin{equation}\llabel{temp_ineq}
\begin{split}
& \textstyle b (\gamma_n)^{3 + 2 \tau}\! + \frac{\fC}{\sum_{j = 1}^{n - 1}\! \gamma_j} - \frac{\fC}{\sum_{j = 1}^{n}\! \gamma_j} - \frac{a \gamma_n \fC^2}{[\sum_{j = 1}^{n}\! \gamma_j]^2} = b (\gamma_n)^{3 + 2 \tau}\! + \frac{\fC \gamma_n}{[\sum_{j = 1}^{n - 1}\! \gamma_j] [\sum_{j = 1}^{n}\! \gamma_j]} - \frac{a \gamma_n \fC^2}{[\sum_{j = 1}^{n}\! \gamma_j]^2} \\
& \textstyle = \frac{\gamma_n}{[\sum_{j = 1}^{n} \gamma_j]^2} \Bigl[b \bigl[(\gamma_n)^{1 + \tau} \sum_{j = 1}^{n}\! \gamma_j\bigr]^2 + \frac{\fC \sum_{j = 1}^{n}\! \gamma_j}{\sum_{j = 1}^{n - 1}\! \gamma_j} - a \fC^2\Bigr] \\
& \textstyle \le \frac{C}{(\gamma_1)^2} \Bigl[b C^2 + \frac{\fC (C + \sum_{j = 1}^{n - 1}\! \gamma_j)}{\sum_{j = 1}^{n - 1}\! \gamma_j} - a \fC^2\Bigr] \le \frac{C}{(\gamma_1)^2} \bigl[\fC + \fC (1 + \frac{C}{\gamma_1}) - a \fC^2\bigr] \\
& \textstyle = \frac{C \fC a}{(\gamma_1)^2} \bigl[\frac{2 \gamma_1 + C}{a \gamma_1} -\fC\bigr] \le 0 \dott
\end{split}
\end{equation}
}
\argument{\lref{temp_ineq}; the assumption that for all $n \in \N$ it holds that $e_n (1 + a \gamma_n e_n) \le e_{n - 1} + b (\gamma_n)^{3 + 2 \tau}$}{for all $n \in \N \cap [2, \infty)$ that
\begin{equation}\llabel{res}
\begin{split}
& \textstyle \bigl(e_n - \frac{\fC}{\sum_{j = 1}^{n} \gamma_j}\bigr) \bigl[1 + a \gamma_n \bigl(e_n + \frac{\fC}{\sum_{j = 1}^{n} \gamma_j}\bigr)\bigr] = e_n (1 + a \gamma_n e_n) - \frac{\fC}{\sum_{j = 1}^{n} \gamma_j} - \frac{a \gamma_n \fC^2}{[\sum_{j = 1}^{n} \gamma_j]^2} \\
& \textstyle \le \bigl(e_{n - 1} - \frac{\fC}{\sum_{j = 1}^{n - 1} \gamma_j}\bigr) + \bigl[b (\gamma_n)^{3 + 2 \tau} + \frac{\fC}{\sum_{j = 1}^{n - 1} \gamma_j} - \frac{\fC}{\sum_{j = 1}^{n} \gamma_j} - \frac{a \gamma_n \fC^2}{[\sum_{j = 1}^{n} \gamma_j]^2}\bigr] \\
& \textstyle \le e_{n - 1} - \frac{\fC}{\sum_{j = 1}^{n - 1} \gamma_j} \dott
\end{split}
\end{equation}
}
\argument{\lref{C_fC}; \lref{res}; the induction}[verbs=e]{\lref{result}. }
\end{aproof}

\cfclear
\begin{savenotes}
\begin{samepage}
\begin{tcolorbox}[colback=white!95!gray,
                  colframe=black,
                  boxrule=0.5pt,
                  sharp corners,
                  enhanced,
                 ]
\begin{athm}{prop}{theorem:gen_convergence_Loja}
Let $\fd \in \N$, $\delta, \tau \in (0, 1]$, let $\cL \colon \R^{\fd} \to \R$ be a\cfadd{def:Loja} \KL\ function, assume that $\nabla \cL$ is locally $\delta$-H\"{o}lder continuous, let $\bbA \colon \N \to \R^{\fd \times \fd}$, $\Theta \allowbreak \colon \allowbreak \N_0 \allowbreak \to \allowbreak \R^\fd$, $\mu \colon \N \to \R^{\fd}$, and $p \colon \allowbreak \N \allowbreak \to \R$ be bounded, let $\alpha \colon \N \to [0, 1]$, $\gamma \colon \N \to (0, \infty)$, and $\m \colon \allowbreak \N_0 \allowbreak \to \R^\fd$ satisfy for all $n \in \N$ that
\begin{gather}\llabel{m}
\textstyle \m_n = \alpha_n \m_{n - 1} + (1 - \alpha_n) \bigl[(\nabla \cL)(\Theta_{n - 1}) + (\gamma_n)^{\delta} \mu_n\bigr] \\
\llabel{Theta}
\textstyle \andqq \Theta_n = \Theta_{n - 1} - \gamma_n \bbA_n \bigl[p_n \m_n + (1 - p_n) (\nabla \cL)(\Theta_{n - 1})\bigr],
\end{gather}
assume for all $n \in \N$ that $\bbA_n - \bbI_{\fd}$ is symmetric positive semi-definite, and assume
\begin{equation}\llabel{bound_2}
\textstyle \limsup_{n \to \infty} \bigl[\frac{(\alpha_{n + 1})^{1 / \delta} \gamma_n}{\gamma_{n + 1}} - 1\bigr] < 0 < \liminf_{n \to \infty} \bigl[\frac{(\gamma_n)^{1 + 2 \tau \delta} - (\gamma_{n + 1})^{1 + 2 \tau \delta}}{(\gamma_n)^{1 + 2 \delta}}\bigr]
\end{equation}
and $\sum_{n = 1}^{\infty} (\gamma_n)^{1 + \tau \delta} < \infty = \sum_{n = 1}^{\infty} \gamma_n$ \cfout. Then there exist $\fC, \rho \in (0, \allowbreak \infty)$, $\vartheta \in \allowbreak \R^{\fd}$ which satisfy for all $n \in \N$ that
\begin{gather}\llabel{res1}
\textstyle (\nabla \cL)(\vartheta) = 0, \qquad \abs{\cL(\Theta_n) - \cL(\vartheta)} \le \fC \bigl[\sum_{j = 1}^{n} \gamma_j\bigr]^{- 1}, \\ \llabel{res2}
\textstyle \andqq \norm{\Theta_n - \vartheta} \le \fC \bigl[\sum_{j = 1}^{n} \gamma_j\bigr]^{- \rho} + \fC \bigl[\sum_{j = n + 1}^{\infty} (\gamma_j)^{1 + \tau \delta}\bigr].
\end{gather}
\end{athm}
\end{tcolorbox}
\end{samepage}
\end{savenotes}
\begin{aproof}
Throughout this proof let $\cG = \nabla \cL$ and for every $n \in \N$ let $\fm_n, \allowbreak \Gamma_n \allowbreak \in \R^{\fd}$, $\sigma_n \in \R$ satisfy
\begin{equation}\llabel{fm_L}
\textstyle \fm_n = (\gamma_n)^{- \delta} \bbA_n p_n (\m_n - \cG(\Theta_{n - 1})), \qquad \Gamma_n = \bbA_n \cG(\Theta_{n - 1}), \qqandqq \sigma_n = (\gamma_n)^{\delta}
\end{equation}
and let $K_n \allowbreak \subseteq \N$ and $G_n \allowbreak \subseteq \N$ satisfy
\begin{equation}\llabel{Gamma_l}
\textstyle K_n = \{m \in \N \cap [n, \infty) \colon \cG(\Theta_{m - 1}) \neq 0\} \qqandqq G_n = \{m \in \N \cap [n, \infty) \colon \Gamma_m \neq 0\}.
\end{equation}
\argument{\lref{Theta}; \lref{fm_L}}{that for all $n \in \N$ it holds that
\begin{equation}\llabel{Theta_new}
\begin{split}
\textstyle \Theta_n & \textstyle = \Theta_{n - 1} - \gamma_n (\bbA_n \cG(\Theta_{n - 1}) + [\bbA_n p_n (\m_n - \cG(\Theta_{n - 1}))]) \\
& \textstyle = \Theta_{n - 1} - \gamma_n (\Gamma_n + (\gamma_n)^{\delta} \fm_n) = \Theta_{n - 1} - \gamma_n (\Gamma_n + \sigma_n \fm_n).
\end{split}
\end{equation}
}
\argument{\lref{bound_2}; \lref{fm_L}; the fact that $0 < \tau \le 1$}{that there exists $N \in \N$ which satisfies $\limsup_{n \to \infty} \gamma_n = 0 = \limsup_{n \to \infty} \sigma_n$ and
\begin{equation}\llabel{gamma_sum}
\begin{split}
& \textstyle 2 \sum_{n = N}^{\infty} (\gamma_n)^{1 + \delta} = \sum_{n = N}^{\infty} \gamma_n [\sigma_n + (\gamma_n)^{\delta}] \\
& \textstyle \le \sum_{n = N}^{\infty} \gamma_n [\sigma_n + (\gamma_n)^{\delta}]^{\tau} = 2^{\tau} \sum_{n = N}^{\infty} (\gamma_n)^{1 + \tau \delta} < \infty.
\end{split}
\end{equation}
}
\argument{\cref{prop:gen_convergence_ohne_Loja} (applied for every $n \in \N$ with $\fd \with \fd$, $\delta \with \delta$, $\cL \with \cL$, $\bbA_n \with \bbA_n$, $\Theta_n \with \Theta_n$, $\mu_n \with \mu_n$, $p_n \with p_n$, $\alpha_n \with \alpha_n$, $\gamma_n \with \gamma_n$, $\m_n \with \m_n$ in the notation of \cref{prop:gen_convergence_ohne_Loja}); \lref{m}; \lref{Theta}; \lref{bound_2}; \lref{fm_L}; \lref{Gamma_l}; \lref{gamma_sum}; the fact that $\cG$ is H\"{o}lder continuous; the assumption that $p$ is bounded; H\"{o}lder's inequality}{that there exists $M \in \allowbreak \N \cap [N, \infty)$ which satisfies for all $k \in \N$ that
\begin{gather}\llabel{Gamma_bound}
\textstyle \sup_{n \in \N} \norm{\Gamma_n} = \sup_{n \in \N} \norm{\bbA_n \cG(\Theta_{n - 1})} \le (\sup_{n \in \N} \norm{\bbA_n}) (\sup_{n \in \N} \norm{\cG(\Theta_{n - 1})}) < \infty, \\ \llabel{fm_bound}
\textstyle \sup_{n \in \N} \norm{\fm_n} = \sup_{n \in \N} (\gamma_n)^{- \delta} \norm{\bbA_n p_n (\m_n - \cG(\Theta_{n - 1}))} < \infty, \\ \llabel{Cesaro}
\textstyle \frac{\gamma_k [\sigma_k + (\gamma_k)^{\delta}]^{2 \tau} - \gamma_{k + 1} [\sigma_{k + 1} + (\gamma_{k + 1})^{\delta}]^{2 \tau}}{\gamma_k [\sigma_k + (\gamma_k)^{\delta}]^2} = \frac{2^{2 \tau - 2} [(\gamma_k)^{1 + 2 \tau \delta} - (\gamma_{k + 1})^{1 + 2 \tau \delta}]}{(\gamma_k)^{1 + 2 \delta}}, \\ \llabel{Gamma_positiv_op}
\textstyle \andqq \spro{\Gamma_{k + M - 1}, \cG(\Theta_{k + M - 2})} \ge \norm{\cG(\Theta_{k + M - 2})}^2\dott
\end{gather}
}
\argument{\lref{bound_2}; \lref{fm_L}; \lref{Gamma_l}; \lref{Gamma_positiv_op}}{that
\begin{gather}
\textstyle \liminf_{n \to \infty} \min\Bigl\{\inf\limits_{m \in K_n} \frac{\spro{\Gamma_m, (\nabla \cL)(\Theta_{m - 1})}}{\norm{(\nabla \cL)(\Theta_{m - 1})}^{2}}, \inf\limits_{m \in G_n} \frac{\norm{(\nabla \cL)(\Theta_{m - 1})}}{\norm{\Gamma_m}}\Bigr\} > 0 \\ \llabel{l_L_Cesaro}
\textstyle \andqq \liminf_{n \to \infty} \bigl[\frac{\gamma_n [\sigma_n + (\gamma_n)^{\delta}]^{2 \tau} - \gamma_{n + 1} [\sigma_{n + 1} + (\gamma_{n + 1})^{\delta}]^{2 \tau}}{\gamma_n [\sigma_n + (\gamma_n)^{\delta}]^2}\bigr] > 0 \dott
\end{gather}
}
\argument{\cref{thm:univ_scheme_convergence_Loja} (applied for every $n \in \N$ with $\fd \with \fd$, $\delta \with \delta$, $\tau \with \tau$, $\kappa \with 2$, $\cL \with \cL$, $\Theta_n \with \Theta_n$, $\Gamma_n \with \Gamma_n$, $\mu_n \with \fm_n$, $K_n \with K_n$, $G_n \with G_n$, $\gamma_n \with \gamma_n$, $\sigma_n \with \sigma_n$ in the notation of \cref{thm:univ_scheme_convergence_Loja}); \lref{fm_L}; \lref{Gamma_l}; \lref{Theta_new}; \lref{gamma_sum}; \lref{Gamma_bound}; \lref{fm_bound}; \lref{l_L_Cesaro}; the fact that $\sum_{n = 1}^{\infty} \gamma_n = \infty$; the fact that $\frac{2 + \delta}{1 + \delta} < 2$}{that there exist $C, \eps \in (0, \infty)$, $\rho \in (0, 1)$, $\vartheta \in \R^{\fd}$ which satisfy for all $n \in \N$, $\theta \in \R^{\fd}$ with $\norm{\theta - \vartheta} < \eps$ that
\begin{gather}\llabel{Loja}
\textstyle \abs{\cL(\theta) - \cL(\vartheta)}^{\rho} \le C \norm{\cG(\theta)}, \\
\llabel{descending}
\textstyle \cL(\Theta_n) + C^{-1} \gamma_n \norm{\cG(\Theta_{n - 1})}^{2} \le \cL(\Theta_{n - 1}) + C (\gamma_n)^{1 + 2 \delta}, \\
\llabel{cL_lower_bound}
\textstyle C^{-1} (\cL(\vartheta) - \cL(\Theta_{n - 1})) \le \sum_{j = n}^{\infty} (\gamma_j)^{1 + 2 \delta} \le C (\gamma_n)^{1 + 2 \tau \delta}, \\
\llabel{Theta_upper_bound}
\textstyle \norm{\Theta_{n - 1} - \vartheta} \le C \bigl[\cL(\Theta_{n - 1}) - \cL(\vartheta) + C \sum_{j = n}^{\infty} (\gamma_j)^{1 + 2 \delta} \bigr]^{1 - \rho} + C \sum_{j = n}^{\infty} (\gamma_j)^{1 + \tau \delta}, \\
\llabel{Theta_convergence}
\textstyle \andqq \norm{(\nabla \cL)(\vartheta)} + \limsup_{n \to \infty} \norm{\Theta_n - \vartheta} = 0 \dott
\end{gather}
}
Let $\bbL \colon \N_0 \to \R$ satisfy for all $n \in \N_0$ that
\begin{equation}\llabel{bbL}
\textstyle \bbL_n = \cL(\Theta_n) - \cL(\vartheta) + C \sum_{j = n + 1}^{\infty} (\gamma_j)^{1 + 2 \delta}.
\end{equation}
\startnewargseq
\argument{\lref{descending}; \lref{cL_lower_bound}; \lref{Theta_upper_bound}; \lref{bbL}}{for all $n \in \N$ that
\begin{gather}\llabel{descend_bbL}
\textstyle \bbL_{n - 1} \ge \bbL_{n} + C^{-1} \gamma_n \norm{\cG(\Theta_{n - 1})}^2 \ge 0 \\
\llabel{Theta_upper_bound_bbL}
\textstyle \andqq \norm{\Theta_{n - 1} - \vartheta} \le C [\bbL_{n - 1}]^{1 - \rho} + C \sum_{j = n}^{\infty} (\gamma_j)^{1 + \tau \delta}.
\end{gather}
}
\argument{\lref{bound_2}; \lref{cL_lower_bound}; \lref{Theta_convergence}}{that there exists $K \in \N \cap [M, \infty)$ which satisfies for all $n \in \N \cap [K, \infty)$ that
\begin{gather}\llabel{K}
\textstyle \gamma_{n + 1} \le \gamma_n, \qquad \norm{\Theta_{n - 1} - \vartheta} + \abs{\cL(\Theta_{n - 1}) - \cL(\vartheta)} < \min\{1, \eps\}, \\
\textstyle \andqq 2 \gamma_n \Abs{\bbL_{n - 1} - C \sum_{j = n}^{\infty} (\gamma_j)^{1 + 2 \delta}} < \min\{1, C^3\} \dott
\end{gather}
}
\argument{\lref{Loja}; \lref{cL_lower_bound}; \lref{bbL}; \lref{descend_bbL}; \lref{K}; the fact that for all $a, b \in \R$ it holds that $\abs{a - b}^2 \ge \frac{1}{2} \abs{a}^2 - \abs{b}^2$}{for all $n \in \N \cap [K, \infty)$ that
\begin{equation}\llabel{bbL_telescope}
\begin{split}
& \textstyle \bbL_{n - 1} - \bbL_n \ge C^{-1} \gamma_n \norm{\cG(\Theta_{n - 1})}^2 \ge C^{-3} \gamma_n \abs{\cL(\Theta_{n - 1}) - \cL(\vartheta)}^{2 \rho} \\
& \textstyle \ge C^{-3} \gamma_n \abs{\cL(\Theta_{n - 1}) - \cL(\vartheta)}^{2} = C^{-3} \gamma_n \Abs{\bbL_{n - 1} - C \sum_{j = n}^{\infty} (\gamma_j)^{1 + 2 \delta}}^2 \\
& \textstyle \ge C^{-3} \gamma_n \bigl(\frac{1}{2} [\bbL_{n - 1}]^2 - C^2 \bigl[\sum_{j = n}^{\infty} (\gamma_j)^{1 + 2 \delta}\bigr]^2\bigr) \ge C^{-3} \gamma_n \bigl(\frac{1}{2} [\bbL_{n - 1}]^2 - C^4 \bigl[(\gamma_n)^{1 + 2 \tau \delta}\bigr]^2\bigr) \\
& \textstyle = \frac{\gamma_n}{2 C^3} [\bbL_{n - 1}]^2 - C (\gamma_n)^{3 + 4 \tau \delta} \dott
\end{split}
\end{equation}
}
\argument{\lref{descend_bbL}; \lref{bbL_telescope}}{for all $n \in \N \cap [K, \infty)$ that
\begin{equation}\llabel{main_ineq_for_rate}
\textstyle \bbL_n + \frac{\gamma_n}{2 C^3} [\bbL_{n}]^2 \le \bbL_n + \frac{\gamma_n}{2 C^3} [\bbL_{n - 1}]^2 \le \bbL_{n - 1} + C (\gamma_n)^{3 + 4 \tau \delta} \dott
\end{equation}
}
\argument{\lref{bound_2}; \cref{lr_for_rate} (applied for every $n \in \N$ with $\delta \with 2 \delta $, $\tau \with 2 \tau \delta$, $\gamma_n \with \gamma_n$ in the notation of \cref{lr_for_rate})}{that
\begin{equation}\llabel{lr_for_rate_res}
\textstyle \sup_{n \in \N} \bigl[(\gamma_n)^{1 + 2 \tau \delta} \sum_{j = 1}^{n} \gamma_j\bigr] < \infty \dott
\end{equation}
}
\argument{\lref{main_ineq_for_rate}; \lref{lr_for_rate_res}; \cref{rate_lemma} (applied for every $n \in \N \cap [K, \infty)$ with $a \with (2 C^3)^{-1}$, $b \with C$, $\tau \with 2 \tau \delta$, $\gamma_n \with \gamma_n$, $e_{n - 1} \with \bbL_{n - 1}$ in the notation of \cref{rate_lemma})}{that
\begin{equation}\llabel{bbL_rate}
\textstyle \sup_{n \in \N} \bigl[\bbL_n \sum_{j = 1}^{n} \gamma_j\bigr] < \infty \dott
\end{equation}
}
\argument{\lref{bound_2}; \lref{cL_lower_bound}; \lref{bbL}; \lref{descend_bbL}; \lref{K}; \lref{lr_for_rate_res}; \lref{bbL_rate}; the triangle inequality}{that there exists $\scrC \in \R$ which satisfies for all $n \in \N \cap [K, \infty)$ that
\begin{equation}\llabel{cL_res}
\begin{split}
& \textstyle \abs{\cL(\Theta_n) - \cL(\vartheta)} = \abs{\bbL_n - C \sum_{j = n + 1}^{\infty} (\gamma_j)^{1 + 2 \delta}} \le \bbL_n + C \sum_{j = n + 1}^{\infty} (\gamma_j)^{1 + 2 \delta} \\
& \textstyle \le \bbL_n + C^2 (\gamma_{n + 1})^{1 + 2 \tau \delta} \le \bbL_n + C^2 (\gamma_n)^{1 + 2 \tau \delta} \le \scrC \bigl[\sum_{j = 1}^{n} \gamma_j\bigr]^{-1} \dott
\end{split}
\end{equation}
}
\argument{\lref{cL_res}}{that
\begin{equation}\llabel{cL_result}
\textstyle \limsup_{n \to \infty} \bigl[\abs{\cL(\Theta_n) - \cL(\vartheta)} \sum_{j = 1}^{n} \gamma_j\bigr] < \infty \dott
\end{equation}
}
\argument{\lref{Theta_convergence}; \lref{bbL}; \lref{Theta_upper_bound_bbL}; \lref{bbL_rate}; \lref{cL_result}}[verbs=e]{\lref{res1,res2}. }
\end{aproof}

\subsection{UGD applied to KL functions with locally Lipschitz continuous gradients}
\label{subsec:convergence_UGD_KL_Lipschitz}

In the next result, \cref{theorem:gen_convergence_Loja_no_delta} below, we specialize \cref{theorem:gen_convergence_Loja} above to the situation where the regularity parameter $\delta$ is equal to $1$. In particular, in \cref{theorem:gen_convergence_Loja_no_delta} we assume that the objective function $\cL \colon \R^{\fd} \to \R$ of the considered minimization problem is locally Lipschitz continuous instead of locally $\delta$-H\"{o}lder continuous for some $\delta \in (0,1]$ as in \cref{theorem:gen_convergence_Loja}.

\cfclear
\begin{savenotes}
\begin{samepage}
\begin{tcolorbox}[colback=white!95!gray,
                  colframe=black,
                  boxrule=0.5pt,
                  sharp corners,
                  enhanced,
                 ]
\begin{athm}{cor}{theorem:gen_convergence_Loja_no_delta}
Let $\fd \in \N$, $\tau \in (0, 1]$, let $\cL \colon \R^{\fd} \to \R$ be a\cfadd{def:Loja} \KL\ function, assume that $\nabla \cL$ is locally Lipschitz continuous, let $\bbA \colon \N \to \R^{\fd \times \fd}$, $\Theta \allowbreak \colon \allowbreak \N_0 \allowbreak \to \allowbreak \R^\fd$, $\mu \colon \N \to \R^{\fd}$, and $p \colon \allowbreak \N \allowbreak \to \R$ be bounded, let $\alpha \colon \N \to [0, 1]$, $\gamma \colon \N \to (0, \infty)$, and $\m \colon \allowbreak \N_0 \allowbreak \to \R^\fd$ satisfy for all $n \in \N$ that
\begin{gather}\llabel{m}
\textstyle \m_n = \alpha_n \m_{n - 1} + (1 - \alpha_n) \bigl[(\nabla \cL)(\Theta_{n - 1}) + \gamma_n \mu_n\bigr] \\
\llabel{Theta}
\textstyle \andqq \Theta_n = \Theta_{n - 1} - \gamma_n \bbA_n \bigl[p_n \m_n + (1 - p_n) (\nabla \cL)(\Theta_{n - 1})\bigr],
\end{gather}
assume for all $n \in \N$ that $\bbA_n - \bbI_{\fd}$ is symmetric positive semi-definite, and assume
\begin{equation}\llabel{bound_2}
\textstyle \limsup_{n \to \infty} \bigl[\frac{\alpha_{n + 1} \gamma_n}{\gamma_{n + 1}} - 1\bigr] < 0 < \liminf_{n \to \infty} \bigl[\frac{(\gamma_n)^{1 + 2 \tau} - (\gamma_{n + 1})^{1 + 2 \tau}}{(\gamma_n)^3}\bigr]
\end{equation}
and $\sum_{n = 1}^{\infty} (\gamma_n)^{1 + \tau} < \infty = \sum_{n = 1}^{\infty} \gamma_n$ \cfout. Then there exist $\fC, \rho \in (0, \allowbreak \infty)$, $\vartheta \in \allowbreak \R^{\fd}$ which satisfy for all $n \in \N$ that
\begin{gather}\llabel{res1}
\textstyle (\nabla \cL)(\vartheta) = 0, \qquad \abs{\cL(\Theta_n) - \cL(\vartheta)} \le \fC \bigl[\sum_{j = 1}^{n} \gamma_j\bigr]^{- 1}, \\ \llabel{res2}
\textstyle \andqq \norm{\Theta_n - \vartheta} \le \fC \bigl[\sum_{j = 1}^{n} \gamma_j\bigr]^{- \rho} + \fC \bigl[\sum_{j = n + 1}^{\infty} (\gamma_j)^{1 + \tau}\bigr].
\end{gather}
\end{athm}
\end{tcolorbox}
\end{samepage}
\end{savenotes}
\begin{aproof}
\argument{\cref{theorem:gen_convergence_Loja} (applied for every $n \in \N$ with $\fd \with \fd$, $\delta \with 1$, $\tau \with \tau$, $\cL \with \cL$, $\bbA_n \with \bbA_n$, $\Theta_n \with \Theta_n$, $\mu_n \with \mu_n$, $p_n \with p_n$, $\alpha_n \with \alpha_n$, $\gamma_n \with \gamma_n$, $\m_n \with \m_n$ in the notation of \cref{theorem:gen_convergence_Loja}); \lref{m}; \lref{Theta}; \lref{bound_2}}[verbs=e]{\lref{res1} and \lref{res2}}.
\end{aproof}

In the following auxiliary and elementary result, \cref{lemma:l_rate}, we provide explicit examples for the sequence $(\gamma_n)_{n \in \N} \subseteq (0, \infty)$ appearing in the optimization method in \cref{theorem:gen_convergence_Loja_no_delta.m}--\cref{theorem:gen_convergence_Loja_no_delta.Theta} in \cref{theorem:gen_convergence_Loja_no_delta} that satisfy all the assumptions of \cref{theorem:gen_convergence_Loja_no_delta}. We prove the main result of this work, \cref{cor:gen_convergence_rate_Loja_l_rate} below, by an application of \cref{theorem:gen_convergence_Loja_no_delta} and in our proof of \cref{cor:gen_convergence_rate_Loja_l_rate} we employ \cref{lemma:l_rate} to verify that the assumptions of \cref{theorem:gen_convergence_Loja_no_delta} are fulfilled.

\startnewargseq
\argument{}{ that for all $\lrexpo \in (\nicefrac{3}{4}, \infty)$ it holds that $\nicefrac{3}{2} < 2 \lrexpo$.
}
\Hence for all $\lrexpo \in (\nicefrac{3}{4}, \infty)$ that $1 - \lrexpo < \lrexpo - \frac{1}{2}$. This shows that for all $\lrexpo \in (\nicefrac{3}{4}, \infty)$ it holds that $\frac{1 - \lrexpo}{\lrexpo} < 1 - \frac{1}{2 \lrexpo}$. \Hence that for all $\lrexpo \in (\nicefrac{3}{4}, \infty)$ it holds that
\begin{equation}\label{subsec:feasible_lr:existence}
\textstyle \max\{\frac{1 - \lrexpo}{\lrexpo}, - \frac{1}{2}\} < 1 - \frac{1}{2 \lrexpo} \dott
\end{equation}
\Hence that the set $(\frac{1 - \lrexpo}{\lrexpo}, 1 - \frac{1}{2 \lrexpo}] \cap (- \nicefrac{1}{2}, \infty)$ appearing below \cref{lemma:l_rate.fns} in \cref{lemma:l_rate} is not the empty set.

\cfclear
\begin{savenotes}
\begin{samepage}
\begin{tcolorbox}[colback=white!95!gray,
                  colframe=black,
                  boxrule=0.5pt,
                  sharp corners,
                  enhanced,
                 ]
\begin{athm}{lemma}{lemma:l_rate}
Let $c \in (0, \infty)$, $\lrexpo \in (\nicefrac{3}{4}, \infty)$, let $\alpha \colon \allowbreak \N \allowbreak \to [0, \allowbreak \infty)$ satisfy $\limsup_{n \to \infty} \alpha_n < 1$, and for every $n \in \N$ let $\gamma_n \in \R$ satisfy
\begin{equation}\llabel{fns}
\textstyle \gamma_n = c n^{-\lrexpo}.
\end{equation}
Then it holds for all $\tau \in (\frac{1 - \lrexpo}{\lrexpo}, 1 - \frac{1}{2 \lrexpo}] \cap (\nicefrac{- 1}{2}, \infty)$, $\eps \in (0, \infty)$ that there exists $\scrc \in \allowbreak (0, \allowbreak \infty)$ which satisfies
\begin{gather}
\textstyle \limsup_{n \to \infty} \bigl[(\alpha_{n + 1})^{\eps} \gamma_n (\gamma_{n + 1})^{-1}\bigr] < 1 = \lim_{n \to \infty} \bigl[\gamma_n (\gamma_{n + 1})^{-1}\bigr], \\
\textstyle \sum_{n = 1}^{\infty} (\gamma_n)^{1 + \tau} < \infty = \sum_{n = 1}^{\infty} (\gamma_n)^{\mathbbm{1}_{(0, 1]}(\lrexpo)}, \\
\textstyle \lim_{n \to \infty} \bigl[\frac{(\gamma_n)^{1 + 2 \tau} - (\gamma_{n + 1})^{1 + 2 \tau}}{(\gamma_n)^3}\bigr] =
\begin{cases}
\infty & \colon \tau < 1 - \frac{1}{2 \lrexpo} \\
\scrc & \colon \tau = 1 - \frac{1}{2 \lrexpo},
\end{cases} \\
\textstyle \andqq \limsup_{n \to \infty} \bigl[n^{\lrexpo (1 + \tau) - 1} \sum_{j = n + 1}^{\infty} (\gamma_j)^{1 + \tau}\bigr] < \infty \dott
\end{gather}
\end{athm}
\end{tcolorbox}
\end{samepage}
\end{savenotes}
\begin{aproof}
Throughout this proof let $\tau \in (\frac{1 - \lrexpo}{\lrexpo}, 1 - \frac{1}{2 \lrexpo}] \cap (\nicefrac{- 1}{2}, \infty)$, $\eps \in (0, \infty)$.
\argument{the fact that $\max\{\frac{1 - \lrexpo}{\lrexpo}, - \frac{1}{2}\} < \tau \le 1 - \frac{1}{2 \lrexpo}$}{that
\begin{equation}\llabel{tau}
\textstyle 2 \lrexpo (1 - \tau) \ge 1 < \lrexpo (1 + \tau) \qqandqq 1 + 2 \tau > 0 \dott
\end{equation}
}
\argument{\lref{fns}; the fact that $\limsup_{n \to \infty} \alpha_n < 1$}{that
\begin{equation}\llabel{eqn1}
\textstyle \lim_{n \to \infty} \bigl(\frac{\gamma_n}{\gamma_{n + 1}}\bigr) = \lim_{n \to \infty} \bigl(\frac{c n^{- \lrexpo}}{c (n + 1)^{- \lrexpo}}\bigr) = \lim_{n \to \infty} \bigl(1 + \frac{1}{n}\bigr)^{\lrexpo} = 1 > \limsup_{n \to \infty} (\alpha_{n + 1})^{\eps} \dott
\end{equation}}
\argument{\lref{eqn1}}[verbs=e]{that
\begin{equation}
\textstyle \limsup_{n \to \infty} \bigl[(\alpha_{n + 1})^{\eps} \gamma_n (\gamma_{n + 1})^{-1}\bigr] < 1 = \lim_{n \to \infty} \bigl[\gamma_n (\gamma_{n + 1})^{-1}\bigr] \dott
\end{equation}
}
\argument{\lref{fns}}{for all $n \in \N$ that
\begin{equation}\llabel{tau_for_summability}
\textstyle (\gamma_n)^{1 + \tau} = c^{1 + \tau} n^{- \lrexpo (1 + \tau)}.
\end{equation}
}
\argument{\lref{tau}; \lref{tau_for_summability}}{that
\begin{equation}
\textstyle \sum_{n = 1}^{\infty} (\gamma_n)^{1 + \tau} < \infty = \sum_{n = 1}^{\infty} (\gamma_n)^{\mathbbm{1}_{(0, 1]}(\lrexpo)} \dott
\end{equation}
}
\argument{\lref{fns}}{for all $n \in \N$ that
\begin{equation}\llabel{gamma_sigma_to_n}
\begin{split}
& \textstyle \frac{(\gamma_n)^{1 + 2 \tau} - (\gamma_{n + 1})^{1 + 2 \tau}}{(\gamma_n)^3} = \frac{c^{1 + 2 \tau}}{c^3} \cdot \frac{n^{- \lrexpo (1 + 2 \tau)} - (n + 1)^{- \lrexpo (1 + 2 \tau)}}{n^{- 3 \lrexpo}} \\
& \textstyle = c^{2 \tau - 2} \cdot \frac{1 - (1 + \frac{1}{n})^{- \lrexpo (1 + 2 \tau)}}{n^{- \lrexpo(2 - 2 \tau)}} = c^{2 \tau - 2} \cdot \frac{1 - (1 + \frac{1}{n})^{- \lrexpo(1 + 2 \tau)}}{(\frac{1}{n})^{2 \lrexpo (1 - \tau)}}.
\end{split}
\end{equation}
}
\argument{\lref{tau}; L'H\^{o}pital's rule}{that
\begin{equation}\llabel{Lhopital}
\begin{split}
& \textstyle \lim\limits_{n \to \infty} \frac{1 - (1 + \frac{1}{n})^{- \lrexpo(1 + 2 \tau)}}{(\frac{1}{n})^{2 \lrexpo (1 - \tau)}} = \lim\limits_{x \to 0} \frac{1 - (1 + x)^{- \lrexpo(1 + 2 \tau)}}{x^{2 \lrexpo (1 - \tau)}} = \lim\limits_{x \to 0} \frac{\lrexpo (1 + 2 \tau) (1 + x)^{- \lrexpo (1 + 2 \tau) - 1}}{2 \lrexpo (1 - \tau) x^{2 \lrexpo (1 - \tau) - 1}} \\
& \textstyle = \frac{1 + 2 \tau}{2 (1 - \tau)} \lim\limits_{x \to 0} \frac{1}{x^{2 \lrexpo (1 - \tau) - 1}} =
\begin{cases}
\infty & \colon 2 \lrexpo (1 - \tau) > 1 \\
\frac{1 + 2 \tau}{2 (1 - \tau)} & \colon 2 \lrexpo (1 - \tau) = 1 \dott
\end{cases}
\end{split}
\end{equation}
}
\argument{\lref{tau}; \lref{gamma_sigma_to_n}; \lref{Lhopital}}{that there exists $\scrc \in (0, \infty)$ which satisfies
\begin{equation}
\textstyle \lim_{n \to \infty} \bigl[\frac{(\gamma_n)^{1 + 2 \tau}  - (\gamma_{n + 1})^{1 + 2 \tau}}{(\gamma_n)^3}\bigr] =
\begin{cases}
\infty &\!\!\! \colon \tau < 1 - \frac{1}{2 \lrexpo} \\
\scrc &\!\!\! \colon \tau = 1 - \frac{1}{2 \lrexpo} \dott
\end{cases}
\end{equation}
}
\argument{\lref{tau}; L'H\^{o}pital's rule}{that
\begin{equation}\llabel{Lhopital_last}
\begin{split}
& \textstyle \lim\limits_{n \to \infty} \frac{n^{- \lrexpo (1 + \tau)}}{(n - 1)^{1 - \lrexpo (1 + \tau)} - n^{1 - \lrexpo (1 + \tau)}} = \lim\limits_{n \to \infty} \frac{\frac{1}{n}}{(1 - \frac{1}{n})^{1 - \lrexpo (1 + \tau)} - 1} = \lim\limits_{x \to 0} \frac{x}{(1 - x)^{1 - \lrexpo (1 + \tau)} - 1} \\
& \textstyle = \lim\limits_{x \to 0} \frac{1}{(\lrexpo (1 + \tau) - 1)(1 - x)^{- \lrexpo (1 + \tau)}} = \frac{1}{\lrexpo (1 + \tau) - 1} \dott
\end{split}
\end{equation}
}
\argument{\lref{Lhopital_last}; the fact that $\lim_{n \to \infty} n^{1 - \lrexpo (1 + \tau)} = 0 = \lim_{n \to \infty} \sum_{j = n + 1}^{\infty} j^{- \lrexpo (1 + \tau)}$; the fact that $(n^{1 - \lrexpo (1 + \tau)})_{n \in \N} \allowbreak \subseteq \R$ is strictly decreasing; Stolz-Ces\'{a}ro theorem}{that
\begin{equation}
\begin{split}
& \textstyle \lim\limits_{n \to \infty} \frac{\sum_{j = n + 1}^{\infty} j^{- \lrexpo (1 + \tau)}}{n^{1 - \lrexpo (1 + \tau)}} = \lim\limits_{n \to \infty} \frac{\sum_{j = n}^{\infty} j^{- \lrexpo (1 + \tau)} - \sum_{j = n + 1}^{\infty} j^{- \lrexpo (1 + \tau)}}{(n - 1)^{1 - \lrexpo (1 + \tau)} - n^{1 - \lrexpo (1 + \tau)}} \\
& \textstyle = \lim\limits_{n \to \infty} \frac{n^{- \lrexpo (1 + \tau)}}{(n - 1)^{1 - \lrexpo (1 + \tau)} - n^{1 - \lrexpo (1 + \tau)}} = \frac{1}{\lrexpo (1 + \tau) - 1} \dott
\end{split}
\end{equation}
}
\end{aproof}

In \cref{cor:gen_convergence_rate_Loja_l_rate.result} in \cref{cor:gen_convergence_rate_Loja_l_rate} below we bound weak and scaled strong approximation errors of the considered optimization process from above by a constant multiplied by the quantity $[\sum_{j = 1}^n j^{- \lrexpo}]^{- 1}$. This quantity can, of course, also be expressed in terms of a polynomial convergence rate depending on $\lrexpo$. Specifically, in the following elementary fact in \cref{lemma:l_rate_specific} we recall optimal lower and upper bounds for this quantity. \cref{lemma:l_rate_specific} is then also used in our proof of \cref{cor:gen_convergence_rate_Loja_l_rate}.

\cfclear
\begin{savenotes}
\begin{samepage}
\begin{tcolorbox}[colback=white!95!gray,
                  colframe=black,
                  boxrule=0.5pt,
                  sharp corners,
                  enhanced,
                 ]
\begin{athm}{lemma}{lemma:l_rate_specific}
Let $\lrexpo \in (- \infty, 1]$. Then there exists $\rho \in (0, \infty)$ which satisfies for all $n \in \N$ that
\begin{equation}\llabel{result}
\textstyle \rho^{-1} [n^{1 - \lrexpo} + \ln(n)]^{- 1} \le \bigl[\sum_{j = 1}^n j^{- \lrexpo}\bigr]^{-1} \le \rho [n^{1 - \lrexpo} + \ln(n)]^{- 1} \dott
\end{equation}
\end{athm}
\end{tcolorbox}
\end{samepage}
\end{savenotes}
\begin{aproof}
\argument{}{ that for all $n \in \N$, $\alpha \in [0, \infty)$, $\beta \in (- \infty, 0]$ it holds that
\begin{gather}
\textstyle \int_n^{n + 1} x^{- \alpha} \d x \le \int_n^{n + 1} n^{- \alpha} \d x = n^{ - \alpha} = \int_{n - 1}^n n^{ - \alpha} \d x \le \int_{n - 1}^n x^{- \alpha} \d x \\ \llabel{bound_harmonic}
\textstyle \andqq \int_{n - 1}^n x^{- \beta} \d x \le \int_{n - 1}^n n^{- \beta} \d x = n^{- \beta} = \int_n^{n + 1} n^{- \beta} \d x \le \int_n^{n + 1} x^{- \beta} \d x \dott
\end{gather}
}
\argument{\lref{bound_harmonic}}{for all $n \in \N$, $\alpha \in [0, \infty)$, $\beta \in (- \infty, 0]$ that
\begin{equation}\llabel{int_bound}
\textstyle \int_1^{n + 1} x^{- \alpha} \d x \le \sum_{j = 1}^n j^{- \alpha} \le 1 + \int_1^n x^{- \alpha} \d x \qandq 1 + \int_1^n x^{- \beta} \d x \le \sum_{j = 1}^n j^{- \beta} \le \int_1^{n + 1} x^{- \beta} \d x \dott
\end{equation}
}
\argument{\lref{int_bound}}{for all $n \in \N$, $\alpha \in [0, \infty) \backslash \{1\}$, $\beta \in (- \infty, 0]$ that
\begin{gather}
\textstyle \frac{(n + 1)^{1 - \alpha} - 1}{1 - \alpha} \le \sum_{j = 1}^n j^{- \alpha} \le 1 + \frac{n^{1 - \alpha} - 1}{1 - \alpha}, \qquad 1 + \frac{n^{1 - \beta} - 1}{1 - \beta} \le \sum_{j = 1}^n j^{- \beta} \le \frac{(n + 1)^{1 - \beta} - 1}{1 - \beta}, \\ \llabel{ln_case}
\textstyle \andqq \ln(n + 1) \le \sum_{j = 1}^n j^{- 1} \le 1 + \ln(n) \dott
\end{gather}
}
\argument{\lref{ln_case}; the fact that for all $n \in \N$ it holds that $1 + \ln(n) \le 2 \ln(n + 1)$; the fact that for all $\alpha \allowbreak \in [0, \allowbreak 1)$, $\beta \in (- \infty, \allowbreak 0]$ there exists $\rho \in (0, \infty)$ which satisfies for all $n \in \N$ that
\begin{gather}
\textstyle \rho^{- 1} [n^{1 - \alpha} + \ln(n)] \le \frac{(n + 1)^{1 - \alpha} - 1}{1 - \alpha}, \quad 1 + \frac{n^{1 - \alpha} - 1}{1 - \alpha} \le \rho [n^{1 - \alpha} + \ln(n)], \\
\textstyle \rho^{- 1} [n^{1 - \beta} + \ln(n)] \le 1 + \frac{n^{1 - \beta} - 1}{1 - \beta}, \qqandqq \frac{(n + 1)^{1 - \beta} - 1}{1 - \beta} \le \rho [n^{1 - \beta} + \ln(n)]
\end{gather}
}[verbs=e]{\lref{result}.
}
\end{aproof}

In the next step we combine \cref{theorem:gen_convergence_Loja_no_delta}, \cref{lemma:l_rate}, and \cref{lemma:l_rate_specific} to establish in the following result, \cref{cor:gen_convergence_rate_Loja_l_rate}, the main result of this work. \cref{thm_intro_strong_convergence_1} in the introduction is a direct consequence of \cref{cor:gen_convergence_rate_Loja_l_rate}. In \cref{section:UGD_applications} we illustrate the generality of \cref{cor:gen_convergence_rate_Loja_l_rate} by applying it and its consequence in \cref{cor:gen_convergence_rate_Loja_gen_l_rate_Lipschitz_alpha} below, respectively, to a number of optimization methods.

\cfclear
\begin{savenotes}
\begin{samepage}
\begin{tcolorbox}[colback=white!95!gray,
                  colframe=black,
                  boxrule=0.5pt,
                  sharp corners,
                  enhanced,
                 ]
\begin{athm}{theorem}{cor:gen_convergence_rate_Loja_l_rate}
Let $\fd \in \N$, $\lrexpo \in (\nicefrac{3}{4}, 1]$, let $\cL \colon \allowbreak \R^{\fd} \allowbreak \to \R$ be a\cfadd{def:Loja} \KL\ function, assume that $\nabla \cL$ is locally Lipschitz continuous, let $\bbA \colon \N \to \R^{\fd \times \fd}$, $\Theta \allowbreak \colon \allowbreak \N_0 \allowbreak \to \allowbreak \R^\fd$, $\mu \colon \N \to \R^{\fd}$, and $p \colon \N \to \R$ be bounded, let $\alpha \colon \N \to [0, 1]$ and $\m \allowbreak \colon \allowbreak \N_0 \allowbreak \to \R^\fd$ satisfy for all $n \in \N$ that
\begin{gather}\llabel{m}
\textstyle \m_n = \alpha_n \m_{n - 1} + (1 - \alpha_n) \bigl[(\nabla \cL)(\Theta_{n - 1}) + n^{- \lrexpo} \mu_n\bigr], \\
\llabel{Theta}
\textstyle \Theta_n = \Theta_{n - 1} - n^{- \lrexpo} \bbA_n \bigl[p_n \m_n + (1 - p_n) (\nabla \cL)(\Theta_{n - 1})\bigr],
\end{gather}
and $\limsup_{k \to \infty} \alpha_k < 1$, and assume for all $n \in \N$ that $\bbA_n - \bbI_{\fd}$ is symmetric positive semi-definite \cfout. Then there exist $\vartheta \allowbreak \in \allowbreak \R^{\fd}$, $\rho \in (0, \infty)$ which satisfy for all $n \in \N$ that
\begin{equation}\llabel{result}
\textstyle \norm{(\nabla \cL) (\vartheta)} + \abs{\cL(\Theta_n) - \cL(\vartheta)} + \norm{\Theta_n - \vartheta}^{\rho} \le \rho \bigl[\sum_{j = 1}^{n} j^{- \lrexpo}\bigr]^{-1} \dott
\end{equation}
\end{athm}
\end{tcolorbox}
\end{samepage}
\end{savenotes}
\begin{aproof}
Throughout this proof let $\tau = 1 - \frac{1}{2 \lrexpo}$ and for every $n \in \N$ let $\gamma_n \in \R$ satisfy $\gamma_n = n^{- \lrexpo}$.
\argument{\lref{m}; \lref{Theta}; the fact that for all $n \in \N$ it holds that $\gamma_n = n^{- \lrexpo}$}{for all $n \in \N$ that
\begin{gather}\llabel{m_new}
\textstyle \m_n = \alpha_n \m_{n - 1} + (1 - \alpha_n) \bigl[(\nabla \cL)(\Theta_{n - 1}) + \gamma_n \mu_n\bigr] \\
\llabel{Theta_new}
\textstyle \andqq \Theta_{n} = \Theta_{n - 1} - \gamma_n \bbA_n \bigl[p_n \m_n + (1 - p_n) (\nabla \cL)(\Theta_{n - 1})\bigr] \dott
\end{gather}
}
\argument{the fact that $\frac{3}{4} < \lrexpo \le 1$}{that
\begin{equation}\llabel{tau}
\textstyle \tau \in (\frac{1 - \lrexpo}{\lrexpo}, 1 - \frac{1}{2 \lrexpo}] \cap (- \frac{1}{2}, \infty) \dott
\end{equation}
}
\argument{\lref{tau}; \cref{lemma:l_rate} (applied for every $n \in \N$ with $c \with 1$, $\lrexpo \allowbreak \with \lrexpo$, $\alpha_n \allowbreak \with \alpha_n$, $\gamma_n \allowbreak \with \gamma_n$ in the notation of \cref{lemma:l_rate}); the fact that $\frac{3}{4} < \lrexpo \le 1$}{that
\begin{gather}\llabel{bound_1}
\textstyle \bigl[\sum_{n = 1}^{\infty} (\gamma_n)^{1 + \tau}\bigr] + \limsup_{n \to \infty} \bigl[n^{\lrexpo (1 + \tau) - 1} \sum_{j = n + 1}^{\infty} j^{- \lrexpo (1 + \tau)}\bigr] < \infty = \sum_{n = 1}^{\infty} \gamma_n \\
\llabel{bound_2}
\textstyle \textstyle \andqq \limsup_{n \to \infty} \bigl[\frac{\alpha_{n + 1} \gamma_n}{\gamma_{n + 1}}\bigr] - 1< 0 < \liminf_{n \to \infty} \bigl[\frac{(\gamma_n)^{1 + 2 \tau} - (\gamma_{n + 1})^{1 + 2 \tau}}{(\gamma_n)^3}\bigr] \dott
\end{gather}
}
\argument{\cref{theorem:gen_convergence_Loja_no_delta} (applied for every $n \in \N$ with $\fd \with \fd$, $\tau \with \tau$, $\cL \with \cL$, $\bbA_n \with \bbA_n$, $\Theta_n \with \Theta_n$, $\mu_n \with \mu_n$, $p_n \with p_n$, $\alpha_n \with \alpha_n$, $\gamma_n \with \gamma_n$, $\m_n \with \m_n$ in the notation of \cref{theorem:gen_convergence_Loja_no_delta}); \lref{m_new}; \lref{Theta_new}; \lref{bound_2}}{that there exist $c, \rho \in (0, \allowbreak \infty)$, $\vartheta \allowbreak \in \allowbreak \R^{\fd}$ which satisfy for all $n \in \N$ that
\begin{gather}\llabel{cL_Theta_rate}
\textstyle (\nabla \cL)(\vartheta) = 0, \qquad \abs{\cL(\Theta_n) - \cL(\vartheta)} \le c \bigl[\sum_{j = 1}^{n} j^{- \lrexpo}\bigr]^{- 1}, \\ \llabel{pre_res}
\textstyle \andqq \norm{\Theta_n - \vartheta} \le c \bigl[\sum_{j = 1}^{n} j^{- \lrexpo}\bigr]^{- \rho} + c \sum_{j = n + 1}^{\infty} j^{- \lrexpo (1 + \tau)} \dott
\end{gather}
}
\argument{\lref{bound_1}; \lref{pre_res}; the fact that
\begin{equation}
\textstyle 1 - \lrexpo (1 + \tau) = 1 - \lrexpo (1 + 1 - \frac{1}{2\lrexpo}) = \frac{3}{2} - 2 \lrexpo
\end{equation}
}{that there exists $\fc \in (c, \infty)$ which satisfies
\begin{equation}\llabel{Theta_rate_spec}
\textstyle \norm{\Theta_n - \vartheta} \le \fc \bigl[\sum_{j = 1}^{n} j^{- \lrexpo}\bigr]^{- \rho} + \fc n^{1 - \lrexpo (1 + \tau)} = \fc \bigl[\sum_{j = 1}^{n} j^{- \lrexpo}\bigr]^{- \rho} + \fc n^{(\nicefrac{3}{2}) - 2 \lrexpo} \dott
\end{equation}
}
\argument{\cref{lemma:l_rate_specific} (applied with $\lrexpo \with \lrexpo$ in the notation of \cref{lemma:l_rate_specific}); \lref{cL_Theta_rate}; \lref{Theta_rate_spec}; the assumption that $\nicefrac{3}{4} < \lrexpo \le 1$}{that there exists $\scrc \in (0, \infty)$ which satisfies for all $n \in \N$ that
\begin{equation}\llabel{pre_resss}
\textstyle \norm{(\nabla \cL) (\vartheta)} + \abs{\cL(\Theta_n) - \cL(\vartheta)} + \norm{\Theta_n - \vartheta}^{\scrc} \le \scrc \bigl[\sum_{j = 1}^{n} j^{- \lrexpo}\bigr]^{- 1} \dott
\end{equation}
}
\argument{\lref{pre_resss}}[verbs=e]{\lref{result}. }
\end{aproof}

In the next result, \cref{cor:gen_convergence_rate_Loja_l_rate_TEMP}, we provide a slight reformulation of \cref{cor:gen_convergence_rate_Loja_l_rate} which makes it a bit easier to apply the convergence result to the concrete optimization methods \ref{intro:momentum}--\ref{intro:Yogi} from \cref{sec:introduction}.

\cfclear
\begin{savenotes}
\begin{samepage}
\begin{tcolorbox}[colback=white!95!gray,
                  colframe=black,
                  boxrule=0.5pt,
                  sharp corners,
                  enhanced,
                 ]
\begin{athm}{cor}{cor:gen_convergence_rate_Loja_l_rate_TEMP}
Let $\fd \in \N$, $\lrexpo \in (\nicefrac{3}{4}, 1]$, let $\cL \colon \allowbreak \R^{\fd} \allowbreak \to \R$ be a\cfadd{def:Loja} \KL\ function, assume that $\nabla \cL$ is locally Lipschitz continuous,  let $\alpha \colon \N \to [0, 1]$, $p \colon \N \to \R$, $\bbA \colon \N \to \R^{\fd \times \fd}$, $\Theta \allowbreak \colon \allowbreak \N_0 \allowbreak \to \allowbreak \R^\fd$, $\m \allowbreak \colon \allowbreak \N_0 \allowbreak \to \R^\fd$, and $\mu \colon \N \to \R^{\fd}$ satisfy for all $n \in \N$ that
\begin{gather}\llabel{m}
\textstyle \m_n = \alpha_n \m_{n - 1} + (1 - \alpha_n) \bigl[(\nabla \cL)(\Theta_{n - 1}) + \mu_n\bigr], \\
\llabel{Theta}
\textstyle \Theta_n = \Theta_{n - 1} - \bbA_n \bigl[p_n \m_n + (1 - p_n) (\nabla \cL)(\Theta_{n - 1})\bigr],
\end{gather}
and $\limsup_{k \to \infty} \alpha_k < 1 \le 1 + \sup_{k \in \N} (\norm{\Theta_k} + \abs{p_k} + k^{\lrexpo} (\norm{\bbA_k} + \norm{\mu_k})) < \infty$, and assume for all $n \in \N$ that $n^{\lrexpo} \bbA_n - \bbI_{\fd}$ is symmetric positive semi-definite \cfout. Then there exist $\vartheta \allowbreak \in \allowbreak \R^{\fd}$, $\rho \in (0, \infty)$ which satisfy for all $n \in \N$ that
\begin{equation}\llabel{result}
\textstyle \norm{(\nabla \cL) (\vartheta)} + \abs{\cL(\Theta_n) - \cL(\vartheta)} + \norm{\Theta_n - \vartheta}^{\rho} \le \rho \bigl[\sum_{j = 1}^{n} j^{- \lrexpo}\bigr]^{-1} \dott
\end{equation}
\end{athm}
\end{tcolorbox}
\end{samepage}
\end{savenotes}
\begin{aproof}
\argument{\cref{cor:gen_convergence_rate_Loja_l_rate} (applied for every $n \in \N$ with $\fd \with \fd$, $\lrexpo \with \lrexpo$, $\cL \with \cL$, $\bbA_n \with n^{\lrexpo} \bbA_n$, $\Theta_n \with \Theta_n$, $\mu_n \with n^{\lrexpo} \mu_n$, $p_n \with p_n$, $\alpha_n \with \alpha_n$, $\m_n \with \m_n$ in the notation of \cref{cor:gen_convergence_rate_Loja_l_rate})}[verbs=e]{\lref{result}. }
\end{aproof}

In the next result, \cref{cor:gen_convergence_rate_Loja_gen_l_rate_Lipschitz_alpha}, we specialize \cref{cor:gen_convergence_rate_Loja_l_rate_TEMP} above to the situation where the momentum decay factors $(\alpha_n)_{n \in \N} \subseteq [0, 1]$ in \cref{cor:gen_convergence_rate_Loja_l_rate_TEMP.m} do not depend on the number of gradient steps $n \in \N$. In \cref{section:UGD_applications} we apply \cref{cor:gen_convergence_rate_Loja_gen_l_rate_Lipschitz_alpha} to a number of optimization methods including each of the optimization methods \ref{intro:momentum}--\ref{intro:Yogi} from \cref{sec:introduction}.

\cfclear
\begin{savenotes}
\begin{samepage}
\begin{tcolorbox}[colback=white!95!gray,
                  colframe=black,
                  boxrule=0.5pt,
                  sharp corners,
                  enhanced,
                 ]
\begin{athm}{cor}{cor:gen_convergence_rate_Loja_gen_l_rate_Lipschitz_alpha}
Let $\fd \in \N$, $\alpha \in [0, 1)$, $\lambda \in (0, \infty)$, $\lrexpo \in (\nicefrac{3}{4}, 1]$, let $\cL \colon \allowbreak \R^{\fd} \allowbreak \to \R$ be a\cfadd{def:Loja} \KL\ func\-tion, assume that $\nabla \cL$ is locally Lipschitz continuous, let $\Theta \allowbreak \colon \allowbreak \N_0 \allowbreak \to \allowbreak \R^\fd$ and $p \colon \N \to \R$ be bounded, let $\gamma \colon \N \allowbreak \to \R$, $\bbA \colon \N \to \R^{\fd \times \fd}$, $\m \allowbreak \colon \allowbreak \N_0 \allowbreak \to \R^{\fd}$, and $\mu \colon \N \allowbreak \to \R^{\fd}$ satisfy for all $n \in \N$ that
\begin{gather}\llabel{m}
\textstyle \m_n = \alpha \m_{n - 1} + (1 - \alpha) \bigl[(\nabla \cL)(\Theta_{n - 1}) + \gamma_n \mu_n\bigr], \\
\llabel{Theta}
\textstyle \Theta_n = \Theta_{n - 1} - \gamma_n \bbA_n \bigl[p_n \m_n + (1 - p_n) (\nabla \cL)(\Theta_{n - 1})\bigr],
\end{gather}
and $\sup_{k \in \N} (k^{\lrexpo} \abs{\gamma_k} (\norm{\bbA_k} + \norm{\mu_k})) < \infty$, and assume for all $n \in \N$ that $\gamma_n n^{\lrexpo} \bbA_n - \lambda \bbI_{\fd}$ is symmetric positive semi-definite \cfout. Then there exist $\vartheta \allowbreak \in \allowbreak \R^{\fd}$, $\rho \in (0, \allowbreak \infty)$ which satisfy for all $n \in \N$ that
\begin{equation}\llabel{res}
\textstyle \norm{(\nabla \cL) (\vartheta)} + \abs{\cL(\Theta_n) - \cL(\vartheta)} + \norm{\Theta_n - \vartheta}^{\rho} \le \rho \bigl[\sum_{j = 1}^{n} j^{- \lrexpo}\bigr]^{-1} \dott
\end{equation}
\end{athm}
\end{tcolorbox}
\end{samepage}
\end{savenotes}
\begin{aproof}
\argument{\cref{cor:gen_convergence_rate_Loja_l_rate} (applied for every $n \in \N$ with $\fd \with \fd$, $\lrexpo \with \lrexpo$, $\cL \with \lambda \cL$, $\bbA_n \with \lambda^{- 1} \gamma_n n^{\lrexpo} \bbA_n$, $\Theta_n \with \Theta_n$, $\mu_n \with \lambda \gamma_n n^{\lrexpo} \mu_n$, $p_n \with p_n$, $\alpha_n \with \alpha$, $\m_n \with \lambda \m_n$ in the notation of \cref{cor:gen_convergence_rate_Loja_l_rate})}[verbs=e]{\lref{res}. }
\end{aproof}

In implementations a popular choice for the learning rate $(\gamma_n)_{n \in \N} \subseteq \R$ in \cref{cor:gen_convergence_rate_Loja_gen_l_rate_Lipschitz_alpha} is the learning rate schedule $\frac{a}{b + n} \in (0, \infty)$, $n \in \N$, for appropriate $a \in (0, \infty)$, $b \in [0, \infty)$ (see, \eg, \cite[Subsection~3.2.1]{DereichJentzenRiekert_adaptive_lr}, \cite[Subsection~3.1]{WangSAdam2020}, and \cite[Section~III]{MR4637063}). In the elementary fact in \cref{lemma:common_l_rate} we briefly demonstrate that this specific learning rate schedule (corresponding to the choice $\lrexpo_1 = \lrexpo_2 = 1$ in \cref{lemma:common_l_rate}) satisfies the condition/property appearing in \cref{cor:gen_convergence_rate_Loja_gen_l_rate_Lipschitz_alpha} (with the choice $\lrexpo = 1$ in \cref{cor:gen_convergence_rate_Loja_gen_l_rate_Lipschitz_alpha}) that $\sup_{k \in \N} (k^{\lrexpo} \abs{\gamma_k}) < \infty$.

\cfclear
\begin{savenotes}
\begin{samepage}
\begin{tcolorbox}[colback=white!95!gray,
                  colframe=black,
                  boxrule=0.5pt,
                  sharp corners,
                  enhanced,
                 ]
\begin{athm}{lemma}{lemma:common_l_rate}
Let $a \in \R \backslash \{0\}$, $b, \lrexpo_1, \lrexpo_2 \in [0, \infty)$ and let $\gamma \colon \N \to \R$ satisfy for all $n \in \N$ that $\gamma_n = \frac{a}{(b + n^{\lrexpo_1})^{\lrexpo_2}}$. Then
\begin{enumerate}[label=(\roman*)]
\item\llabel{item1} it holds for all $n \in \N$ that $0 < \abs{a} (b + 1)^{- \lrexpo_2} \le \abs{\gamma_n} n^{\lrexpo_1 \lrexpo_2} \le \abs{a}$ and

\item\llabel{item2} it holds that $\sup_{n \in \N} (\sum_{k = - 1}^1 \abs{\gamma_n n^{\lrexpo_1 \lrexpo_2}}^k) < \infty$.
\end{enumerate}
\end{athm}
\end{tcolorbox}
\end{samepage}
\end{savenotes}
\begin{aproof}
\argument{the fact that $a \neq 0$; the fact that for all $n \in \N$ it holds that
\begin{equation}
\textstyle n^{\lrexpo_1 \lrexpo_2} \le (b + n^{\lrexpo_1})^{\lrexpo_2} \le (b + 1)^{\lrexpo_2} n^{\lrexpo_1 \lrexpo_2}
\end{equation}
}{that $0 < \abs{a} (b + 1)^{- \lrexpo_2} \le \abs{\gamma_n} n^{\lrexpo_1 \lrexpo_2} \le \abs{a}$. \llabel{item1_proof}}
\argument{\lref{item1_proof}}[verbs=e]{\lref{item1}. }%
\argument{\lref{item1}}[verbs=e]{\lref{item2}. }
\end{aproof}

\subsection{UGD in the training of DNNs}
\label{subsec:dnns}

In this subsection we show in \cref{cor:dnns_1} that the assumptions on the objective function $\cL \colon \R^{\fd} \to \R$ in \cref{cor:gen_convergence_rate_Loja_gen_l_rate_Lipschitz_alpha} are general enough to cover the training of fully-connected feedforward \DNNs\ with an analytic activation function. In \cref{cor:dnns_2} we specialize \cref{cor:dnns_1} to the situation where the activation function for the \DNNs\ is nothing else but the softplus activation $\R \ni x \mapsto \ln(1 + \exp(x)) \in \R$.

\cfclear
\begin{savenotes}
\begin{samepage}
\begin{tcolorbox}[colback=white!95!gray,
                  colframe=black,
                  boxrule=0.5pt,
                  sharp corners,
                  enhanced,
                 ]
\begin{athm}{cor}{cor:dnns_1}[\resname{\DNNs\ with general activations}]
Let $\fd, K, L \in \N$, $\ell_0, \allowbreak \ell_1, \allowbreak \dots, \allowbreak \ell_L \allowbreak \in \N$, $\fx_1, \allowbreak \fx_2, \allowbreak \dots, \allowbreak \fx_K \allowbreak \in \R^{\ell_0}$ satisfy $\fd = \allowbreak \sum_{i = 1}^{L} \ell_i (\ell_{i - 1} + 1)$, let $\act \colon \R \to \R$ be analytic, for every $k \in \N$ let $H_k \colon \R^{\ell_L} \times \R^{\fd} \to \R$ be analytic, for every $\theta = \allowbreak (\theta_1, \allowbreak \dots, \allowbreak \theta_{\fd}) \allowbreak \in \R^{\fd}$ let $\cN^{k, \theta} = \allowbreak (\cN^{k, \theta}_1, \allowbreak \dots, \allowbreak \cN^{k, \theta}_{\ell_k}) \colon \allowbreak \R^{\ell_0} \allowbreak \to \R^{\ell_k}$, $k \in \allowbreak \{0, \allowbreak 1, \allowbreak \dots, \allowbreak L\}$, satisfy for all $k \in \allowbreak \{0, \allowbreak 1, \allowbreak \dots, \allowbreak L - 1\}$, $x = (x_1, \allowbreak \dots, \allowbreak x_{\ell_0}) \in \allowbreak \R^{\ell_0}$, $i \in \allowbreak \{1, \allowbreak 2, \allowbreak \dots, \allowbreak \ell_{k + 1}\}$ that
\begin{multline}\llabel{cN}
\textstyle \cN^{k + 1, \theta}_i(x) = \theta_{\ell_{k + 1} \ell_{k} + i + \sum_{h = 1}^{k} \ell_h (\ell_{h - 1} + 1)} \\
\textstyle + \sum_{j = 1}^{\ell_{k}} \theta_{(i - 1) \ell_{k} + j + \sum_{h = 1}^{k} \ell_h (\ell_{h - 1} + 1)} \bigl[x_j \mathbbm{1}_{\{0\}}(k) + \act(\cN^{k, \theta}_j(x)) \mathbbm{1}_{\N}(k)\bigr],
\end{multline}
let $\cL \colon \R^{\fd} \to \R$ satisfy for all $\theta \in \R^{\fd}$ that $\cL(\theta) = \sum_{k = 1}^{K} H_k (\cN^{L, \theta}(\fx_k), \theta)$, let $\alpha \in [0, \allowbreak 1)$, $\lrexpo \in (\nicefrac{3}{4}, \allowbreak 1]$, let $\Theta \allowbreak \colon \allowbreak \N_0 \allowbreak \to \allowbreak \R^\fd$ and $p \colon \N \to \R$ be bounded, let $\gamma \colon \N \allowbreak \to \R$, $\bbA \colon \N \to \R^{\fd \times \fd}$, $\m \allowbreak \colon \allowbreak \N_0 \allowbreak \to \R^\fd$, and $\mu \colon \N \to \R^{\fd}$ satisfy for all $n \in \N$ that
\begin{gather}\llabel{m}
\textstyle \m_n = \alpha \m_{n - 1} + (1 - \alpha) \bigl[(\nabla \cL)(\Theta_{n - 1}) + \gamma_n \mu_n\bigr], \\
\llabel{Theta}
\textstyle \Theta_n = \Theta_{n - 1} - \gamma_n \bbA_n \bigl[p_n \m_n + (1 - p_n) (\nabla \cL)(\Theta_{n - 1})\bigr],
\end{gather}
and $\sup_{k \in \N} (k^{\lrexpo} \abs{\gamma_k} (\norm{\bbA_k} + \norm{\mu_k})) < \infty$, and assume for all $n \in \N$ that $\gamma_n n^{\lrexpo} \bbA_n - \bbI_{\fd}$ is symmetric positive semi-definite \cfout. Then there exist $\vartheta \allowbreak \in \allowbreak \R^{\fd}$, $\rho \allowbreak \in (0, \allowbreak \infty)$ which satisfy for all $n \in \N$ that
\begin{equation}\llabel{result}
\textstyle \norm{(\nabla \cL) (\vartheta)} + \abs{\cL(\Theta_n) - \cL(\vartheta)} + \norm{\Theta_n - \vartheta}^{\rho} \le \rho \bigl[\sum_{j = 1}^{n} j^{- \lrexpo}\bigr]^{-1} \dott
\end{equation}
\end{athm}
\end{tcolorbox}
\end{samepage}
\end{savenotes}
\begin{aproof}
\argument{the assumption that $\act$ and $(H_k)_{k \in \N}$ are analytic; \lref{cN}}{that $\cL$ is analytic (cf., \eg, \cite[Corollary~9.14.5 in Subsection~9.14]{JentzenBookDeepLearning2023}). \llabel{cL_analytic}}
\argument{\lref{cL_analytic}}{that $\cL$ is a \KL\ function (cf., \eg, \cite{MR2274510} and the references therein). }
\argument{the fact that $\cL$ is analytic}{that $\nabla \cL$ is locally Lipschitz continuous. \llabel{cL_KL}}
\argument{\cref{cor:gen_convergence_rate_Loja_gen_l_rate_Lipschitz_alpha} (applied for every $n \in \N$ with $\fd \with \fd$, $\alpha \with \alpha$, $\lambda \with 1$, $\lrexpo \with \lrexpo$, $\cL \with \cL$, $\Theta_n \with \Theta_n$, $p_n \with p_n$, $\gamma_n \with \gamma_n$, $\bbA_n \with \bbA_n$, $\m_n \with \m_n$, $\mu_n \with \mu_n$ in the notation of \cref{cor:gen_convergence_rate_Loja_gen_l_rate_Lipschitz_alpha}); \lref{m}; \lref{Theta}; \lref{cL_KL}; the fact that $\cL$ is a \KL\ function}[verbs=e]{\lref{result}. }
\end{aproof}

\cfclear
\begin{savenotes}
\begin{samepage}
\begin{tcolorbox}[colback=white!95!gray,
                  colframe=black,
                  boxrule=0.5pt,
                  sharp corners,
                  enhanced,
                 ]
\begin{athm}{cor}{cor:dnns_2}[\resname{\DNNs\ with softplus activation}]
Let $\fd, L, M \in \N$, $\ell_0, \allowbreak \ell_1, \allowbreak \dots, \allowbreak \ell_L \allowbreak \in \N$, $\fx_1, \allowbreak \fx_2, \allowbreak \dots, \allowbreak \fx_M \in \allowbreak \R^{\ell_0}$, $\fy_1, \allowbreak \fy_2, \allowbreak \dots, \allowbreak \fy_M \in \allowbreak \R^{\ell_L}$ satisfy $\fd = \allowbreak \sum_{i = 1}^{L} \ell_i (\ell_{i - 1} + 1)$, for every $\theta = \allowbreak (\theta_1, \allowbreak \dots, \allowbreak \theta_{\fd}) \allowbreak \in \R^{\fd}$ let $\cN^{k, \theta} = \allowbreak (\cN^{k, \theta}_1, \allowbreak \dots, \allowbreak \cN^{k, \theta}_{\ell_k}) \colon \allowbreak \R^{\ell_0} \allowbreak \to \R^{\ell_k}$, $k \in \allowbreak \{0, \allowbreak 1, \allowbreak \dots, \allowbreak L\}$, satisfy for all $k \in \allowbreak \{0, \allowbreak 1, \allowbreak \dots, \allowbreak L - 1\}$, $x = (x_1, \allowbreak \dots, \allowbreak x_{\ell_0}) \in \allowbreak \R^{\ell_0}$, $i \in \allowbreak \{1, \allowbreak 2, \allowbreak \dots, \allowbreak \ell_{k + 1}\}$ that
\begin{multline}\llabel{cN}
\textstyle \cN^{k + 1, \theta}_i(x) = \theta_{\ell_{k + 1} \ell_{k} + i + \sum_{h = 1}^{k} \ell_h (\ell_{h - 1} + 1)} \\
\textstyle + \sum_{j = 1}^{\ell_{k}} \theta_{(i - 1) \ell_{k} + j + \sum_{h = 1}^{k} \ell_h (\ell_{h - 1} + 1)} \bigl[x_j \mathbbm{1}_{\{0\}}(k) + \ln(1 + \exp(\cN^{k, \theta}_j(x))) \mathbbm{1}_{\N}(k)\bigr],
\end{multline}
let $\cL \colon \R^{\fd} \to \R$ satisfy for all $\theta \in \R^{\fd}$ that
\begin{equation}\llabel{loss}
\textstyle \cL(\theta) = \frac{1}{M} \sum_{m = 1}^{M} \norm{\cN^{L, \theta}(\fx_m) - \fy_m}^2,
\end{equation}
let $\alpha \in [0, 1)$, $\lrexpo \in (\nicefrac{3}{4}, 1]$, let $\Theta \allowbreak \colon \allowbreak \N_0 \allowbreak \to \allowbreak \R^\fd$ and $p \colon \N \to \R$ be bounded, let $\gamma \colon \N \allowbreak \to \R$, $\bbA \colon \N \to \R^{\fd \times \fd}$, $\m \allowbreak \colon \allowbreak \N_0 \allowbreak \to \R^\fd$, and $\mu \colon \N \to \R^{\fd}$ satisfy for all $n \in \N$ that
\begin{gather}\llabel{m}
\textstyle \m_n = \alpha \m_{n - 1} + (1 - \alpha) \bigl[(\nabla \cL)(\Theta_{n - 1}) + \gamma_n \mu_n\bigr], \\
\llabel{Theta}
\textstyle \Theta_n = \Theta_{n - 1} - \gamma_n \bbA_n \bigl[p_n \m_n + (1 - p_n) (\nabla \cL)(\Theta_{n - 1})\bigr],
\end{gather}
and $\sup_{k \in \N} (k^{\lrexpo} \abs{\gamma_k} (\norm{\bbA_k} + \norm{\mu_k})) < \infty$, and assume for all $n \in \N$ that $\gamma_n n^{\lrexpo} \bbA_n - \bbI_{\fd}$ is symmetric positive semi-definite \cfout. Then there exist $\vartheta \allowbreak \in \allowbreak \R^{\fd}$, $\rho \allowbreak \in (0, \allowbreak \infty)$ which satisfy for all $n \in \N$ that
\begin{equation}\llabel{result}
\textstyle \norm{(\nabla \cL) (\vartheta)} + \abs{\cL(\Theta_n) - \cL(\vartheta)} + \norm{\Theta_n - \vartheta}^{\rho} \le \rho \bigl[\sum_{j = 1}^{n} j^{- \lrexpo}\bigr]^{-1} \dott
\end{equation}
\end{athm}
\end{tcolorbox}
\end{samepage}
\end{savenotes}
\begin{aproof}
\argument{}{ that $\R \ni x \mapsto \ln(1 + \exp(x)) \in \R$ is analytic. \llabel{softplus_analytic}}
\argument{\cref{cor:dnns_1}; \lref{softplus_analytic}}[verbs=e]{\lref{result}. }
\end{aproof}

\section{Strong convergence for specific GD optimization me\-thods}
\label{section:UGD_applications}

In \cref{cor:gen_convergence_rate_Loja_l_rate} and its consequences in \cref{cor:gen_convergence_rate_Loja_l_rate_TEMP} and \cref{cor:gen_convergence_rate_Loja_gen_l_rate_Lipschitz_alpha} in \cref{subsec:convergence_UGD_KL_Lipschitz} we provide a general unified convergence analysis for \GD\ optimization methods. To demonstrate the generality of these findings, we briefly apply in this section \cref{cor:gen_convergence_rate_Loja_gen_l_rate_Lipschitz_alpha} to each of the optimizers in \ref{intro:momentum}--\ref{intro:Yogi} from \cref{sec:introduction}.

\subsection{Classical momentum}
\label{subsec:momentum}

In the following result, \cref{cor:UGD_to_momentum}, we apply \cref{cor:gen_convergence_rate_Loja_gen_l_rate_Lipschitz_alpha} from \cref{subsec:convergence_UGD_KL_Lipschitz} above in the special situation where the optimization method is the momentum optimizer \cite{POLYAK19641} (cf., \eg, \cite[Subsection~6.3]{JentzenBookDeepLearning2023}).

\cfclear
\begin{savenotes}
\begin{samepage}
\begin{tcolorbox}[colback=white!95!gray,
                  colframe=black,
                  boxrule=0.5pt,
                  sharp corners,
                  enhanced,
                 ]
\begin{athm}{cor}{cor:UGD_to_momentum}[\resname{Momentum}]
Let $\fd \in \N$, $\alpha \in [0, 1)$, $\lrexpo \in (\nicefrac{3}{4}, 1]$, let $\cL \colon \R^{\fd} \to \R$ be a\cfadd{def:Loja} \KL\ function, assume that $\nabla \cL$ is locally Lipschitz continuous, let $\gamma \colon \N \allowbreak \to (0, \allowbreak \infty)$, $\Theta \allowbreak \colon \allowbreak \N_0 \allowbreak \to \allowbreak \R^\fd$, and $\m \allowbreak \colon \allowbreak \N_0 \allowbreak \to \allowbreak \R^\fd$ satisfy for all $n \in \N$ that
\begin{equation}\llabel{m_Theta}
\textstyle \m_n = \alpha \m_{n - 1} + (1 - \alpha) (\nabla \cL)(\Theta_{n - 1}) \qqandqq \Theta_n = \Theta_{n - 1} - \gamma_n \m_n,
\end{equation}
and assume $\sup_{n \in \N} (\norm{\Theta_n} + \sum_{k = -1}^{1} (\gamma_n n^{\lrexpo})^k) < \infty$ \cfout. Then there exist $\vartheta \in \allowbreak \R^{\fd}$, $\rho \in (0, \allowbreak \infty)$ which satisfy for all $n \in \N$ that
\begin{equation}\llabel{result}
\textstyle \norm{(\nabla \cL)(\vartheta)} + \abs{\cL(\Theta_n) - \cL(\vartheta)} + \norm{\Theta_n - \vartheta}^{\rho} \le \rho \bigl[\sum_{j = 1}^{n} j^{- \lrexpo}\bigr]^{-1} \dott
\end{equation}
\end{athm}
\end{tcolorbox}
\end{samepage}
\end{savenotes}
\begin{aproof}
\argument{the assumption that $\sup_{n \in \N} (\norm{\Theta_n} + \sum_{k = -1}^{1} (\gamma_n n^{\lrexpo})^k) < \infty$; \cref{cor:gen_convergence_rate_Loja_gen_l_rate_Lipschitz_alpha} (applied for every $n \in \N$ with $\fd \with \fd$, $\alpha \with \alpha$, $\lambda \with \inf_{m \in \N} (\gamma_m m^{\lrexpo})$, $\lrexpo \with \lrexpo$, $\cL \with \cL$, $\Theta_n \with \Theta_n$, $p_n \with 1$, $\gamma_n \with \gamma_n$, $\bbA_n \with \bbI_{\fd}$, $\m_n \with \m_n$, $\mu_n \with 0$ in the notation of \cref{cor:gen_convergence_rate_Loja_gen_l_rate_Lipschitz_alpha}); \lref{m_Theta}}[verbs=e]{\lref{result}. }
\end{aproof}

\subsection{Nesterov accelerated gradient (NAG)}
\label{subsec:NAG}

In the following result, \cref{cor:UGD_to_Nesterov_momentum}, we apply \cref{cor:gen_convergence_rate_Loja_gen_l_rate_Lipschitz_alpha} from \cref{subsec:convergence_UGD_KL_Lipschitz} above in the special situation where the optimization method is the \NAG\ optimizer \cite{Nesterov1983AMF} (cf., \eg, \cite[Subsection~6.4]{JentzenBookDeepLearning2023}).

\cfclear
\begin{savenotes}
\begin{samepage}
\begin{tcolorbox}[colback=white!95!gray,
                  colframe=black,
                  boxrule=0.5pt,
                  sharp corners,
                  enhanced,
                 ]
\begin{athm}{cor}{cor:UGD_to_Nesterov_momentum}[\NAG]
Let $\fd \in \N$, $\alpha \in [0, 1)$, $\lrexpo \in (\nicefrac{3}{4}, 1]$, let $\cL \colon \R^{\fd} \to \R$ be a\cfadd{def:Loja} \KL\ function, assume that $\nabla \cL$ is locally Lipschitz continuous, let $\gamma \colon \N \allowbreak \to (0, \allowbreak \infty)$, $\Theta \allowbreak \colon \allowbreak \N_0 \allowbreak \to \allowbreak \R^\fd$, and $\m \allowbreak \colon \allowbreak \N_0 \allowbreak \to \allowbreak \R^\fd$ satisfy for all $n \in \N$ that
\begin{equation}\llabel{m_Theta}
\textstyle \m_n = \alpha \m_{n - 1} + (1 - \alpha) (\nabla \cL)(\Theta_{n - 1} - \alpha \gamma_n \m_{n - 1}) \qandq \Theta_n = \Theta_{n - 1} - \gamma_n \m_n,
\end{equation}
and assume $\sup_{n \in \N} (\norm{\Theta_n} + \sum_{k = -1}^{1} (\gamma_n n^{\lrexpo})^k) < \infty$ \cfout. Then there exist $\vartheta \in \allowbreak \R^{\fd}$, $\rho \in (0, \allowbreak \infty)$ which satisfy for all $n \in \N$ that
\begin{equation}\llabel{result}
\textstyle \norm{(\nabla \cL)(\vartheta)} + \abs{\cL(\Theta_n) - \cL(\vartheta)} + \norm{\Theta_n - \vartheta}^{\rho} \le \rho \bigl[\sum_{j = 1}^{n} j^{- \lrexpo}\bigr]^{-1} \dott
\end{equation}
\end{athm}
\end{tcolorbox}
\end{samepage}
\end{savenotes}
\begin{aproof}
Throughout this proof for every $n \in \N$ let $\mu_n \in \R^{\fd}$ satisfy
\begin{equation}\llabel{mu}
\textstyle \mu_n = (\gamma_n)^{-1} \bigl[(\nabla \cL)(\Theta_{n - 1} - \alpha \gamma_n \m_{n - 1}) - (\nabla \cL)(\Theta_{n - 1})\bigr] \dott
\end{equation}
\argument{\lref{m_Theta}; \lref{mu}}{for all $n \in \N$ that
\begin{equation}\llabel{m_mu_Theta}
\textstyle \m_n = \alpha \m_{n - 1} + (1 - \alpha) \bigl[(\nabla \cL)(\Theta_{n - 1}) + \gamma_n \mu_n\bigr] \qandq \Theta_n = \Theta_{n - 1} - \gamma_n \m_n \dott
\end{equation}
}
\argument{\lref{m_Theta}; the assumption that $\sup_{n \in \N} (\norm{\Theta_n} \allowbreak + \sum_{k = -1}^{1} (\gamma_n n^{\lrexpo})^k) \allowbreak < \infty$; the fact that $\sup_{n \in \N} [1 + \frac{1}{n}]^{- \lrexpo} < \infty$; the triangle inequality}{that
\begin{equation}\llabel{n_m_bounded}
\textstyle \sup_{n \in \N} [\gamma_{n + 1} \norm{\m_n}] = \sup_{n \in \N} \bigl([\gamma_{n + 1} (n + 1)^{\lrexpo}] [\gamma_n n^{\lrexpo}]^{- 1} [1 + \frac{1}{n}]^{- \lrexpo} \Norm{\Theta_n - \Theta_{n - 1}}\bigr) < \infty \dott
\end{equation}
}
\argument{\lref{n_m_bounded}; the triangle inequality}{that $\sup_{n \in \N} \norm{\Theta_{n - 1} - \alpha \gamma_n \m_{n - 1}} < \infty$. }
Let $B, L \in (0, \infty)$ satisfy for all $n \in \N$, $u, w \in \{x \in \R^{\fd} \colon \norm{x} \le B\}$ that
\begin{equation}\llabel{B_L}
\textstyle \norm{\Theta_{n - 1} - \alpha \gamma_n \m_{n - 1}} + \norm{\Theta_{n - 1}} \le B \qqandqq \norm{(\nabla \cL)(u) - (\nabla \cL)(w)} \le L \norm{u - w} \dott
\end{equation}
\startnewargseq
\argument{\lref{B_L}; \cref{lemma:momentum_priori_bound_2} (applied for every $n \in \N$ with $\alpha_n \with \alpha$, $A_{n - 1} \with \norm{(\nabla \cL)(\Theta_{n - 1} - \alpha \gamma_n \m_{n - 1})}$, $x_{n - 1} \with \norm{\m_{n - 1}}$ in the notation of \cref{lemma:momentum_priori_bound_2}); the triangle inequality}{that $\sup_{n \in \N} \norm{\m_{n - 1}} < \infty$. \llabel{m_bound}}
\argument{\lref{mu}; \lref{B_L}; \lref{m_bound}}{for all $n \in \N$ that
\begin{equation}\llabel{mu_bound}
\textstyle \norm{\mu_n} \le (\gamma_n)^{- 1} L \norm{\alpha \gamma_n \m_{n - 1}} = L \alpha \norm{\m_{n - 1}} \le L \alpha \sup_{k \in \N}\norm{\m_{k - 1}} < \infty \dott
\end{equation}
}
\argument{the assumption that $\sup_{n \in \N} (\norm{\Theta_n} + \sum_{k = -1}^{1} (\gamma_n n^{\lrexpo})^k) < \infty$; \cref{cor:gen_convergence_rate_Loja_gen_l_rate_Lipschitz_alpha} (applied for every $n \in \N$ with $\fd \with \fd$, $\alpha \with \alpha$, $\lambda \with \inf_{m \in \N} (\gamma_m m^{\lrexpo})$, $\lrexpo \with \lrexpo$, $\cL \with \cL$, $\Theta_n \with \Theta_n$, $p_n \with 1$, $\gamma_n \with \gamma_n$, $\bbA_n \with \bbI_{\fd}$, $\m_n \with \m_n$, $\mu_n \with \mu_n$ in the notation of \cref{cor:gen_convergence_rate_Loja_gen_l_rate_Lipschitz_alpha}); \lref{m_mu_Theta}; \lref{mu_bound}}[verbs=e]{\lref{result}. }
\end{aproof}

\subsection{Root mean square propagation (RMSprop)}
\label{subsec:RMSprop}

In the following result, \cref{cor:UGD_to_RMSProp}, we apply \cref{cor:gen_convergence_rate_Loja_gen_l_rate_Lipschitz_alpha} from \cref{subsec:convergence_UGD_KL_Lipschitz} above in the special situation where the optimization method is the \RMSprop\ optimizer \cite{HintonRMSProp} (cf., \eg, \cite[Subsection~6.6]{JentzenBookDeepLearning2023}).

\cfclear
\begin{savenotes}
\begin{samepage}
\begin{tcolorbox}[colback=white!95!gray,
                  colframe=black,
                  boxrule=0.5pt,
                  sharp corners,
                  enhanced,
                 ]
\begin{athm}{cor}{cor:UGD_to_RMSProp}[\RMSprop]
Let $\fd \in \N$, $\beta \in [0, 1]$, $\eps \in (0, \infty)$, $\lrexpo \in (\nicefrac{3}{4}, 1]$, let $\cL \colon \R^{\fd} \to \R$ be a\cfadd{def:Loja} \KL\ function, assume that $\nabla \cL$ is locally Lipschitz continuous, let $\gamma \colon \N \allowbreak \to (0, \allowbreak \infty)$, $\Theta = \allowbreak (\Theta^1, \allowbreak \dots, \allowbreak \Theta^{\fd}) \allowbreak \colon \allowbreak \N_0 \allowbreak \to \allowbreak \R^\fd$, and $\bbM = (\bbM^1, \allowbreak \dots, \allowbreak \bbM^{\fd}) \colon \allowbreak \N_0 \allowbreak \to \R^\fd$ satisfy for all $n \in \N$, $i \in \allowbreak \{1, \allowbreak 2, \allowbreak \dots, \allowbreak \fd\}$ that
\begin{gather}\llabel{M}
\textstyle \bbM_n^i = \beta \bbM_{n - 1}^i + (1 - \beta) \abs{(\nabla \cL)^i(\Theta_{n - 1})}^2 \\ \llabel{Theta}
\textstyle \andqq \Theta_n^i = \Theta_{n - 1}^i - \gamma_n \bigl[\eps + \abs{\bbM_n^i}\bigr]^{- \nicefrac{1}{2}} (\nabla \cL)^i(\Theta_{n - 1}),
\end{gather}
and assume $\sup_{n \in \N} (\norm{\Theta_n} + \sum_{k = -1}^{1} (\gamma_n n^{\lrexpo})^k) < \infty$ \cfout. Then there exist $\vartheta \in \allowbreak \R^{\fd}$, $\rho \in (0, \allowbreak \infty)$ which satisfy for all $n \in \N$ that
\begin{equation}\llabel{result}
\textstyle \norm{(\nabla \cL)(\vartheta)} + \abs{\cL(\Theta_n) - \cL(\vartheta)} + \norm{\Theta_n - \vartheta}^{\rho} \le \rho \bigl[\sum_{j = 1}^{n} j^{- \lrexpo}\bigr]^{-1} \dott
\end{equation}
\end{athm}
\end{tcolorbox}
\end{samepage}
\end{savenotes}
\begin{aproof}
Throughout this proof for every $n \in \N$ let $\bbA_n \in \R^{\fd \times \fd}$ satisfy
\begin{equation}\llabel{bbA}
\textstyle \bbA_n = \bigl(\bigl[\eps + \abs{\bbM_n^i}\bigr]^{- \nicefrac{1}{2}} \mathbbm{1}_{\{j\}}(i)\bigr)_{(i, j) \in \{1, 2, \dots, \fd\}^2} \dott
\end{equation}
\startnewargseq
\argument{\lref{Theta}; \lref{bbA}}{for all $n \in \N$ that
\begin{equation}\llabel{Theta_new}
\textstyle \Theta_{n} = \Theta_{n - 1} - \gamma_n \bbA_n (\nabla \cL)(\Theta_{n - 1}) \dott
\end{equation}
}
\argument{the assumption that $\sup_{n \in \N} \norm{\Theta_n} < \infty$; the assumption that $\nabla \cL$ is locally Lipschitz continuous}{that $\sup_{n \in \N} \norm{(\nabla \cL)(\Theta_n)} < \infty$. \llabel{cG_bound}}
\argument{\lref{M}; \lref{cG_bound}; \cref{lemma:momentum_priori_bound_2} (applied for every $n \in \N$, $i \in \{1, 2, \dots, \fd\}$ with $A_{n - 1} \with \abs{(\nabla \cL)^i(\Theta_{n - 1})}^2$, $\alpha_n \with \beta$, $x_{n - 1} \with \abs{\bbM_{n - 1}^i}$ in the notation of \cref{lemma:momentum_priori_bound_2}); the triangle inequality}{that $\sup_{n \in \N} \norm{\bbM_n} < \infty$. \llabel{bbM_bound}}
\argument{the assumption that $\sup_{n \in \N} (\sum_{k = -1}^{1} (\gamma_n n^{\lrexpo})^k) < \infty$; \lref{bbA}; \lref{bbM_bound}}{that there exists $\lambda \in (0, \infty)$ which satisfies for all $n \in \N$ that $\gamma_n n^{\lrexpo} \bbA_n - \lambda \bbI_{\fd}$ is symmetric positive semi-definite and $\sup_{k \in \N} \norm{\gamma_k k^{\lrexpo} \bbA_k} < \infty$. \llabel{bbA_bound}}
\argument{\lref{Theta_new}; \lref{bbA_bound}; the assumption that $\sup_{n \in \N} \norm{\Theta_n} < \infty$; \cref{cor:gen_convergence_rate_Loja_gen_l_rate_Lipschitz_alpha} (applied for every $n \in \N$ with $\fd \with \fd$, $\alpha \with 0$, $\lambda \with \lambda$, $\lrexpo \with \lrexpo$, $\cL \with \cL$, $\Theta_n \with \Theta_n$, $p_n \with 1$, $\gamma_n \with \gamma_n$, $\bbA_n \with \bbA_n$, $\m_n \with (\nabla \cL)(\Theta_{n - 1})$, $\mu_n \with 0$ in the notation of \cref{cor:gen_convergence_rate_Loja_gen_l_rate_Lipschitz_alpha})}[verbs=e]{\lref{result}. }
\end{aproof}

\subsection{Adaptive moment estimation (Adam)}
\label{subsec:Adam}

In the following result, \cref{cor:UGD_to_Adam}, we apply \cref{cor:gen_convergence_rate_Loja_gen_l_rate_Lipschitz_alpha} from \cref{subsec:convergence_UGD_KL_Lipschitz} above in the special situation where the optimization method is the \Adam\ optimizer \cite{KingmaBaAdam} (cf., \eg, \cite[Subsection~6.8]{JentzenBookDeepLearning2023}).

\cfclear
\begin{savenotes}
\begin{samepage}
\begin{tcolorbox}[colback=white!95!gray,
                  colframe=black,
                  boxrule=0.5pt,
                  sharp corners,
                  enhanced,
                 ]
\begin{athm}{cor}{cor:UGD_to_Adam}[\Adam]
Let $\fd \in \N$, $\beta_1, \beta_2 \in [0, 1)$, $\eps \in (0, \infty)$, $\lrexpo \in (\nicefrac{3}{4}, 1]$, let $\cL \colon \allowbreak \R^{\fd} \allowbreak \to \R$ be a\cfadd{def:Loja} \KL\ function, assume that $\nabla \cL$ is locally Lipschitz continuous, let $\gamma \colon \N \to (0, \infty)$ and $\Theta\Index{k}{} \allowbreak = \allowbreak (\Theta\Index{k}[1]{}, \allowbreak \dots, \allowbreak \Theta\Index{k}[\fd]{}) \allowbreak \colon \allowbreak \N_0 \allowbreak \to \allowbreak \R^\fd$, $k \in \N_0$, satisfy for all $n \in \N$, $k \in \allowbreak \{1, \allowbreak 2\}$, $i \in \allowbreak \{1, \allowbreak 2, \allowbreak \dots, \allowbreak \fd\}$ that
\begin{gather}\llabel{m_M}
\textstyle \Theta\Index{k}[i]{n} = \beta_k \Theta\Index{k}[i]{n - 1} + (1 - \beta_k) \bigl[(\nabla \cL)^i(\Theta\Index{0}{n - 1})\bigr]^k \\ \llabel{Theta}
\textstyle \andqq \Theta\Index{0}[i]{n} = \Theta\Index{0}[i]{n - 1} - \gamma_n \bigl[\frac{\Theta\Index{1}[i]{n}}{1 - (\beta_1)^n}\bigr] \Bigl[\eps + \bigl[\frac{\abs{\Theta\Index{2}[i]{n}}}{1 - (\beta_2)^n}\bigr]^{\nicefrac{1}{2}}\Bigr]^{-1},
\end{gather}
and assume $\sup_{n \in \N} (\norm{\Theta\Index{0}{n}} + \sum_{k = -1}^{1} (\gamma_n n^{\lrexpo})^k) < \infty$ \cfout. Then there exist $\vartheta  \allowbreak \in \allowbreak \R^{\fd}$, $\rho \in (0, \allowbreak \infty)$ which satisfy for all $n \in \N$ that
\begin{equation}\llabel{result}
\textstyle  \norm{(\nabla \cL)(\vartheta)} + \abs{\cL(\Theta\Index{0}{n}) - \cL(\vartheta)} + \norm{\Theta\Index{0}{n} - \vartheta}^{\rho} \le \rho \bigl[\sum_{j = 1}^{n} j^{- \lrexpo}\bigr]^{-1} \dott
\end{equation}
\end{athm}
\end{tcolorbox}
\end{samepage}
\end{savenotes}
\begin{aproof}
Throughout this proof for every $n \in \N$ let $\bbA_n \in \R^{\fd \times \fd}$ satisfy
\begin{equation}\llabel{bbA}
\textstyle \bbA_n = \Bigl(\frac{\mathbbm{1}_{\{j\}}(i)}{1 - (\beta_1)^n} \Bigl[\eps + \bigl[\frac{\abs{\Theta\Index{2}[i]{n}}}{1 - (\beta_2)^n}\bigr]^{\nicefrac{1}{2}}\Bigr]^{-1}\Bigr)_{(i, j) \in \{1, 2, \dots, \fd\}^2} \dott
\end{equation}
\startnewargseq
\argument{\lref{Theta}; \lref{bbA}}{for all $n \in \N$ that
\begin{equation}\llabel{Theta_new}
\textstyle \Theta\Index{0}{n} = \Theta\Index{0}{n - 1} - \gamma_n \bbA_n \Theta\Index{1}{n} \dott
\end{equation}
}
\argument{the assumption that $\sup_{n \in \N} \norm{\Theta\Index{0}{n}} < \infty$; the assumption that $\nabla \cL$ is locally Lipschitz continuous}{that $\sup_{n \in \N} \norm{(\nabla \cL)(\Theta\Index{0}{n})} < \infty$. \llabel{cG_bound}}
\argument{\lref{m_M}; \lref{cG_bound}; \cref{lemma:momentum_priori_bound_2} (applied for every $n \in \N$, $i \in \{1, 2, \dots, \fd\}$ with $A_{n - 1} \with \abs{(\nabla \cL)^i(\Theta\Index{0}{n - 1})}^2$, $\alpha_n \with \beta_2$, $x_{n - 1} \with \abs{\Theta\Index{2}[i]{n - 1}}$ in the notation of \cref{lemma:momentum_priori_bound_2}); the triangle inequality}{that $\sup_{n \in \N} \norm{\Theta\Index{2}{n}} \allowbreak < \infty$. \llabel{bbM_bound}}
\argument{the assumption that $\sup_{n \in \N} (\sum_{k = -1}^{1} (\gamma_n n^{\lrexpo})^k) < \infty$; \lref{bbA}; \lref{bbM_bound}; the fact that for all $n \in \N$ it holds that
\begin{equation}\llabel{alpha_beta_bound}
\textstyle 1 \le \frac{1}{1 - (\beta_1)^n} \le \frac{1}{1 - \beta_1} < \infty \qqandqq 1 \le \frac{1}{1 - (\beta_2)^n} \le \frac{1}{1 - \beta_2} < \infty
\end{equation}
}{that there exists $\lambda \in (0, \infty)$ which satisfies for all $n \in \N$ that $\gamma_n n^{\lrexpo} \bbA_n - \lambda \bbI_{\fd}$ is symmetric positive semi-definite and $\sup_{k \in \N} \norm{\gamma_k k^{\lrexpo} \bbA_k} < \infty$. \llabel{bbA_bound}}
\argument{\lref{m_M}; \lref{Theta_new}; \lref{bbA_bound}; the assumption that $\sup_{n \in \N} \norm{\Theta\Index{0}{n}} < \infty$; \cref{cor:gen_convergence_rate_Loja_gen_l_rate_Lipschitz_alpha} (applied for every $n \in \N$ with $\fd \with \fd$, $\alpha \with \beta_1$, $\lambda \with \lambda$, $\lrexpo \with \lrexpo$, $\cL \with \cL$, $\Theta_n \with \Theta\Index{0}{n}$, $p_n \with 1$, $\gamma_n \with \gamma_n$, $\bbA_n \with \bbA_n$, $\m_n \with \Theta\Index{1}{n}$, $\mu_n \with 0$ in the notation of \cref{cor:gen_convergence_rate_Loja_gen_l_rate_Lipschitz_alpha})}[verbs=e]{\lref{result}. }
\end{aproof}

We note that in the case $\lrexpo < 1$ we have that \cref{cor:UGD_to_Adam} exploits that the learning rates $\gamma_n \in (0,\infty)$, $n \in \N$, converge to zero with an appropriate speed of convergence in the sense that $\sum_{k = - 1}^1 (\sup_{n \in \N} \abs{\gamma_n n^{\lrexpo}}^k) < \infty$ to conclude that there exist $\rho$, $\fC \in (0, \infty)$ such that for all $n \in \N$ it holds that $\norm{\Theta\Index{0}{n} - \vartheta} \le \fC n^{- \rho}$ (the error of \Adam\ decays with polynomial speed of convergence $n^{- \rho}$). We point out that in the situation where the learning rates $\gamma_n \in (0, \infty)$, $n \in \N$, do not decay to zero but stay constant a much quicker speed of convergence (exponential speed of convergence $\rho^n$ for some $\rho \in (0, 1)$) can be obtained under suitable additional assumptions \cite[Item~(iv) in Theorem~1.2]{DereichJentzenRiekert_sharpconvergencerates_Adam}.

We also observe that \cref{cor:UGD_to_Adam.result} in \cref{cor:UGD_to_Adam} establishes convergence of the optimization process $(\Theta\Index{0}{n})_{n \in \N_0}$ to a zero $\vartheta \in \R^{\fd}$ of the gradient $\nabla \cL \colon \R^{\fd} \to \R^{\fd}$ of the objective function (a criticial point of the objective function $\cL \colon \R^{\fd} \to \R$) for \Adam\ applied to a \emph{deterministic} optimization problem. However, when \Adam\ is applied to certain stochastic optimization problems, then it is known that \Adam\ does typically \emph{not} converge to a zero of the gradient of the objective function (see \cite[Theorem~1.1]{DereichDoJentzenPhilippeAdamSymmetry} for details) but instead converges to a zero of the \Adam\ vector field \cite[Definition~2.4]{DereichJentzenAdamRates}; see \cite[Corollary~1.10]{DereichJentzenKassing_ODE_Adam}.

In the next result, \cref{cor:Adam_for_dnns_2}, we specialize \cref{cor:UGD_to_Adam} to the situation where the underlying optimization problem is nothing else but empirical risk minimization for fully-connected feedforward \DNNs\ with analytic activation functions. In \cref{cor:Adam_for_dnns_2.cN} in \cref{cor:Adam_for_dnns_2} the realization functions for the considered fully-connected feedforward \DNNs\ are formulated.

\cfclear
\begin{savenotes}
\begin{samepage}
\begin{tcolorbox}[colback=white!95!gray,
                  colframe=black,
                  boxrule=0.5pt,
                  sharp corners,
                  enhanced,
                 ]
\begin{athm}{cor}{cor:Adam_for_dnns_2}[\resname{\Adam\ training \DNNs}]
Let $\fd, L, M \in \N$, $\ell_0, \allowbreak \ell_1, \allowbreak \dots, \allowbreak \ell_L \allowbreak \in \N$, $\fx_1, \allowbreak \fx_2, \allowbreak \dots, \allowbreak \fx_M \in \allowbreak \R^{\ell_0}$, $\fy_1, \allowbreak \fy_2, \allowbreak \dots, \allowbreak \fy_M \in \allowbreak \R^{\ell_L}$ satisfy $\fd = \allowbreak \sum_{i = 1}^{L} \ell_i (\ell_{i - 1} + 1)$, let $\act \colon \R \to \R$ be analytic, for every $\theta = \allowbreak (\theta_1, \allowbreak \dots, \allowbreak \theta_{\fd}) \allowbreak \in \R^{\fd}$ let $\cN^{k, \theta} = \allowbreak (\cN^{k, \theta}_1, \allowbreak \dots, \allowbreak \cN^{k, \theta}_{\ell_k}) \colon \allowbreak \R^{\ell_0} \allowbreak \to \R^{\ell_k}$, $k \in \allowbreak \{0, \allowbreak 1, \allowbreak \dots, \allowbreak L\}$, satisfy for all $k \in \allowbreak \{0, \allowbreak 1, \allowbreak \dots, \allowbreak L - 1\}$, $x = (x_1, \allowbreak \dots, \allowbreak x_{\ell_0}) \in \allowbreak \R^{\ell_0}$, $i \in \allowbreak \{1, \allowbreak 2, \allowbreak \dots, \allowbreak \ell_{k + 1}\}$ that
\begin{multline}\llabel{cN}
\textstyle \cN^{k + 1, \theta}_i(x) = \theta_{\ell_{k + 1} \ell_{k} + i + \sum_{h = 1}^{k} \ell_h (\ell_{h - 1} + 1)} \\
\textstyle + \sum_{j = 1}^{\ell_{k}} \theta_{(i - 1) \ell_{k} + j + \sum_{h = 1}^{k} \ell_h (\ell_{h - 1} + 1)} \bigl[x_j \mathbbm{1}_{\{0\}}(k) + \act(\cN^{k, \theta}_j(x)) \mathbbm{1}_{\N}(k)\bigr],
\end{multline}
let $\cL \colon \R^{\fd} \to \R$ satisfy for all $\theta \in \R^{\fd}$ that
\begin{equation}\llabel{loss}
\textstyle \cL(\theta) = \frac{1}{M} \sum_{m = 1}^{M} \norm{\cN^{L, \theta}(\fx_m) - \fy_m}^2,
\end{equation}
let $\beta_1, \beta_2 \in [0, 1)$, $\eps \in (0, \infty)$, $\lrexpo \in (\nicefrac{3}{4}, 1]$, let $\gamma \colon \N \to (0, \infty)$ and $\Theta\Index{k}{} \allowbreak = (\Theta\Index{k}[1]{}, \allowbreak \dots, \allowbreak \Theta\Index{k}[\fd]{}) \allowbreak \colon \allowbreak \N_0 \allowbreak \to \allowbreak \R^\fd$, $k \in \N_0$, satisfy for all $k \in \allowbreak \{1, \allowbreak 2\}$, $n \in \N$, $i \in \allowbreak \{1, \allowbreak 2, \allowbreak \dots, \allowbreak \fd\}$ that
\begin{gather}\llabel{m_M}
\textstyle \Theta\Index{k}[i]{n} = \beta_k \Theta\Index{k}[i]{n - 1} + (1 - \beta_k) \bigl[(\nabla \cL)^i(\Theta\Index{0}{n - 1})\bigr]^k \\ \llabel{Theta}
\textstyle \andqq \Theta\Index{0}[i]{n} = \Theta\Index{0}[i]{n - 1} - \gamma_n \bigl[\frac{\Theta\Index{1}[i]{n}}{1 - (\beta_1)^n}\bigr] \Bigl[\eps + \bigl[\frac{\abs{\Theta\Index{2}[i]{n}}}{1 - (\beta_2)^n}\bigr]^{\nicefrac{1}{2}}\Bigr]^{-1},
\end{gather}
and assume $\sup_{n \in \N} (\norm{\Theta\Index{0}{n}} + \sum_{k = -1}^{1} (\gamma_n n^{\lrexpo})^k) < \infty$ \cfout. Then there exist $\vartheta  \allowbreak \in \allowbreak \R^{\fd}$, $\rho \in (0, \allowbreak \infty)$ which satisfy for all $n \in \N$ that
\begin{equation}\llabel{result}
\textstyle  \norm{(\nabla \cL)(\vartheta)} + \abs{\cL(\Theta\Index{0}{n}) - \cL(\vartheta)} + \norm{\Theta\Index{0}{n} - \vartheta}^{\rho} \le \rho \bigl[\sum_{j = 1}^{n} j^{- \lrexpo}\bigr]^{-1} \dott
\end{equation}
\end{athm}
\end{tcolorbox}
\end{samepage}
\end{savenotes}
\begin{aproof}
\argument{the fact that $\act$ is analytic; \lref{cN}; \lref{loss}}{that $\cL$ is analytic (cf., \eg, \cite[Corollary~9.14.5 in Subsection~9.14]{JentzenBookDeepLearning2023}). \llabel{cL_analytic}}
\argument{\lref{cL_analytic}}{that $\cL$ is a \KL\ function (cf., \eg, \cite{MR2274510} and the references therein). }
\argument{the fact that $\cL$ is analytic}{that $\nabla \cL$ is locally Lipschitz continuous. \llabel{cL_KL}}
\argument{the fact that $\cL$ is a \KL\ function; the assumption that $\sup_{n \in \N} (\norm{\Theta\Index{0}{n}} + \sum_{k = -1}^{1} (\gamma_n n^{\lrexpo})^k) < \infty$; \cref{cor:UGD_to_Adam} (applied for every $k \in \N_0$ with $\fd \with \fd$, $\beta_1 \with \beta_1$, $\beta_2 \with \beta_2$, $\eps \with \eps$, $\lrexpo \with \lrexpo$, $\cL \with \cL$, $\gamma \with \gamma$, $\Theta\Index{k}{} \with \Theta\Index{k}{}$ in the notation of \cref{cor:UGD_to_Adam}); \lref{m_M}; \lref{Theta}; \lref{cL_KL}}[verbs=e]{\lref{result}. }
\end{aproof}

\subsection{Adaptive moment estimation maximum (Adamax)}
\label{subsec:Adamax}

In the following result, \cref{cor:UGD_to_Adamax}, we apply \cref{cor:gen_convergence_rate_Loja_gen_l_rate_Lipschitz_alpha} from \cref{subsec:convergence_UGD_KL_Lipschitz} above in the special situation where the optimization method is the \Adamax\ optimizer \cite{KingmaBaAdam} (cf., \eg, \cite[Subsection~6.8.1]{JentzenBookDeepLearning2023}).

\cfclear
\begin{savenotes}
\begin{samepage}
\begin{tcolorbox}[colback=white!95!gray,
                  colframe=black,
                  boxrule=0.5pt,
                  sharp corners,
                  enhanced,
                 ]
\begin{athm}{cor}{cor:UGD_to_Adamax}[\Adamax]
Let $\fd \in \N$, $\alpha \in [0, 1)$, $\beta \in [0, 1]$, $\eps \in (0, \infty)$, $\lrexpo \in (\nicefrac{3}{4}, 1]$, let $\cL \colon \R^{\fd} \allowbreak \to \allowbreak \R$ be a\cfadd{def:Loja} \KL\ function, assume that $\nabla \cL$ is locally Lipschitz continuous, let $\gamma \colon \N \allowbreak \to (0, \allowbreak \infty)$, $\Theta = (\Theta^1, \dots, \Theta^{\fd}) \allowbreak \colon \allowbreak \N_0 \allowbreak \to \allowbreak \R^\fd$, $\m = (\m^1, \allowbreak \dots, \allowbreak \m^{\fd}) \colon \allowbreak \N_0 \allowbreak \to \R^\fd$, and $\bbM = (\bbM^1, \allowbreak \dots, \allowbreak \bbM^{\fd}) \colon \allowbreak \N_0 \allowbreak \to \R^\fd$ satisfy for all $n \in \N$, $i \in \{1, 2, \dots, \fd\}$ that
\begin{gather}\llabel{m_M}
\textstyle \m_n = \alpha \m_{n - 1} + (1 - \alpha) (\nabla \cL)(\Theta_{n - 1}), \qquad \bbM_n^i = \max\{\beta \bbM_{n - 1}^i,  \abs{(\nabla \cL)^i(\Theta_{n - 1})}\}, \\ \llabel{Theta}
\textstyle \andqq \Theta_n^i = \Theta_{n - 1}^i - \gamma_n \bigl[\frac{\m_n^i}{1 - \alpha^n}\bigr] \bigl[\eps + \abs{\bbM_n^i}\bigr]^{-1},
\end{gather}
and assume $\sup_{n \in \N} (\norm{\Theta_n} + \sum_{k = -1}^{1} (\gamma_n n^{\lrexpo})^k) < \infty$ \cfout. Then there exist $\vartheta \in \allowbreak \R^{\fd}$, $\rho \in (0, \allowbreak \infty)$ which satisfy for all $n \in \N$ that
\begin{equation}\llabel{result}
\textstyle \norm{(\nabla \cL)(\vartheta)} + \abs{\cL(\Theta_n) - \cL(\vartheta)} + \norm{\Theta_n - \vartheta}^{\rho} \le \rho \bigl[\sum_{j = 1}^{n} j^{- \lrexpo}\bigr]^{-1} \dott
\end{equation}
\end{athm}
\end{tcolorbox}
\end{samepage}
\end{savenotes}
\begin{aproof}
Throughout this proof for every $n \in \N$ let $\bbA_n \in \R^{\fd \times \fd}$ satisfy
\begin{equation}\llabel{bbA}
\textstyle \bbA_n = \Bigl(\frac{\mathbbm{1}_{\{j\}}(i)}{1 - \alpha^n} \bigl[\eps + \abs{\bbM_n^i}\bigr]^{-1}\Bigr)_{(i, j) \in \{1, 2, \dots, \fd\}^2} \dott
\end{equation}
\startnewargseq
\argument{\lref{Theta}; \lref{bbA}}{for all $n \in \N$ that
\begin{equation}\llabel{Theta_new}
\textstyle \Theta_{n} = \Theta_{n - 1} - \gamma_n \bbA_n \m_n \dott
\end{equation}
}
\argument{the assumption that $\sup_{n \in \N} \norm{\Theta_n} < \infty$; the assumption that $\nabla \cL$ is locally Lipschitz continuous}{that $\sup_{n \in \N} \norm{(\nabla \cL)(\Theta_n)} < \infty$. \llabel{cG_bound}}
\argument{\lref{m_M}; \lref{cG_bound}; the induction; the triangle inequality}{that $\sup_{n \in \N} \norm{\bbM_n} < \infty$. \llabel{bbM_bound}}
\argument{the assumption that $\sup_{n \in \N} (\sum_{k = -1}^{1} (\gamma_n n^{\lrexpo})^k) < \infty$; \lref{bbA}; \lref{bbM_bound}; the fact that for all $n \in \N$ it holds that
\begin{equation}\llabel{alpha_beta_bound}
\textstyle 0 < 1 - \alpha \le 1 - \alpha^n \le 1
\end{equation}
}{that there exists $\lambda \in (0, \infty)$ which satisfies for all $n \in \N$ that $\gamma_n n^{\lrexpo} \bbA_n - \lambda \bbI_{\fd}$ is symmetric positive semi-definite and $\sup_{k \in \N} \norm{\gamma_k k^{\lrexpo} \bbA_k} < \infty$. \llabel{bbA_bound}}
\argument{\lref{m_M}; \lref{Theta_new}; \lref{bbA_bound}; the assumption that $\sup_{n \in \N} \norm{\Theta_n} < \infty$; \cref{cor:gen_convergence_rate_Loja_gen_l_rate_Lipschitz_alpha} (applied for every $n \in \N$ with $\fd \with \fd$, $\alpha \with \alpha$, $\lambda \with \lambda$, $\lrexpo \with \lrexpo$, $\cL \with \cL$, $\Theta_n \with \Theta_n$, $p_n \with 1$, $\gamma_n \with \gamma_n$, $\bbA_n \with \bbA_n$, $\m_n \with \m_n$, $\mu_n \with 0$ in the notation of \cref{cor:gen_convergence_rate_Loja_gen_l_rate_Lipschitz_alpha})}[verbs=e]{\lref{result}. }
\end{aproof}

\subsection{Nesterov accelerated adaptive moment estimation (Nadam)}
\label{subsec:Nadam}

In the following result, \cref{cor:UGD_to_Nadam}, we apply \cref{cor:gen_convergence_rate_Loja_gen_l_rate_Lipschitz_alpha} from \cref{subsec:convergence_UGD_KL_Lipschitz} above in the special situation where the optimization method is the \Nadam\ optimizer \cite{DozatNadam2} (cf., \eg, \cite[Subsection~6.9]{JentzenBookDeepLearning2023}).

\cfclear
\begin{savenotes}
\begin{samepage}
\begin{tcolorbox}[colback=white!95!gray,
                  colframe=black,
                  boxrule=0.5pt,
                  sharp corners,
                  enhanced,
                 ]
\begin{athm}{cor}{cor:UGD_to_Nadam}[\Nadam]
Let $\fd \in \N$, $\alpha, \beta \in [0, 1)$, $\eps \in (0, \infty)$, $\lrexpo \in (\nicefrac{3}{4}, 1]$, let $\cL \colon \R^{\fd} \to \R$ be a\cfadd{def:Loja} \KL\ function, assume that $\nabla \cL$ is locally Lipschitz continuous, let $\gamma \colon \N \allowbreak \to (0, \allowbreak \infty)$, $\Theta = (\Theta^1, \dots, \Theta^{\fd}) \allowbreak \colon \allowbreak \N_0 \allowbreak \to \allowbreak \R^\fd$, $\m = (\m^1, \allowbreak \dots, \allowbreak \m^{\fd}) \colon \allowbreak \N_0 \allowbreak \to \R^\fd$, $\bbM = (\bbM^1, \allowbreak \dots, \allowbreak \bbM^{\fd}) \colon \allowbreak \N_0 \allowbreak \to \R^\fd$, and $\cM = (\cM^1, \allowbreak \dots, \allowbreak \cM^{\fd}) \colon \N \to \R^\fd$ satisfy for all $n \in \N$, $i \in \{1, 2, \dots, \fd\}$ that
\begin{gather}\llabel{m_M}
\textstyle \m_n = \alpha \m_{n - 1} + (1 - \alpha) (\nabla \cL)(\Theta_{n - 1}), \qquad \bbM_n^i = \beta \bbM_{n - 1}^i + (1 - \beta) \abs{(\nabla \cL)^i(\Theta_{n - 1})}^2, \\ \llabel{fm_Theta}
\textstyle \cM_n = \frac{(1 - \alpha) (\nabla \cL)(\Theta_{n - 1})}{1 - \alpha^n} + \frac{\alpha \m_n}{1 - \alpha^{n + 1}}, \qandq \Theta_n^i = \Theta_{n-1}^i - \gamma_n \bigl(\eps + \bigl[\frac{\abs{\bbM_n^i}}{1 - \beta^n}\bigr]^{\nicefrac{1}{2}}\bigr)^{- 1} \cM_n^i,
\end{gather}
and assume $\sup_{n \in \N} (\norm{\Theta_n} + \sum_{k = -1}^{1} (\gamma_n n^{\lrexpo})^k) < \infty$ \cfout. Then there exist $\vartheta \in \allowbreak \R^{\fd}$, $\rho \in (0, \allowbreak \infty)$ which satisfy for all $n \in \N$ that
\begin{equation}\llabel{result}
\textstyle \norm{(\nabla \cL)(\vartheta)} + \abs{\cL(\Theta_n) - \cL(\vartheta)} + \norm{\Theta_n - \vartheta}^{\rho} \le \rho \bigl[\sum_{j = 1}^{n} j^{- \lrexpo}\bigr]^{-1} \dott
\end{equation}
\end{athm}
\end{tcolorbox}
\end{samepage}
\end{savenotes}
\begin{aproof}
Throughout this proof for every $n \in \N$ let $p_n \in \R$, $\bbA_n \in \R^{\fd \times \fd}$ satisfy
\begin{equation}\llabel{bbA_p}
\textstyle p_n = \frac{\frac{\alpha}{1 - \alpha^{n + 1}}}{\frac{1 - \alpha}{1 - \alpha^n} + \frac{\alpha}{1 - \alpha^{n + 1}}} \qandq \bbA_n = \Bigl(\mathbbm{1}_{\{j\}}(i) \bigl(\eps + \bigl[\frac{\abs{\bbM_n^i}}{1 - \beta^n}\bigr]^{\nicefrac{1}{2}}\bigr)^{- 1} \bigl[\frac{1 - \alpha}{1 - \alpha^n} + \frac{\alpha}{1 - \alpha^{n + 1}}\bigr]\Bigr)_{(i, j) \in \{1, 2, \dots, \fd\}^2} \dott
\end{equation}
\startnewargseq
\argument{\lref{fm_Theta}; \lref{bbA_p}}{for all $n \in \N$ that
\begin{equation}\llabel{Theta_new}
\textstyle \Theta_{n} = \Theta_{n - 1} - \gamma_n \bbA_n \bigl[p_n \m_n + (1 - p_n) (\nabla \cL)(\Theta_{n - 1})\bigr] \dott
\end{equation}
}
\argument{the assumption that $\sup_{n \in \N} \norm{\Theta_n} < \infty$; the assumption that $\nabla \cL$ is locally Lipschitz continuous}{that $\sup_{n \in \N} \norm{(\nabla \cL)(\Theta_n)} < \infty$. \llabel{cG_bound}}
\argument{\lref{m_M}; \lref{cG_bound}; \cref{lemma:momentum_priori_bound_2} (applied for every $n \in \N$, $i \in \{1, 2, \dots, \fd\}$ with $A_{n - 1} \with \abs{(\nabla \cL)^i(\Theta_{n - 1})}^2$, $\alpha_n \with \beta$, $x_{n - 1} \with \abs{\bbM_{n - 1}^i}$ in the notation of \cref{lemma:momentum_priori_bound_2}); the triangle inequality}{that $\sup_{n \in \N} \norm{\bbM_n} < \infty$. \llabel{bbM_bound}}
\argument{the assumption that $\sup_{n \in \N} (\sum_{k = -1}^{1} (\gamma_n n^{\lrexpo})^k) < \infty$; \lref{bbA_p}; \lref{bbM_bound}; the fact that for all $n \in \N$ it holds that
\begin{equation}\llabel{alpha_beta_bound}
\textstyle 0 \le \alpha \le \frac{\alpha}{1 - \alpha^{n + 1}} \le \frac{\alpha}{1 - \alpha} \qqandqq 0 < 1 - \alpha \le \frac{1 - \alpha}{1 - \alpha^n} \le 1 \le \frac{1}{1 - \beta^n} \le \frac{1}{1 - \beta}
\end{equation}
}{that there exists $\lambda \in (0, \infty)$ which satisfies for all $n \in \N$ that $\gamma_n n^{\lrexpo} \bbA_n - \lambda \bbI_{\fd}$ is symmetric positive semi-definite and $\sup_{k \in \N} (\norm{\gamma_k k^{\lrexpo} \bbA_k} + \abs{p_k}) < \infty$. \llabel{bbA_bound}}
\argument{the assumption that $\sup_{n \in \N} \norm{\Theta_n} < \infty$; \lref{m_M}; \lref{Theta_new}; \lref{bbA_bound}; \cref{cor:gen_convergence_rate_Loja_gen_l_rate_Lipschitz_alpha} (applied for every $n \in \N$ with $\fd \allowbreak \with \fd$, $\alpha \allowbreak \with \alpha$, $\lambda \allowbreak \with \lambda$, $\lrexpo \allowbreak \with \lrexpo$, $\cL \allowbreak \with \cL$, $\Theta_n \allowbreak \with \Theta_n$, $p_n \allowbreak \with p_n$, $\gamma_n \allowbreak \with \gamma_n$, $\bbA_n \allowbreak \with \bbA_n$, $\m_n \allowbreak \with \m_n$, $\mu_n \allowbreak \with 0$, in the notation of \cref{cor:gen_convergence_rate_Loja_gen_l_rate_Lipschitz_alpha})}[verbs=e]{\lref{result}. }
\end{aproof}

\subsection{Nesterov accelerated adaptive moment estimation maximum (Nadamax)}
\label{subsec:Nadamax}

In the following result, \cref{cor:UGD_to_Nadamax}, we apply \cref{cor:gen_convergence_rate_Loja_gen_l_rate_Lipschitz_alpha} from \cref{subsec:convergence_UGD_KL_Lipschitz} above in the special situation where the optimization method is the \Nadamax\ optimizer \cite{DozatNadam2} (cf., \eg, \cite[Subsection~6.9.2]{JentzenBookDeepLearning2023}).

\cfclear
\begin{savenotes}
\begin{samepage}
\begin{tcolorbox}[colback=white!95!gray,
                  colframe=black,
                  boxrule=0.5pt,
                  sharp corners,
                  enhanced,
                 ]
\begin{athm}{cor}{cor:UGD_to_Nadamax}[\Nadamax]
Let $\fd \in \N$, $\alpha \in [0, 1)$, $\beta \in [0, 1]$, $\eps \in (0, \infty)$, $\lrexpo \in (\nicefrac{3}{4}, 1]$, let $\cL \colon \R^{\fd} \to \R$ be a\cfadd{def:Loja} \KL\ function, assume that $\nabla \cL$ is locally Lipschitz continuous, let $\gamma \colon \N \allowbreak \to (0, \allowbreak \infty)$, $\Theta = \allowbreak (\Theta^1, \allowbreak \dots, \allowbreak \Theta^{\fd}) \allowbreak \colon \allowbreak \N_0 \allowbreak \to \allowbreak \R^\fd$, $\m = (\m^1, \allowbreak \dots, \allowbreak \m^{\fd}) \colon \allowbreak \N_0 \allowbreak \to \R^\fd$, $\bbM = (\bbM^1, \allowbreak \dots, \allowbreak \bbM^{\fd}) \colon \allowbreak \N_0 \allowbreak \to \R^\fd$, and $\cM = (\cM^1, \allowbreak \dots, \allowbreak \cM^{\fd}) \colon \N \to \R^\fd$ satisfy for all $n \in \N$, $i \in \{1, 2, \dots, \fd\}$ that
\begin{gather}\llabel{m_M}
\textstyle \m_n = \alpha \m_{n - 1} + (1 - \alpha) (\nabla \cL)(\Theta_{n - 1}), \qquad \bbM_n^i = \max\{\beta \bbM_{n-1}^i, \abs{(\nabla \cL)^i(\Theta_{n - 1})}\}, \\ \llabel{fm_Theta}
\textstyle \cM_n = \frac{(1 - \alpha) (\nabla \cL)(\Theta_{n - 1})}{1 - \alpha^n} + \frac{\alpha \m_n}{1 - \alpha^{n + 1}}, \qqandqq \Theta_n^i = \Theta_{n-1}^i - \gamma_n [\eps + \abs{\bbM_n^i}]^{- 1} \cM_n^i,
\end{gather}
and assume $\sup_{n \in \N} (\norm{\Theta_n} + \sum_{k = -1}^{1} (\gamma_n n^{\lrexpo})^k) < \infty$ \cfout. Then there exist $\vartheta \in \allowbreak \R^{\fd}$, $\rho \in (0, \allowbreak \infty)$ which satisfy for all $n \in \N$ that
\begin{equation}\llabel{result}
\textstyle \norm{(\nabla \cL)(\vartheta)} + \abs{\cL(\Theta_n) - \cL(\vartheta)} + \norm{\Theta_n - \vartheta}^{\rho} \le \rho \bigl[\sum_{j = 1}^{n} j^{- \lrexpo}\bigr]^{-1} \dott
\end{equation}
\end{athm}
\end{tcolorbox}
\end{samepage}
\end{savenotes}
\begin{aproof}
Throughout this proof for every $n \in \N$ let $p_n \in \R$, $\bbA_n \in \R^{\fd \times \fd}$ satisfy
\begin{equation}\llabel{bbA_p}
\textstyle p_n = \frac{\frac{\alpha}{1 - \alpha^{n + 1}}}{\frac{1 - \alpha}{1 - \alpha^n} + \frac{\alpha}{1 - \alpha^{n + 1}}} \qandq \bbA_n = \bigl(\mathbbm{1}_{\{j\}}(i) [\eps + \abs{\bbM_n^i}]^{- 1} \bigl[\frac{1 - \alpha}{1 - \alpha^n} + \frac{\alpha}{1 - \alpha^{n + 1}}\bigr]\bigr)_{(i, j) \in \{1, 2, \dots, \fd\}^2} \dott
\end{equation}
\startnewargseq
\argument{\lref{fm_Theta}; \lref{bbA_p}}{for all $n \in \N$ that
\begin{equation}\llabel{Theta_new}
\textstyle \Theta_{n} = \Theta_{n - 1} - \gamma_n \bbA_n \bigl[p_n \m_n + (1 - p_n) (\nabla \cL)(\Theta_{n - 1})\bigr] \dott
\end{equation}
}
\argument{the assumption that $\sup_{n \in \N} \norm{\Theta_n} < \infty$; the assumption that $\nabla \cL$ is locally Lipschitz continuous}{that $\sup_{n \in \N} \norm{(\nabla \cL)(\Theta_n)} < \infty$. \llabel{cG_bound}}
\argument{\lref{m_M}; \lref{cG_bound}; the induction; the triangle inequality}{that $\sup_{n \in \N} \norm{\bbM_n} < \infty$. \llabel{bbM_bound}}
\argument{the assumption that $\sup_{n \in \N} (\sum_{k = -1}^{1} (\gamma_n n^{\lrexpo})^k) \allowbreak < \infty$; \lref{bbA_p}; \lref{bbM_bound}; the fact that for all $n \in \N$ it holds that
\begin{equation}\llabel{alpha_beta_bound}
\textstyle 0 \le \alpha \le \frac{\alpha}{1 - \alpha^{n + 1}} \le \frac{\alpha}{1 - \alpha} < \infty \qqandqq 0 < 1 - \alpha \le \frac{1 - \alpha}{1 - \alpha^n} \le 1
\end{equation}
}{that there exists $\lambda \in (0, \infty)$ which satisfies for all $n \in \N$ that $\gamma_n n^{\lrexpo} \bbA_n - \lambda \bbI_{\fd}$ is symmetric positive semi-definite and $\sup_{k \in \N} (\norm{\gamma_k k^{\lrexpo} \bbA_k} + \abs{p_k}) < \infty$. \llabel{bbA_bound}}
\argument{the assumption that $\limsup_{n \to \infty} \norm{\Theta_n} < \infty$; \lref{m_M}; \lref{Theta_new}; \lref{bbA_bound}; \cref{cor:gen_convergence_rate_Loja_gen_l_rate_Lipschitz_alpha} (applied for every $n \in \N$ with $\fd \with \fd$, $\alpha \with \alpha$, $\lambda \with \lambda$, $\lrexpo \with \lrexpo$, $\cL \with \cL$, $\Theta_n \with \Theta_n$, $p_n \with p_n$, $\gamma_n \with \gamma_n$, $\bbA_n \with \bbA_n$, $\m_n \with \m_n$, $\mu_n \with 0$ in the notation of \cref{cor:gen_convergence_rate_Loja_gen_l_rate_Lipschitz_alpha})}[verbs=e]{\lref{result}. }
\end{aproof}

\subsection{Adaptive Nesterov momentum (Adan)}
\label{subsec:Adan}

In the following result, \cref{cor:UGD_to_vanilla_Adan}, we apply \cref{cor:gen_convergence_rate_Loja_gen_l_rate_Lipschitz_alpha} from \cref{subsec:convergence_UGD_KL_Lipschitz} above in the special situation where the optimization method is the (vanilla) \Adan\ optimizer \cite[Subsection~3.2]{Adan}.

\cfclear
\begin{savenotes}
\begin{samepage}
\begin{tcolorbox}[colback=white!95!gray,
                  colframe=black,
                  boxrule=0.5pt,
                  sharp corners,
                  enhanced,
                 ]
\begin{athm}{cor}{cor:UGD_to_vanilla_Adan}[\resname{\Adan}]
Let $\fd \in \N$, $\beta_1 \in [0, 1)$, $\beta_2 \in [0, 1]$, $\eps \in (0, \infty)$, $\lrexpo \in (\nicefrac{3}{4}, \allowbreak 1]$, let $\cL \colon \allowbreak \R^{\fd} \allowbreak \to \allowbreak \R$ be a\cfadd{def:Loja} \KL\ function, assume that $\nabla \cL$ is locally Lipschitz continuous, let $\gamma \colon \N \allowbreak \to (0, \allowbreak \infty)$ and $\Theta\Index{k}{} \allowbreak = (\Theta\Index{k}[1]{}, \allowbreak \dots, \allowbreak \Theta\Index{k}[\fd]{}) \allowbreak \colon \allowbreak \N_0 \allowbreak \to \allowbreak \R^\fd$, $k \in \N_0$, satisfy for all $n \in \N$, $k \in \allowbreak \{1, \allowbreak 2\}$, $i \in \allowbreak \{1, \allowbreak 2, \allowbreak \dots, \allowbreak \fd\}$ that
\begin{gather}\llabel{m_M}
\textstyle \Theta\Index{k}[i]{n + 1} = \beta_k \Theta\Index{k}[i]{n} + (1 - \beta_k) \bigl[(\nabla \cL)^i(\Theta\Index{0}{n}) + \beta_1 [(\nabla \cL)^i(\Theta\Index{0}{n}) - (\nabla \cL)^i(\Theta\Index{0}{n - 1})]\bigr]^k \\ \llabel{Theta}
\textstyle \andqq \Theta\Index{0}[i]{n} = \Theta\Index{0}[i]{n - 1} - \gamma_n \Theta\Index{1}[i]{n} \bigl[\eps + \abs{\Theta\Index{2}[i]{n}}^{\nicefrac{1}{2}}\bigr]^{-1},
\end{gather}
and assume $\sup_{n \in \N} (\norm{\Theta\Index{0}{n}} + \sum_{k = -1}^{1} (\gamma_n n^{\lrexpo})^k) < \infty$ \cfout. Then there exist $\vartheta  \allowbreak \in \allowbreak \R^{\fd}$, $\rho \in (0, \allowbreak \infty)$ which satisfy for all $n \in \N$ that
\begin{equation}\llabel{result}
\textstyle  \norm{(\nabla \cL)(\vartheta)} + \abs{\cL(\Theta\Index{0}{n}) - \cL(\vartheta)} + \norm{\Theta\Index{0}{n} - \vartheta}^{\rho} \le \rho \bigl[\sum_{j = 1}^{n} j^{- \lrexpo}\bigr]^{-1} \dott
\end{equation}
\end{athm}
\end{tcolorbox}
\end{samepage}
\end{savenotes}
\begin{aproof}
Throughout this proof for every $n \in \N$ let $\mu_n \in \R^{\fd}$, $\bbA_n \in \R^{\fd \times \fd}$ satisfy
\begin{equation}\llabel{bbA_mu}
\textstyle \mu_{n + 1} = (\gamma_{n + 1})^{- 1} \beta_1 [(\nabla \cL)(\Theta\Index{0}{n}) - (\nabla \cL)(\Theta\Index{0}{n - 1})] \qandq \bbA_n = \Bigl(\frac{\mathbbm{1}_{\{j\}}(i)}{\eps + \abs{\Theta\Index{2}[i]{n}}^{\nicefrac{1}{2}}}\Bigr)_{(i, j) \in \{1, 2, \dots, \fd\}^2} \dott
\end{equation}
\argument{\lref{m_M}; \lref{Theta}; \lref{bbA_mu}}{for all $n \in \N$ that
\begin{equation}\llabel{m_Theta_new}
\textstyle \Theta\Index{1}{n + 1} = \beta_1 \Theta\Index{1}{n} + (1 - \beta_1) \bigl[(\nabla \cL)(\Theta\Index{0}{n}) + \gamma_{n + 1} \mu_{n + 1}\bigr] \qandq \Theta\Index{0}{n} = \Theta\Index{0}{n - 1} - \gamma_n \bbA_n \Theta\Index{1}{n} \dott
\end{equation}
}
\argument{the assumption that $\sup_{n \in \N} \norm{\Theta\Index{0}{n}} < \infty$; the assumption that $\nabla \cL$ is locally Lipschitz continuous; the triangle inequality}{for all $\alpha \in \R$ that
\begin{equation}\llabel{G_bound}
\textstyle \sup_{n \in \N} \norm{(\nabla \cL)(\Theta\Index{0}{n}) + \alpha [(\nabla \cL)(\Theta\Index{0}{n}) - (\nabla \cL)(\Theta\Index{0}{n - 1})]} < \infty \dott
\end{equation}
}
\argument{\lref{m_M}; \lref{bbA_mu}; \lref{G_bound}; \cref{lemma:momentum_priori_bound_2} (applied for every $n \in \N$ with $\alpha_n \with \beta_1$, $A_{n - 1} \allowbreak \with \allowbreak \norm{(\nabla \cL)(\Theta\Index{0}{n}) + \gamma_{n + 1} \mu_{n + 1}}$, $x_{n - 1} \allowbreak \with \allowbreak \norm{\Theta\Index{1}{n}}$ in the notation of \cref{lemma:momentum_priori_bound_2}); \cref{lemma:momentum_priori_bound_2} (applied for every $n \in \N$, $i \in \{1, \allowbreak 2, \allowbreak \dots, \allowbreak \fd\}$ with $\alpha_n \with \allowbreak \beta_2$, $A_{n - 1} \allowbreak \with \allowbreak \abs{(\nabla \cL)^i(\Theta\Index{0}{n}) + \beta_2 [(\nabla \cL)^i(\Theta\Index{0}{n}) - (\nabla \cL)^i(\Theta\Index{0}{n - 1})]}^2$, $x_{n - 1} \with \allowbreak \abs{\Theta\Index{2}[i]{n}}$ in the notation of \cref{lemma:momentum_priori_bound_2}); the triangle inequality}{that
\begin{equation}\llabel{m_M_Bound}
\textstyle \sup_{n \in \N} (\norm{\Theta\Index{1}{n}} + \norm{\Theta\Index{2}{n}}) < \infty \dott
\end{equation}
}
\argument{the assumption that $\sup_{n \in \N} (\sum_{k = -1}^{1} (\gamma_n n^{\lrexpo})^k) < \infty$; \lref{bbA_mu}; \lref{m_M_Bound}}{that there exists $\lambda \in (0, \infty)$ which satisfies for all $n \in \N$ that $\gamma_n n^{\lrexpo} \bbA_n - \lambda \bbI_{\fd}$ is symmetric positive semi-definite and $\sup_{k \in \N} (\norm{\gamma_k k^{\lrexpo} \bbA_k} + \norm{\bbA_k}) < \infty$. \llabel{bbA_bound}}
\argument{\lref{bbA_mu}; \lref{m_Theta_new}; \lref{m_M_Bound}; \lref{bbA_bound}; the assumption that $\nabla \cL$ is locally Lipschitz continuous; the assumption that $\sup_{n \in \N} (\norm{\Theta\Index{0}{n}} + \sum_{k = -1}^{1} (\gamma_n n^{\lrexpo})^k) < \infty$; the triangle inequality}{that there exists $L \in (0, \infty)$ which satisfies for all $n \in \N$ that
\begin{multline}\llabel{mu_Bound}
\textstyle \norm{\mu_{n + 1}} \le (\gamma_{n + 1})^{- 1} \beta_1 L \norm{\Theta\Index{0}{n} - \Theta\Index{0}{n - 1}} = (\gamma_{n + 1})^{- 1} \beta_1 L \gamma_n \norm{\bbA_n \Theta\Index{1}{n}} \\
\textstyle = \beta_1 L (\gamma_{n + 1} (n + 1)^{\lrexpo})^{- 1} (\gamma_n n^{\lrexpo}) [1 + \frac{1}{n}]^{\lrexpo} \norm{\bbA_n \Theta\Index{1}{n}} \le L^2 \dott
\end{multline}
}
\argument{\lref{m_Theta_new}; \lref{mu_Bound}; the assumption that $\sup_{n \in \N} (\norm{\Theta\Index{0}{n}} + \sum_{k = -1}^{1} (\gamma_n n^{\lrexpo})^k) < \infty$; the fact that for all $n \in \N$ it holds that $\gamma_n n^{\lrexpo} \bbA_n - \lambda \bbI_{\fd}$ is symmetric positive semi-definite; the fact that $\sup_{k \in \N} \norm{\gamma_k k^{\lrexpo} \bbA_k} < \infty$; \cref{cor:gen_convergence_rate_Loja_gen_l_rate_Lipschitz_alpha} (applied for every $n \in \N$ with $\fd \with \fd$, $\alpha \with \beta_1$, $\lambda \with \lambda$, $\lrexpo \with \lrexpo$, $\cL \with \cL$, $\Theta_n \with \Theta\Index{0}{n}$, $p_n \with 1$, $\gamma_n \with \gamma_n$, $\bbA_n \with \bbA_n$, $\m_n \with \Theta\Index{1}{n}$, $\mu_{n + 1} \with \mu_{n + 1}$ in the notation of \cref{cor:gen_convergence_rate_Loja_gen_l_rate_Lipschitz_alpha})}[verbs=e]{\lref{result}. }
\end{aproof}

\subsection{Adaptive belief (AdaBelief)}
\label{subsec:AdaBelief}

In the following result, \cref{cor:UGD_to_AdaBelief_new}, we apply \cref{cor:gen_convergence_rate_Loja_gen_l_rate_Lipschitz_alpha} from \cref{subsec:convergence_UGD_KL_Lipschitz} above in the special situation where the optimization method is the \AdaBelief\ optimizer \cite{ZhuangAdaBelief2020}.

\cfclear
\begin{savenotes}
\begin{samepage}
\begin{tcolorbox}[colback=white!95!gray,
                  colframe=black,
                  boxrule=0.5pt,
                  sharp corners,
                  enhanced,
                 ]
\begin{athm}{cor}{cor:UGD_to_AdaBelief_new}[\AdaBelief]
Let $\fd \in \N$, $\beta_1, \beta_2 \in [0, 1)$, $\eps, \epsilon \in (0, \infty)$, $\lrexpo \in (\nicefrac{3}{4}, 1]$, let $\cL \colon \R^{\fd} \to \R$ be a\cfadd{def:Loja} \KL\ function, assume that $\nabla \cL$ is locally Lipschitz continuous, let $\gamma \colon \N \allowbreak \to (0, \allowbreak \infty)$ and $\Theta\Index{k}{} \allowbreak = (\Theta\Index{k}[1]{}, \allowbreak \dots, \allowbreak \Theta\Index{k}[\fd]{}) \allowbreak \colon \allowbreak \N_0 \allowbreak \to \allowbreak \R^\fd$, $k \in \N_0$, satisfy for all $n \in \N$, $k \in \allowbreak \{1, \allowbreak 2\}$, $i \in \allowbreak \{1, \allowbreak 2, \allowbreak \dots, \allowbreak \fd\}$ that
\begin{gather}\llabel{m_M}
\textstyle \Theta\Index{k}[i]{n} = \beta_k \Theta\Index{k}[i]{n - 1} + (1 - \beta_k) \bigl[(\nabla \cL)^i(\Theta\Index{0}{n - 1}) - \Theta\Index{1}[i]{n} \mathbbm{1}_{\{2\}}(k)\bigr]^k \\ \llabel{Theta}
\textstyle \andq \Theta\Index{0}[i]{n} = \Theta\Index{0}[i]{n - 1} - \gamma_n \bigl[\frac{\Theta\Index{1}[i]{n}}{1 - (\beta_1)^n}\bigr] \Bigl[\eps + \bigl[\frac{\epsilon + \abs{\Theta\Index{2}[i]{n}}}{1 - (\beta_2)^n}\bigr]^{\nicefrac{1}{2}}\Bigr]^{-1},
\end{gather}
and assume $\sup_{n \in \N} (\norm{\Theta\Index{0}{n}} + \sum_{k = -1}^{1} (\gamma_n n^{\lrexpo})^k) < \infty$ \cfout. Then there exist $\vartheta  \allowbreak \in \allowbreak \R^{\fd}$, $\rho \in (0, \allowbreak \infty)$ which satisfy for all $n \in \N$ that
\begin{equation}\llabel{result}
\textstyle  \norm{(\nabla \cL)(\vartheta)} + \abs{\cL(\Theta\Index{0}{n}) - \cL(\vartheta)} + \norm{\Theta\Index{0}{n} - \vartheta}^{\rho} \le \rho \bigl[\sum_{j = 1}^{n} j^{- \lrexpo}\bigr]^{-1} \dott
\end{equation}
\end{athm}
\end{tcolorbox}
\end{samepage}
\end{savenotes}
\begin{aproof}
Throughout this proof for every $n \in \N$ let $\bbA_n \in \R^{\fd \times \fd}$ satisfy
\begin{equation}\llabel{bbA}
\textstyle \bbA_n = \Bigl(\frac{\mathbbm{1}_{\{j\}}(i)}{1 - (\beta_1)^n} \bigl(\eps + \bigl[\frac{\epsilon + \abs{\Theta\Index{2}[i]{n}}}{1 - (\beta_2)^n}\bigr]^{\nicefrac{1}{2}}\bigr)^{- 1}\Bigr)_{(i, j) \in \{1, 2, \dots, \fd\}^2} \dott
\end{equation}
\startnewargseq
\argument{\lref{Theta}; \lref{bbA}}{for all $n \in \N$ that
\begin{equation}\llabel{Theta_new}
\textstyle \Theta\Index{0}{n} = \Theta\Index{0}{n - 1} - \gamma_n \bbA_n \Theta\Index{1}{n} \dott
\end{equation}
}
\argument{the assumption that $\sup_{n \in \N} \norm{\Theta\Index{0}{n}} < \infty$; the assumption that $\nabla \cL$ is locally Lipschitz continuous}{that $\sup_{n \in \N} \norm{(\nabla \cL)(\Theta\Index{0}{n})} < \infty$. \llabel{cG_bound}}
\argument{\lref{m_M}; \lref{cG_bound}; \cref{lemma:momentum_priori_bound_2} (applied for every $n \in \N$ with $\alpha_n \with \beta_1$, $A_{n - 1} \with \norm{(\nabla \cL)(\Theta\Index{0}{n - 1})}$, $x_{n - 1} \allowbreak \with \allowbreak \norm{\Theta\Index{1}{n - 1}}$ in the notation of \cref{lemma:momentum_priori_bound_2}); the triangle inequality}[verbs=d]{that $\sup_{n \in \N} \norm{\Theta\Index{1}{n}} < \infty$. \llabel{m_bound}}
\argument{the fact that $\sup_{n \in \N} \norm{(\nabla \cL)(\Theta\Index{0}{n})} < \infty$; \lref{m_M}; \lref{m_bound}; \cref{lemma:momentum_priori_bound_2} (applied for every $n \in \N$, $i \in \{1, 2, \dots, \fd\}$ with $\alpha_n \with \beta_2$, $A_{n - 1} \with \abs{(\nabla \cL)^i(\Theta\Index{0}{n - 1}) - \Theta\Index{1}[i]{n}}^2$, $x_{n - 1} \allowbreak \with \allowbreak \abs{\Theta\Index{2}[i]{n - 1}}$ in the notation of \cref{lemma:momentum_priori_bound_2}); the triangle inequality}{that $\sup_{n \in \N} \norm{\Theta\Index{2}{n}} < \infty$. \llabel{bbM_bound}}
\argument{the assumption that $\sup_{n \in \N} (\sum_{k = -1}^{1} (\gamma_n n^{\lrexpo})^k) < \infty$; \lref{bbA}; \lref{bbM_bound}; the fact that for all $n \in \N$ it holds that
\begin{equation}\llabel{alpha_beta_bound}
\textstyle 1 \le \frac{1}{1 - (\beta_1)^n} \le \frac{1}{1 - \beta_1} < \infty \qqandqq 1 \le \frac{1}{1 - (\beta_2)^n} \le \frac{1}{1 - \beta_2} < \infty
\end{equation}
}{that there exists $\lambda \in (0, \infty)$ which satisfies for all $n \in \N$ that $\gamma_n n^{\lrexpo} \bbA_n - \lambda \bbI_{\fd}$ is symmetric positive semi-definite and $\sup_{k \in \N} \norm{\gamma_k k^{\lrexpo} \bbA_k} < \infty$. \llabel{bbA_bound}}
\argument{the assumption that $\sup_{n \in \N} \norm{\Theta\Index{0}{n}} < \infty$; \lref{m_M}; \lref{Theta_new}; \lref{bbA_bound}; \cref{cor:gen_convergence_rate_Loja_gen_l_rate_Lipschitz_alpha} (applied for every $n \in \N$ with $\fd \with \fd$, $\alpha \with \beta_1$, $\lambda \with \lambda$, $\lrexpo \with \lrexpo$, $\cL \with \cL$, $\Theta_n \with \Theta\Index{0}{n}$, $p_n \with 1$, $\gamma_n \with \gamma_n$, $\bbA_n \with \bbA_n$, $\m_n \with \Theta\Index{1}{n}$, $\mu_n \with 0$ in the notation of \cref{cor:gen_convergence_rate_Loja_gen_l_rate_Lipschitz_alpha})}[verbs=e]{\lref{result}. }
\end{aproof}

\subsection{AMSGrad}
\label{subsec:AMSGrad}

In the following result, \cref{cor:UGD_to_AMSGrad}, we apply \cref{cor:gen_convergence_rate_Loja_gen_l_rate_Lipschitz_alpha} from \cref{subsec:convergence_UGD_KL_Lipschitz} above in the special situation where the optimization method is the AMSGrad optimizer \cite{ReddiKale2019} (cf., \eg, \cite[Subsection~6.13]{JentzenBookDeepLearning2023}).

\cfclear
\begin{savenotes}
\begin{samepage}
\begin{tcolorbox}[colback=white!95!gray,
                  colframe=black,
                  boxrule=0.5pt,
                  sharp corners,
                  enhanced,
                 ]
\begin{athm}{cor}{cor:UGD_to_AMSGrad}[\resname{AMSGrad}]
Let $\fd \in \N$, $\alpha \in [0, 1)$, $\beta \in [0, 1]$, $\eps \in (0, \infty)$, $\lrexpo \in (\nicefrac{3}{4}, 1]$, let $\cL \colon \R^{\fd} \to \R$ be a\cfadd{def:Loja} \KL\ function, assume that $\nabla \cL$ is locally Lipschitz continuous, let $\gamma \colon \N \allowbreak \to (0, \allowbreak \infty)$, $\Theta = (\Theta^1, \dots, \Theta^{\fd}) \allowbreak \colon \allowbreak \N_0 \allowbreak \to \allowbreak \R^\fd$, $\m = (\m^1, \allowbreak \dots, \allowbreak \m^{\fd}) \colon \allowbreak \N_0 \allowbreak \to \R^\fd$, $\bbM = (\bbM^1, \allowbreak \dots, \allowbreak \bbM^{\fd}) \colon \allowbreak \N_0 \allowbreak \to \R^\fd$, and $\fM = (\fM^1, \allowbreak \dots, \allowbreak \fM^{\fd}) \colon \N_0 \to \R^\fd$ satisfy for all $n \in \N$, $i \in \{1, 2, \dots, \fd\}$ that
\begin{gather}\llabel{m_M}
\textstyle \m_n = \alpha \m_{n - 1} + (1 - \alpha) (\nabla \cL)(\Theta_{n - 1}), \qquad \bbM_n^i = \beta \bbM_{n - 1}^i + (1 - \beta) \abs{(\nabla \cL)^i(\Theta_{n - 1})}^2, \\ \llabel{fM_Theta}
\textstyle \fM_n^i = \max\{\fM_{n - 1}^i, \abs{\bbM_n^i}\}, \qqandqq \Theta_n^i = \Theta_{n-1}^i - \gamma_n \m_n^i \bigl[\eps + \abs{\fM_n^i}^{\nicefrac{1}{2}}\bigr]^{- 1} \dott
\end{gather}
and assume $\sup_{n \in \N} (\norm{\Theta_n} + \sum_{k = -1}^{1} (\gamma_n n^{\lrexpo})^k) < \infty$ \cfout. Then there exist $\vartheta \in \allowbreak \R^{\fd}$, $\rho \in (0, \allowbreak \infty)$ which satisfy for all $n \in \N$ that
\begin{equation}\llabel{result}
\textstyle \norm{(\nabla \cL)(\vartheta)} + \abs{\cL(\Theta_n) - \cL(\vartheta)} + \norm{\Theta_n - \vartheta}^{\rho} \le \rho \bigl[\sum_{j = 1}^{n} j^{- \lrexpo}\bigr]^{- 1} \dott
\end{equation}
\end{athm}
\end{tcolorbox}
\end{samepage}
\end{savenotes}
\begin{aproof}
Throughout this proof for every $n \in \N$ let $\bbA_n \in \R^{\fd \times \fd}$ satisfy
\begin{equation}\llabel{bbA}
\textstyle \bbA_n = \bigl(\bigl[\eps + \abs{\fM_n^i}^{\nicefrac{1}{2}}\bigr]^{- 1} \mathbbm{1}_{\{j\}}(i)\bigr)_{(i, j) \in \{1, 2, \dots, \fd\}^2} \dott
\end{equation}
\argument{\lref{fM_Theta}; \lref{bbA}}{for all $n \in \N$ that
\begin{equation}\llabel{Theta_new}
\textstyle \Theta_n = \Theta_{n - 1} - \gamma_n \bbA_n \m_n \dott
\end{equation}
}
\argument{the assumption that $\sup_{n \in \N} \norm{\Theta_n} < \infty$; the assumption that $\nabla \cL$ is locally Lipschitz continuous}{that $\sup_{n \in \N} \norm{(\nabla \cL)(\Theta_n)} < \infty$. \llabel{cG_bound}}
\argument{\cref{lemma:momentum_priori_bound_2} (applied for every $n \in \N$, $i \in \{1, 2, \dots, \fd\}$ with $\alpha_n \with \beta$, $A_{n - 1} \with \abs{(\nabla \cL)^i(\Theta_{n - 1})}^2$, $x_{n - 1} \with \abs{\bbM_{n - 1}^i}$ in the notation of \cref{lemma:momentum_priori_bound_2}); \lref{m_M}; \lref{cG_bound}; the triangle inequality}{that $\sup_{n \in \N} \norm{\bbM_n} < \infty$. \llabel{bbM_bound}}
\argument{\lref{fM_Theta}; \lref{bbM_bound}; the triangle inequality; the induction}{that $\sup_{n \in \N} \norm{\fM_n} < \infty$. \llabel{fM_bound}}
\argument{the assumption that $\sup_{n \in \N} (\sum_{k = -1}^{1} (\gamma_n n^{\lrexpo})^k) < \infty$; \lref{bbA}; \lref{fM_bound}}{that there exists $\lambda \in (0, \infty)$ which satisfies for all $n \in \N$ that $\gamma_n n^{\lrexpo} \bbA_n - \lambda \bbI_{\fd}$ is symmetric positive semi-definite and $\sup_{k \in \N} \norm{\gamma_k k^{\lrexpo} \bbA_k} < \infty$. \llabel{bbA_bound}}
\argument{\cref{cor:gen_convergence_rate_Loja_gen_l_rate_Lipschitz_alpha} (applied for every $n \in \N$ with $\fd \with \fd$, $\alpha \with \alpha$, $\lambda \with \lambda$, $\lrexpo \with \lrexpo$, $\cL \with \cL$, $\Theta_n \with \Theta_n$, $p_n \with 1$, $\gamma_n \with \gamma_n$, $\bbA_n \with \bbA_n$, $\m_n \with \m_n$, $\mu_n \with 0$ in the notation of \cref{cor:gen_convergence_rate_Loja_gen_l_rate_Lipschitz_alpha}); \lref{m_M}; \lref{Theta_new}; \lref{bbA_bound}; the assumption that $\sup_{n \in \N} \norm{\Theta_n} < \infty$}[verbs=e]{\lref{result}. }
\end{aproof}

\subsection{Yogi}
\label{subsec:Yogi}

In the following result, \cref{cor:UGD_to_Yogi}, we apply \cref{cor:gen_convergence_rate_Loja_gen_l_rate_Lipschitz_alpha} from \cref{subsec:convergence_UGD_KL_Lipschitz} above in the special situation where the optimization method is the Yogi optimizer \cite{Yogi_NEURIPS2018}.

\cfclear
\begin{savenotes}
\begin{samepage}
\begin{tcolorbox}[colback=white!95!gray,
                  colframe=black,
                  boxrule=0.5pt,
                  sharp corners,
                  enhanced,
                 ]
\begin{athm}{cor}{cor:UGD_to_Yogi}[\resname{Yogi}]
Let $\fd \in \N$, $\alpha \in [0, 1)$, $\beta \in [0, 1]$, $\eps \in (0, \infty)$, $\lrexpo \in (\nicefrac{3}{4}, 1]$, let $\cL \colon \R^{\fd} \to \R$ be a\cfadd{def:Loja} \KL\ function, assume that $\nabla \cL$ is locally Lipschitz continuous, let $\gamma \colon \N \allowbreak \to (0, \allowbreak \infty)$, $\Theta = (\Theta^1, \dots, \Theta^{\fd}) \allowbreak \colon \allowbreak \N_0 \allowbreak \to \allowbreak \R^\fd$, $\m = (\m^1, \allowbreak \dots, \allowbreak \m^{\fd}) \colon \allowbreak \N_0 \allowbreak \to \R^\fd$, and $\bbM = (\bbM^1, \allowbreak \dots, \allowbreak \bbM^{\fd}) \colon \allowbreak \N_0 \allowbreak \to \R^\fd$ satisfy for all $n \in \N$, $i \in \{1, 2, \dots, \fd\}$ that
\begin{gather}\llabel{m}
\textstyle \m_n = \alpha \m_{n - 1} + (1 - \alpha) (\nabla \cL)(\Theta_{n - 1}), \\ \llabel{M}
\textstyle \bbM_n^i = \bbM_{n - 1}^i - (1 - \beta) \abs{(\nabla \cL)^i(\Theta_{n - 1})}^2 {\operatorfont{sign}}(\bbM_{n - 1}^i - \abs{(\nabla \cL)^i(\Theta_{n - 1})}^2), \\ \llabel{Theta}
\textstyle \andqq \Theta_n^i = \Theta_{n - 1}^i - \gamma_n \m_n^i \bigl[\eps + \abs{\bbM_n^i}^{\nicefrac{1}{2}}\bigr]^{- 1},
\end{gather}
and assume $\sup_{n \in \N} (\norm{\Theta_n} + \sum_{k = -1}^{1} (\gamma_n n^{\lrexpo})^k) < \infty$ \cfout. Then there exist $\vartheta \in \allowbreak \R^{\fd}$, $\rho \in (0, \allowbreak \infty)$ which satisfy for all $n \in \N$ that
\begin{equation}\llabel{result}
\textstyle \norm{(\nabla \cL)(\vartheta)} + \abs{\cL(\Theta_n) - \cL(\vartheta)} + \norm{\Theta_n - \vartheta}^{\rho} \le \rho \bigl[\sum_{j = 1}^{n} j^{- \lrexpo}\bigr]^{- 1} \dott
\end{equation}
\end{athm}
\end{tcolorbox}
\end{samepage}
\end{savenotes}
\begin{aproof}
Throughout this proof for every $n \in \N$ let $\bbA_n \in \R^{\fd \times \fd}$ satisfy
\begin{equation}\llabel{bbA}
\textstyle \bbA_n = \bigl(\bigl[\eps + \abs{\bbM_n^i}^{\nicefrac{1}{2}}\bigr]^{-1} \mathbbm{1}_{\{j\}}(i)\bigr)_{(i, j) \in \{1, 2, \dots, \fd\}^2} \dott
\end{equation}
\argument{\lref{Theta}; \lref{bbA}}{for all $n \in \N$ that
\begin{equation}\llabel{Theta_new}
\textstyle \Theta_n = \Theta_{n - 1} - \gamma_n \bbA_n \m_n \dott
\end{equation}
}
\argument{\lref{M}}{for all $n \in \N$, $i \in \{1, 2, \dots, \fd\}$ that
\begin{equation}\llabel{M_G_ineq}
\textstyle \abs{\bbM_n^i} \le \max\{\abs{\bbM_{n - 1}^i}, 2 \abs{(\nabla \cL)^i(\Theta_{n - 1})}^2\} \dott
\end{equation}
}
\argument{the assumption that $\sup_{n \in \N} \norm{\Theta_n} < \infty$; the assumption that $\nabla \cL$ is locally Lipschitz continuous}{that $\sup_{n \in \N} \norm{(\nabla \cL)(\Theta_{n})} < \infty$. \llabel{G_bound}}
\argument{\lref{M_G_ineq}; \lref{G_bound}; the induction}{that $\sup_{n \in \N} \norm{\bbM_n} < \infty$. \llabel{M_bound}}
\argument{the assumption that $\sup_{n \in \N} (\sum_{k = -1}^{1} (\gamma_n n^{\lrexpo})^k) < \infty$; \lref{bbA}; \lref{M_bound}}{that there exists $\lambda \in (0, \infty)$ which satisfies for all $n \in \N$ that $\gamma_n n^{\lrexpo} \bbA_n - \lambda \bbI_{\fd}$ is symmetric positive semi-definite and $\sup_{k \in \N} \norm{\gamma_k k^{\lrexpo} \bbA_k} < \infty$. \llabel{bbA_bound}}
\argument{the assumption that $\sup_{n \in \N} \norm{\Theta_n} < \infty$; \lref{m}; \lref{Theta_new}; \lref{bbA_bound}; \cref{cor:gen_convergence_rate_Loja_gen_l_rate_Lipschitz_alpha} (applied for every $n \in \N$ with $\fd \with \fd$, $\alpha \with \alpha$, $\lambda \with \lambda$, $\lrexpo \with \lrexpo$, $\cL \with \cL$, $\Theta_n \with \Theta_n$, $p_n \with 1$, $\gamma_n \with \gamma_n$, $\bbA_n \with \bbA_n$, $\m_n \with \m_n$, $\mu_n \with 0$ in the notation of \cref{cor:gen_convergence_rate_Loja_gen_l_rate_Lipschitz_alpha})}[verbs=e]{\lref{result}. }
\end{aproof}

\subsection{Explicit midpoint}

In the following result, \cref{cor:UGD_to_Explicit_midpoint_GD}, we apply \cref{cor:gen_convergence_rate_Loja_gen_l_rate_Lipschitz_alpha} from \cref{subsec:convergence_UGD_KL_Lipschitz} above in the special situation where the optimization method is the explicit midpoint method (cf., \eg, \cite[Subsection~6.2]{JentzenBookDeepLearning2023}).

\cfclear
\begin{savenotes}
\begin{samepage}
\begin{tcolorbox}[colback=white!95!gray,
                  colframe=black,
                  boxrule=0.5pt,
                  sharp corners,
                  enhanced,
                 ]
\begin{athm}{cor}{cor:UGD_to_Explicit_midpoint_GD}[\resname{Explicit midpoint}]
Let $\fd \in \N$, $\lrexpo \in (\nicefrac{3}{4}, 1]$, let $\cL \colon \R^{\fd} \to \R$ be a\cfadd{def:Loja} \KL\ function, assume that $\nabla \cL$ is locally Lipschitz continuous, let $\gamma \colon \N \allowbreak \to (0, \allowbreak \infty)$ and $\Theta \allowbreak \colon \allowbreak \N_0 \allowbreak \to \allowbreak \R^\fd$ satisfy for all $n \in \N$ that
\begin{equation}\llabel{Theta}
\textstyle \Theta_n = \Theta_{n - 1} - \gamma_n (\nabla \cL)\bigl(\Theta_{n - 1} - 2^{- 1} \gamma_n (\nabla \cL)(\Theta_{n - 1})\bigr),
\end{equation}
and assume $\sup_{n \in \N} (\norm{\Theta_n} + \sum_{k = -1}^{1} (\gamma_n n^{\lrexpo})^k) < \infty$ \cfout. Then there exist $\vartheta \in \allowbreak \R^{\fd}$, $\rho \in (0, \allowbreak \infty)$ which satisfy for all $n \in \N$ that
\begin{equation}\llabel{result}
\textstyle \norm{(\nabla \cL)(\vartheta)} + \abs{\cL(\Theta_n) - \cL(\vartheta)} + \norm{\Theta_n - \vartheta}^{\rho} \le \rho \bigl[\sum_{j = 1}^{n} j^{- \lrexpo}\bigr]^{- 1} \dott
\end{equation}
\end{athm}
\end{tcolorbox}
\end{samepage}
\end{savenotes}
\begin{aproof}
Throughout this proof for every $n \in \N$ let $\mu_n \in \R^{\fd}$ satisfy
\begin{equation}\llabel{mu}
\textstyle \mu_n = (\gamma_n)^{- 1} \bigl[(\nabla \cL)\bigl(\Theta_{n - 1} - 2^{- 1} \gamma_n (\nabla \cL)(\Theta_{n - 1})\bigr) - (\nabla \cL)(\Theta_{n - 1})\bigr] \dott
\end{equation}
\argument{\lref{Theta}; \lref{mu}}{for all $n \in \N$ that
\begin{equation}\llabel{Theta_new}
\textstyle \Theta_n = \Theta_{n - 1} - \gamma_n \bigl[(\nabla \cL)(\Theta_{n - 1}) + \gamma_n \mu_n\bigr] \dott
\end{equation}
}
\argument{the fact that $\nabla \cL$ is locally Lipschitz continuous; the triangle inequality; the assumption that $\sup_{n \in \N} \norm{\Theta_n} \allowbreak < \infty$}{that there exist $B, L \in (0, \infty)$ which satisfy for all $n \in \N$, $u, w \in \{x \in \R^{\fd} \colon \norm{x} \le B\}$ that
\begin{equation}\llabel{B_L}
\textstyle \norm{\Theta_{n - 1} - \frac{\gamma_n}{2} (\nabla \cL)(\Theta_{n - 1})} + \norm{\Theta_{n - 1}} \le B \qandq \norm{(\nabla \cL)(u) - (\nabla \cL)(w)} \le L \norm{u - w} \dott
\end{equation}
}
\argument{\lref{mu}; \lref{B_L}}{for all $n \in \N$ that
\begin{equation}\llabel{mu_bound}
\textstyle \norm{\mu_n} \le (\gamma_n)^{- 1} L \norm{2^{- 1} \gamma_n (\nabla \cL)(\Theta_{n - 1})} \le 2^{- 1} L \sup_{k \in \N} \norm{(\nabla \cL)(\Theta_{k - 1})} < \infty \dott
\end{equation}
}
\argument{the assumption that $\sup_{n \in \N} (\norm{\Theta_n} + \sum_{k = -1}^{1} (\gamma_n n^{\lrexpo})^k) < \infty$; \lref{Theta_new}; \lref{mu_bound}; \cref{cor:gen_convergence_rate_Loja_gen_l_rate_Lipschitz_alpha} (applied for every $n \in \N$ with $\fd \with \fd$, $\alpha \with 0$, $\lambda \with \inf_{m \in \N} (\gamma_m m^{\lrexpo})$, $\lrexpo \with \lrexpo$, $\cL \with \cL$, $\Theta_n \with \Theta_n$, $p_n \with 1$, $\gamma_n \with \gamma_n$, $\bbA_n \with \bbI_{\fd}$, $\m_n \with (\nabla \cL)(\Theta_{n - 1}) + \gamma_n \mu_n$, $\mu_n \with \mu_n$ in the notation of \cref{cor:gen_convergence_rate_Loja_gen_l_rate_Lipschitz_alpha})}[verbs=e]{\lref{result}. }
\end{aproof}

\subsubsection*{Acknowledgements}
This work has been partially funded by the European Union (ERC, MONTECARLO, 101\-0\-4\-5\-8\-1\-1). The views and the opinions expressed in this work are however those of the authors only and do not necessarily reflect those of the European Union or the European Research Council (ERC). Neither the European Union nor the granting authority can be held responsible for them. Furthermore, this work has been supported by the Ministry of Culture and Science NRW as part of the Lamarr Fellow Network. In addition, we also gratefully acknowledge the Cluster of Excellence EXC 2044/2-390685587, Mathematics Münster: Dynamics-Geometry-Structure funded by the Deutsche Forschungsgemeinschaft (DFG, German Research Foundation). Most of the specific formulations in the proofs of this work have been created using \cite{Bennoargumentcommand}.

\subsubsection*{Use of large language models}
{\sc Google Gemini} has significantly supported us in creating the literature review in \cref{subsec:literature}. The specific formulations in \cref{subsec:literature} are due to the authors and all formulations and references in \cref{subsec:literature} and the entire article have been carefully verified by the authors.

\newpage

{\small
\bibliography{ref}

\begin{thebibliography}{10}

\bibitem{MR2197994}
{\sc Absil, P.-A., Mahony, R., and Andrews, B.}
\newblock Convergence of the iterates of descent methods for analytic cost
  functions.
\newblock {\em SIAM J. Optim. 16}, 2 (2005), 531--547.

\bibitem{MR3785672}
{\sc Arag\'on~Artacho, F.~J., Fleming, R. M.~T., and Vuong, P.~T.}
\newblock Accelerating the {DC} algorithm for smooth functions.
\newblock {\em Math. Program. 169}, 1 (2018), 95--118.

\bibitem{pmlr-v129-barakat20a}
{\sc Barakat, A., and Bianchi, P.}
\newblock {Convergence Analysis of a Momentum Algorithm with Adaptive Step Size
  for Non Convex Optimization}.
\newblock {\em \href{https://arxiv.org/abs/1911.07596}{arXiv:1911.07596}\/}
  (2020).

\bibitem{Barakat2021}
{\sc Barakat, A., and Bianchi, P.}
\newblock Convergence and dynamical behavior of the {A}dam algorithm for
  nonconvex stochastic optimization.
\newblock {\em SIAM J. Optim. 31}, 1 (2021), 244--274.

\bibitem{MR4298987}
{\sc Barakat, A., Bianchi, P., Hachem, W., and Schechtman, S.}
\newblock {Stochastic optimization with momentum: convergence, fluctuations,
  and traps avoidance}.
\newblock {\em Electron. J. Stat. 15}, 2 (2021), 3892--3947.

\bibitem{MR1082341}
{\sc Benveniste, A., M\'etivier, M., and Priouret, P.}
\newblock {\em {Adaptive algorithms and stochastic approximations}}, vol.~22 of
  {\em Applications of Mathematics (New York)}.
\newblock Springer-Verlag, Berlin, 1990.
\newblock Translated from the French by Stephen S. Wilson.

\bibitem{SignSGD18}
{\sc Bernstein, J., Wang, Y.-X., Azizzadenesheli, K., and Anandkumar, A.}
\newblock {signSGD: Compressed Optimisation for Non-Convex Problems}.
\newblock {\em \href{https://arxiv.org/abs/1802.04434}{arXiv:1802.04434}\/}
  (2018).

\bibitem{PMIHES_1988_67_5_0}
{\sc Bierstone, E., and Milman, P.~D.}
\newblock Semianalytic and subanalytic sets.
\newblock {\em Inst. Hautes \'Etudes Sci. Publ. Math. 67\/} (1988), 5--42.

\bibitem{BockAdamLocalConvergence2019}
{\sc Bock, S., and Weiß, M.}
\newblock {A Proof of Local Convergence for the Adam Optimizer}.
\newblock In {\em 2019 International Joint Conference on Neural Networks
  (IJCNN)\/} (2019), pp.~1--8.

\bibitem{MR2274510}
{\sc Bolte, J., Daniilidis, A., and Lewis, A.}
\newblock {The \L ojasiewicz inequality for nonsmooth subanalytic functions
  with applications to subgradient dynamical systems}.
\newblock {\em SIAM J. Optim. 17}, 4 (2006), 1205--1223.

\bibitem{ChenLiuSunHong2018}
{\sc Chen, X., Liu, S., Sun, R., and Hong, M.}
\newblock {On the Convergence of A Class of Adam-Type Algorithms for Non-Convex
  Optimization}.
\newblock {\em \href{https://arxiv.org/abs/1808.02941}{arXiv:1808.02941}\/}
  (2018).

\bibitem{Crawshaw2022}
{\sc Crawshaw, M., Liu, M., Orabona, F., Zhang, W., and Zhuang, Z.}
\newblock {Robustness to Unbounded Smoothness of Generalized SignSGD}.
\newblock {\em \href{https://arxiv.org/abs/2208.11195}{arXiv:2208.11195}\/}
  (2022).

\bibitem{Defossez2022}
{\sc D{\'e}fossez, A., Bottou, L., Bach, F., and Usunier, N.}
\newblock {A Simple Convergence Proof of Adam and Adagrad}.
\newblock {\em Transactions on Machine Learning Research\/} (2022).

\bibitem{DereichDoJentzen_uniformbound_Adam}
{\sc Dereich, S., Do, T., and Jentzen, A.}
\newblock {Uniform a priori bounds and error analysis for the Adam stochastic
  gradient descent optimization method}.
\newblock {\em \href{https://arxiv.org/abs/2603.18899}{arXiv:2603.18899}\/}
  (2026).

\bibitem{DereichDoJentzenPhilippeAdamSymmetry}
{\sc Dereich, S., Do, T., Jentzen, A., and {von Wurstemberger}, P.}
\newblock {Adam symmetry theorem: characterization of the convergence of the
  stochastic Adam optimizer}.
\newblock {\em \href{https://arxiv.org/abs/2511.06675}{arXiv:2511.06675}\/}
  (2025).

\bibitem{DereichGraeberJentzen_Adam_nonconvergence}
{\sc Dereich, S., Graeber, R., and Jentzen, A.}
\newblock {Non-convergence of Adam and other adaptive stochastic gradient
  descent optimization methods for non-vanishing learning rates}.
\newblock {\em \href{https://arxiv.org/abs/2407.08100}{arXiv:2407.08100}\/}
  (2024).

\bibitem{DereichGraeberJentzenRiekert_asymtoticstability_Adam}
{\sc Dereich, S., Graeber, R., Jentzen, A., and Riekert, A.}
\newblock {Asymptotic stability properties and a priori bounds for Adam and
  other gradient descent optimization methods}.
\newblock {\em \href{https://arxiv.org/abs/2509.10476}{arXiv:2509.10476}\/}
  (2025).

\bibitem{DereichJentzenAdamRates}
{\sc Dereich, S., and Jentzen, A.}
\newblock {Convergence rates for the Adam optimizer}.
\newblock {\em \href{https://arxiv.org/abs/2407.21078}{arXiv:2407.21078}\/}
  (2024).

\bibitem{DereichJentzenKassing_ODE_Adam}
{\sc Dereich, S., Jentzen, A., and Kassing, S.}
\newblock {ODE approximation for the Adam algorithm: General and
  overparametrized setting}.
\newblock {\em \href{https://arxiv.org/abs/2511.04622}{arXiv:2511.04622}\/}
  (2025).

\bibitem{DereichJentzenRiekert_adaptive_lr}
{\sc Dereich, S., Jentzen, A., and Riekert, A.}
\newblock {Learning rate adaptive stochastic gradient descent optimization
  methods: numerical simulations for deep learning methods for partial
  differential equations and convergence analyses}.
\newblock {\em \href{https://arxiv.org/abs/2406.14340}{arXiv:2406.14340}\/}
  (2024).

\bibitem{DereichJentzenRiekert_sharpconvergencerates_Adam}
{\sc Dereich, S., Jentzen, A., and Riekert, A.}
\newblock {Sharp higher order convergence rates for the Adam optimizer}.
\newblock {\em \href{https://arxiv.org/abs/2504.19426}{arXiv:2504.19426}\/}
  (2025).

\bibitem{DereichKassing}
{\sc Dereich, S., and Kassing, S.}
\newblock {Convergence of stochastic gradient descent schemes for
  Lojasiewicz-landscapes}.
\newblock {\em \href{https://arxiv.org/abs/2102.09385}{arXiv:2102.09385}\/}
  (2024).

\bibitem{DozatNadam2}
{\sc Dozat, T.}
\newblock {Incorporating Nesterov Momentum into Adam}.
\newblock \url{https://openreview.net/forum?id=OM0jvwB8jIp57ZJjtNEZ}, 2016.
\newblock [Accessed 14-November-2025].

\bibitem{AdaNorm22}
{\sc Dubey, S.~R., Singh, S.~K., and Chaudhuri, B.~B.}
\newblock {AdaNorm: Adaptive Gradient Norm Correction based Optimizer for
  CNNs}.
\newblock {\em \href{https://arxiv.org/abs/2210.06364}{arXiv:2210.06364}\/}
  (2022).

\bibitem{FawTziotis2022}
{\sc Faw, M., Tziotis, I., Caramanis, C., Mokhtari, A., Shakkottai, S., and
  Ward, R.}
\newblock {The Power of Adaptivity in SGD: Self-Tuning Step Sizes with
  Unbounded Gradients and Affine Variance}.
\newblock {\em \href{https://arxiv.org/abs/2202.05791}{arXiv:2202.05791}\/}
  (2022).

\bibitem{GarrigosGowerSGD}
{\sc Garrigos, G., and Gower, R.~M.}
\newblock {Handbook of Convergence Theorems for (Stochastic) Gradient Methods}.
\newblock {\em \href{https://arxiv.org/abs/2301.11235}{arXiv:2301.11235}\/}
  (2023).

\bibitem{GessKassingMomentum2023}
{\sc Gess, B., and Kassing, S.}
\newblock {Convergence rates for momentum stochastic gradient descent with
  noise of machine learning type}.
\newblock {\em \href{https://arxiv.org/abs/2302.03550}{arXiv:2302.03550}\/}
  (2023).

\bibitem{7330562}
{\sc Ghadimi, E., Feyzmahdavian, H.~R., and Johansson, M.}
\newblock {Global convergence of the Heavy-ball method for convex
  optimization}.
\newblock In {\em 2015 European Control Conference (ECC)\/} (2015),
  pp.~310--315.

\bibitem{GodichonBaggioni2023}
{\sc Godichon-Baggioni, A., and Tarrago, P.}
\newblock {Non asymptotic analysis of Adaptive stochastic gradient algorithms
  and applications}.
\newblock {\em \href{https://arxiv.org/abs/2303.01370}{arXiv:2303.01370}\/}
  (2023).

\bibitem{HeAdamConvergence2023}
{\sc He, M., Liang, Y., Liu, J., and Xu, D.}
\newblock {Convergence of Adam for Non-convex Objectives: Relaxed
  Hyperparameters and Non-ergodic Case}.
\newblock {\em \href{https://arxiv.org/abs/2307.11782}{arXiv:2307.11782}\/}
  (2023).

\bibitem{HintonRMSProp}
{\sc Hinton, G., Srivastava, N., and Swersky, K.}
\newblock {\emph{Lecture 6e: RMSprop: Divide the gradient by a running average
  of its recent magnitude}}.
\newblock
  \url{https://www.cs.toronto.edu/~tijmen/csc321/slides/lecture_slides_lec6.pdf}.
\newblock [Accessed 14-November-2025].

\bibitem{HongLin2024}
{\sc Hong, Y., and Lin, J.}
\newblock {Revisiting Convergence of AdaGrad with Relaxed Assumptions}.
\newblock {\em \href{https://arxiv.org/abs/2402.13794}{arXiv:2402.13794}\/}
  (2024).

\bibitem{JentzenBookDeepLearning2023}
{\sc Jentzen, A., Kuckuck, B., and von Wurstemberger, P.}
\newblock {Mathematical Introduction to Deep Learning: Methods,
  Implementations, and Theory}.
\newblock {\em \href{https://arxiv.org/abs/2310.20360}{arXiv:2310.20360}\/}
  (2023).

\bibitem{MR4716376}
{\sc Jentzen, A., and Riekert, A.}
\newblock On the existence of global minima and convergence analyses for
  gradient descent methods in the training of deep neural networks.
\newblock {\em J. Mach. Learn. 1}, 2 (2022), 141--246.

\bibitem{JENTZEN2020101438}
{\sc Jentzen, A., and von Wurstemberger, P.}
\newblock Lower error bounds for the stochastic gradient descent optimization
  algorithm: sharp convergence rates for slowly and fast decaying learning
  rates.
\newblock {\em J. Complexity 57\/} (2020), 101438, 16.

\bibitem{Jiang2023UAdam}
{\sc Jiang, Y., Liu, J., Xu, D., and Mandic, D.~P.}
\newblock {UAdam: Unified Adam-Type Algorithmic Framework for Non-Convex
  Stochastic Optimization}.
\newblock {\em \href{https://arxiv.org/abs/2305.05675}{arXiv:2305.05675}\/}
  (2023).

\bibitem{MR4663539}
{\sc Josz, C., Lai, L., and Li, X.}
\newblock Convergence of the momentum method for semialgebraic functions with
  locally {L}ipschitz gradients.
\newblock {\em SIAM J. Optim. 33}, 4 (2023), 3012--3037.

\bibitem{KingmaBaAdam}
{\sc Kingma, D.~P., and Ba, J.}
\newblock {Adam: A Method for Stochastic Optimization}.
\newblock {\em
  \href{https://arxiv.org/abs/1412.6980}{arXiv:\allowbreak1412.6980}\/} (2014).

\bibitem{Bennoargumentcommand}
{\sc Kuckuck, B.}
\newblock {Some useful LATEX commands}.
\newblock {\em
  \href{https://latex.bennokuckuck.de}{https://latex.bennokuckuck.de} [Accessed
  April-2026]\/} (2025).

\bibitem{Li2023}
{\sc Li, H., Rakhlin, A., and Jadbabaie, A.}
\newblock {Convergence of Adam Under Relaxed Assumptions}.
\newblock {\em \href{https://arxiv.org/abs/2304.13972}{arXiv:2304.13972}\/}
  (2023).

\bibitem{LIU2022300}
{\sc Liu, J., Kong, J., Xu, D., Qi, M., and Lu, Y.}
\newblock {Convergence analysis of AdaBound with relaxed bound functions for
  non-convex optimization}.
\newblock {\em Neural Networks 145\/} (2022), 300--307.

\bibitem{pmlr-v178-liu22d}
{\sc Liu, J., and Yuan, Y.}
\newblock {On Almost Sure Convergence Rates of Stochastic Gradient Methods}.
\newblock {\em \href{https://arxiv.org/abs/2202.04295}{arXiv:2202.04295}\/}
  (2022).

\bibitem{LiuJiang2019}
{\sc Liu, L., Jiang, H., He, P., Chen, W., Liu, X., Gao, J., and Han, J.}
\newblock {On the Variance of the Adaptive Learning Rate and Beyond}.
\newblock {\em \href{https://arxiv.org/abs/1908.03265}{arXiv:1908.03265}\/}
  (2019).

\bibitem{NEURIPS2020_d3f5d4de}
{\sc Liu, Y., Gao, Y., and Yin, W.}
\newblock {An Improved Analysis of Stochastic Gradient Descent with Momentum}.
\newblock {\em \href{https://arxiv.org/abs/2007.07989}{arXiv:2007.07989}\/}
  (2020).

\bibitem{LoizouSHB2017}
{\sc Loizou, N., and Richt\'{a}rik, P.}
\newblock {Linearly convergent stochastic heavy ball method for minimizing
  generalization error}.
\newblock {\em \href{https://arxiv.org/abs/1710.10737}{arXiv:1710.10737}\/}
  (2017).

\bibitem{MR4174637}
{\sc Loizou, N., and Richt\'{a}rik, P.}
\newblock Momentum and stochastic momentum for stochastic gradient, {N}ewton,
  proximal point and subspace descent methods.
\newblock {\em Comput. Optim. Appl. 77}, 3 (2020), 653--710.

\bibitem{Lojasiewicz}
{\sc {\L}ojasiewicz, S.}
\newblock {\emph{Ensembles semi-analytiques}. Unpublished lecture notes.
  Institut des Hautes Études Scientifiques. 1965}.
\newblock {\sc url:}
  \url{https://perso.univ-rennes1.fr/michel.coste/Lojasiewicz.pdf}.

\bibitem{LuoAdaBoundAMSBound2019}
{\sc Luo, L., Xiong, Y., Liu, Y., and Sun, X.}
\newblock {Adaptive Gradient Methods with Dynamic Bound of Learning Rate}.
\newblock {\em \href{https://arxiv.org/abs/1902.09843}{arXiv:1902.09843}\/}
  (2019).

\bibitem{Nesterov1983AMF}
{\sc Nesterov, Y.}
\newblock {A method of solving a convex programming problem with convergence
  rate $O(1/k^2)$}.
\newblock {\em Soviet Mathematics Doklady 27\/} (1983), 372--376.

\bibitem{POLYAK19641}
{\sc Polyak, B.}
\newblock Some methods of speeding up the convergence of iteration methods.
\newblock {\em USSR Computational Mathematics and Mathematical Physics 4}, 5
  (1964), 1--17.

\bibitem{Polyak_accelerating}
{\sc Polyak, B.~T., and Juditsky, A.~B.}
\newblock Acceleration of stochastic approximation by averaging.
\newblock {\em SIAM J. Control Optim. 30}, 4 (1992), 838--855.

\bibitem{QiuMaMilzarek_momentum_KL_24}
{\sc Qiu, J., Ma, B., and Milzarek, A.}
\newblock {Convergence of SGD with momentum in the nonconvex case: A time
  window-based analysis}.
\newblock {\em \href{https://arxiv.org/abs/2405.16954}{arXiv:2405.16954}\/}
  (2024).

\bibitem{ReddiKale2019}
{\sc Reddi, S.~J., Kale, S., and Kumar, S.}
\newblock {On the Convergence of Adam and Beyond}.
\newblock {\em \href{https://arxiv.org/abs/1904.09237}{arXiv:1904.09237}\/}
  (2019).

\bibitem{SRuderGDOverview}
{\sc Ruder, S.}
\newblock {An overview of gradient descent optimization algorithms}.
\newblock {\em \href{https://arxiv.org/abs/1609.04747}{arXiv:1609.04747}\/}
  (2016).

\bibitem{pmlr-v134-sebbouh21a}
{\sc Sebbouh, O., Gower, R.~M., and Defazio, A.}
\newblock {Almost sure convergence rates for Stochastic Gradient Descent and
  Stochastic Heavy Ball}.
\newblock In {\em Proceedings of Thirty Fourth Conference on Learning Theory\/}
  (15--19 Aug 2021), M.~Belkin and S.~Kpotufe, Eds., vol.~134 of {\em
  Proceedings of Machine Learning Research}, PMLR, pp.~3935--3971.

\bibitem{SunHBconvergence2018}
{\sc Sun, T., Yin, P., Li, D., Huang, C., Guan, L., and Jiang, H.}
\newblock {Non-ergodic Convergence Analysis of Heavy-Ball Algorithms}.
\newblock {\em \href{https://arxiv.org/abs/1811.01777}{arXiv:1811.01777}\/}
  (2018).

\bibitem{MR3315611}
{\sc Tadi\'c, V.~B.}
\newblock {Convergence and convergence rate of stochastic gradient search in
  the case of multiple and non-isolated extrema}.
\newblock {\em Stochastic Process. Appl. 125}, 5 (2015), 1715--1755.

\bibitem{WangZhang2023}
{\sc Wang, B., Zhang, H., Ma, Z., and Chen, W.}
\newblock {Convergence of AdaGrad for Non-convex Objectives: Simple Proofs and
  Relaxed Assumptions}.
\newblock {\em \href{https://arxiv.org/abs/2305.18471}{arXiv:2305.18471}\/}
  (2023).

\bibitem{WangSAdam2020}
{\sc Wang, G., Lu, S., Tu, W., and Zhang, L.}
\newblock {SAdam: A Variant of Adam for Strongly Convex Functions}.
\newblock {\em \href{https://arxiv.org/abs/1905.02957}{arXiv:1905.02957}\/}
  (2019).

\bibitem{WuAdaLoss2021}
{\sc Wu, X., Xie, Y., Du, S., and Ward, R.}
\newblock {AdaLoss: A computationally-efficient and provably convergent
  adaptive gradient method}.
\newblock {\em \href{https://arxiv.org/abs/2109.08282}{arXiv:2109.08282}\/}
  (2021).

\bibitem{XiaAdamL2023}
{\sc Xia, L., and Massei, S.}
\newblock {AdamL: A fast adaptive gradient method incorporating loss function}.
\newblock {\em \href{https://arxiv.org/abs/2312.15295}{arXiv:2312.15295}\/}
  (2023).

\bibitem{MR4723898}
{\sc Xiao, N., Hu, X., Liu, X., and Toh, K.-C.}
\newblock Adam-family methods for nonsmooth optimization with convergence
  guarantees.
\newblock {\em J. Mach. Learn. Res. 25\/} (2024), Paper No. [48], 53.

\bibitem{Adan}
{\sc Xie, X., Zhou, P., Li, H., Lin, Z., and Yan, S.}
\newblock {Adan: Adaptive Nesterov Momentum Algorithm for Faster Optimizing
  Deep Models}.
\newblock {\em \href{https://arxiv.org/abs/2208.06677}{arXiv:2208.06677}\/}
  (2022).

\bibitem{YangMomentum2016}
{\sc Yang, T., Lin, Q., and Li, Z.}
\newblock {Unified Convergence Analysis of Stochastic Momentum Methods for
  Convex and Non-convex Optimization}.
\newblock {\em \href{https://arxiv.org/abs/1604.03257}{arXiv:1604.03257}\/}
  (2016).

\bibitem{YuanEAdam2020}
{\sc Yuan, W., and Gao, K.-X.}
\newblock {EAdam Optimizer: How $\varepsilon$ Impact Adam}.
\newblock {\em \href{https://arxiv.org/abs/2011.02150}{arXiv:2011.02150}\/}
  (2020).

\bibitem{YuanNonAdaptiveOPT2022}
{\sc Yuan, W., Hu, F., and Lu, L.}
\newblock {A new non-adaptive optimization method: Stochastic gradient descent
  with momentum and difference}.
\newblock {\em Applied Intelligence 52}, 4 (2022), 3939--3953.

\bibitem{Yogi_NEURIPS2018}
{\sc Zaheer, M., Reddi, S., Sachan, D., Kale, S., and Kumar, S.}
\newblock Adaptive methods for nonconvex optimization.
\newblock In {\em Advances in Neural Information Processing Systems\/} (2018),
  S.~Bengio, H.~Wallach, H.~Larochelle, K.~Grauman, N.~Cesa-Bianchi, and
  R.~Garnett, Eds., vol.~31, Curran Associates, Inc.

\bibitem{ZhangFang2020}
{\sc Zhang, B., Jin, J., Fang, C., and Wang, L.}
\newblock {Improved Analysis of Clipping Algorithms for Non-convex
  Optimization}.
\newblock {\em \href{https://arxiv.org/abs/2010.02519}{arXiv:2010.02519}\/}
  (2020).

\bibitem{ZhangChen2022}
{\sc Zhang, Y., Chen, C., Shi, N., Sun, R., and Luo, Z.-Q.}
\newblock {Adam Can Converge Without Any Modification On Update Rules}.
\newblock {\em \href{https://arxiv.org/abs/2208.09632}{arXiv:2208.09632}\/}
  (2022).

\bibitem{MR4637063}
{\sc Zhou, Y., Huang, K., Cheng, C., Wang, X., Hussain, A., and Liu, X.}
\newblock Fast{A}da{B}elief: improving convergence rate for belief-based
  adaptive optimizers by exploiting strong convexity.
\newblock {\em IEEE Trans. Neural Netw. Learn. Syst. 34}, 9 (2023), 6515--6529.

\bibitem{ZhuangAdaBelief2020}
{\sc Zhuang, J., Tang, T., Ding, Y., Tatikonda, S., Dvornek, N., Papademetris,
  X., and Duncan, J.~S.}
\newblock {AdaBelief Optimizer: Adapting Stepsizes by the Belief in Observed
  Gradients}.
\newblock {\em \href{https://arxiv.org/abs/2010.07468}{arXiv:2010.07468}\/}
  (2020).

\bibitem{ZouShen2019}
{\sc Zou, F., Shen, L., Jie, Z., Zhang, W., and Liu, W.}
\newblock {A Sufficient Condition for Convergences of Adam and RMSProp}.
\newblock {\em \href{https://arxiv.org/abs/1811.09358}{arXiv:1811.09358}\/}
  (2018).

\end{thebibliography}
\bibliographystyle{acm}}

\end{document}